\documentclass[12pt,a4paper]{article}
\usepackage[T1]{fontenc}
\usepackage{amssymb}
\usepackage{amsmath}
\usepackage{stmaryrd}
\usepackage{mathrsfs}
\usepackage{mathtools}
\usepackage{bm}
\usepackage{graphicx}
\usepackage{here}
\usepackage[margin=1in]{geometry}
\usepackage{cite}
\usepackage{sectsty}
\allsectionsfont{\sffamily}
\DeclareMathAlphabet{\mathscrbf}{OMS}{mdugm}{b}{n}

\usepackage{amsbsy}
\usepackage{enumitem}

\usepackage{tikz}
\usepackage[leftcaption]{sidecap}
\usepackage{caption}

\let\OLDthebibliography\thebibliography
\renewcommand\thebibliography[1]{
  \OLDthebibliography{#1}
  \setlength{\parskip}{0pt}
  \setlength{\itemsep}{.5pt plus 0.0ex}
}

\usepackage{xcolor}
\definecolor{blueMath}{rgb}{0.368417,0.506779,0.709798}
\definecolor{yellowMath}{rgb}{0.880722,0.611041,0.142051}
\definecolor{greenMath}{rgb}{0.560181,0.691569,0.194885}
\definecolor{redMath}{rgb}{0.922526,0.185626,0.109179}

\usepackage[colorlinks,allcolors=blue]{hyperref}
\newcommand{\be}{\begin{equation}}
\newcommand{\ee}{\end{equation}}
\newcommand{\ben}{\begin{enumerate}}
\newcommand{\een}{\end{enumerate}}
\newcommand{\bi}{\begin{itemize}}
\newcommand{\ei}{\end{itemize}}
\newcommand{\bmm}{\begin{pmatrix}}
\newcommand{\emm}{\end{pmatrix}}

\newcommand{\Ad}{\text{Ad}}

\newcommand{\arccot}{\,\text{arccot}}

\newcommand{\cn}{\,\text{cn}}
\newcommand{\cotanh}{\,\text{cotanh}}
\newcommand{\dd}{\text{d}}
\newcommand{\demi}{\frac{1}{2}}
\newcommand{\der}{\partial}
\newcommand{\Diff}{\text{Diff}\,S^1}
\newcommand{\dn}{\,\text{dn}}

\newcommand{\eg}{e.g.\ }

\newcommand{\ie}{i.e.\ }

\newcommand{\nn}{\nonumber}

\newcommand{\SL}{\text{SL}(2,\mathbb{R})}

\newcommand{\sn}{\,\text{sn}}

\newcommand{\Vect}{\text{Vect}\,S^1}

\newcommand{\bg}{\textbf g}

\newcommand{\bu}{\textbf u}

\newcommand{\cE}{{\cal E}}
\newcommand{\cF}{{\cal F}}

\newcommand{\cO}{{\cal O}}

\newcommand{\sfC}{\mathsf{C}}

\newcommand{\sfM}{\mathsf{M}}
\newcommand{\sfS}{\mathsf{S}}

\newcommand{\CC}{\mathbb{C}}

\newcommand{\II}{\mathbb{I}}

\newcommand{\RR}{\mathbb{R}}
\newcommand{\ZZ}{\mathbb{Z}}

\newcommand\Tstrut{\rule{0pt}{2.6ex}}
\newcommand\Bstrut{\rule[-1.3ex]{0pt}{0pt}}

\begin{document}

\hrule
\begin{center}
\Large{\bfseries{\textsf{Orbital Bifurcations and Shoaling of Cnoidal Waves}}}
\end{center}
\hrule
~\\

\begin{center}
\large{\textsf{Blagoje Oblak}}
~\\
~\\
\begin{minipage}{.9\textwidth}\small\it
\begin{center}
{\tt boblak@phys.ethz.ch}\\
Institut f\"ur Theoretische Physik\\
ETH Z\"urich\\
CH-8093 Z\"urich, Switzerland
\end{center}
\end{minipage}
\end{center}

\vspace{1cm}

\begin{center}
\begin{minipage}{.9\textwidth}
\begin{center}{\bfseries{\textsf{Abstract}}}\end{center}
We study the parameter space of cnoidal waves --- the periodic solitons of the Korteweg-de Vries equation --- from the perspective of Virasoro coadjoint orbits. The monodromy method familiar from inverse scattering implies that many, but not all, of these solitons are conformally equivalent to uniform field configurations (constant coadjoint vectors). The profiles that have no uniform representative lie in Lam\'e band gaps and are separated from the others by bifurcation lines along which the corresponding orbits change from elliptic to hyperbolic. We show that such bifurcations can be produced by shoaling: wave profiles become non-uniformizable once their pointedness parameter crosses a certain critical value (which we compute). As a by-product, we also derive asymptotic relations between the pointedness and velocity of cnoidal waves along orbital level curves.
\end{minipage}
\end{center}

\thispagestyle{empty}
\newpage

\textsf{\tableofcontents}
~\\

\newpage

\section{Introduction}

The Korteweg-de Vries (KdV) equation \cite{Korteweg} is an essential tool in physics thanks to its wide range of applications \cite{Crighton} --- most notably in fluid dynamics \cite{Dejager,Ockendon} --- and its powerful integrability properties \cite{Gardner,Zabusky,Airault}. Cnoidal waves are prominent solutions of KdV: they are periodic solitons \cite{Novikov,NovikovBook} and are stable against perturbations \cite{Benjamin}. In this work, we study these waves from the point of view of Virasoro coadjoint orbits \cite{Lazutkin,Balog:1997zz,Guieu,Oblak:2016eij}. Indeed, the relation between the KdV equation and the Virasoro algebra \cite{Khesin} of two-dimensional conformal field theories (CFTs) allows one to think of KdV wave profiles as CFT stress tensors. It is then natural to ask how cnoidal waves are affected by conformal transformations, \ie circle diffeomorphisms. This question is partly motivated by the rich symplectic structure of coadjoint orbits of Lie groups \cite{KirillovLectures}; sharp bifurcations occurring in this structure, as one moves on the space of cnoidal waves, suggest the existence of observable qualitative transitions in certain properties of the wave-carrying medium. In fact, this paper stems from an attempt to describe Stokes drift in shallow water \cite{Longuet53} as a group-theoretic Berry phase \cite{OblakKozy}, where it is essential to know the Virasoro orbits of cnoidal waves: different orbits generally lead to noticeably different Berry phases. Transitions between orbits then occur as a result of shoaling \cite{Svendsen,Dingemans,Dean}, where waves propagate in a fluid with slowly varying depth.\\

Abstractly, the description of orbits of cnoidal waves rests on well established mathematical tools that feature prominently in the inverse scattering method \cite{Gelfand} for solving the KdV equation \cite{Gardner,Novikov,NovikovBook}. The first is the classification of Virasoro coadjoint orbits in terms of invariant conjugacy classes of monodromy matrices for solutions of Hill's equation \cite{Magnus,Lazutkin,Balog:1997zz}. The second is the fact that, for cnoidal waves, Hill's equation reduces to the simplest Lam\'e equation, whose exact solutions are known \cite{Ince,Ince2,Whittaker,Arscott}. This allows us to address, analytically, the following question:
\be
\text{{\it When is a cnoidal wave conformally equivalent to a uniform field configuration?}}
\nn
\ee
When this is possible, one can think of the constant profile as the `rest-frame' counterpart of the wave. Thanks to the Lam\'e equation, we obtain the closed-form expression (\ref{ss159}) for this uniform orbit representative (when it exists). This result is summarized by fig.\ \ref{s71b}: a bifurcation diagram in the two-dimensional parameter space of cnoidal waves. The picture exhibits three distinct regions separated by two bifurcation lines, corresponding to three qualitatively different families of Virasoro orbits.\footnote{Technically, fig.\ \ref{s71b} displays {\it four} regions and {\it three} bifurcation lines; but as we explain at the end of section \ref{sec32}, the upper line $V=\tfrac{2{-}m}{3}$ is not, in fact, a bifurcation.} In particular, the interior of the bifurcation wedge (the green region in fig.\ \ref{s71b}) consists of cnoidal waves whose orbit has {\it no} uniform representative at all: such waves have no rest frame. By contrast, all waves outside of that `forbidden wedge' do have a rest frame and can be mapped on a constant. Some physical implications of this fact are discussed in sections \ref{sec4} and \ref{sec5}. In particular, we will show that shoaling generally allows wave profiles to enter the bifurcation wedge when their `pointedness' parameter crosses the critical value $m^*\simeq0.8261$.\\

For the record, the mathematical tools we shall rely on have been known for a long time, and our observations mostly follow from basic aspects of elliptic functions. In particular, the bifurcation diagram of fig.\ \ref{s71b} is closely related to the well known band structure of the $N=1$ Lam\'e equation \cite{Ince,Ince2,Whittaker,Arscott}, since the `forbidden wedge' is the standard band gap between `valence' and `conduction' bands. This is also intimately related to the spectrum of the cnoidal Lax operator used in the inverse scattering method \cite{Novikov,NovikovBook}. Our goal is not to rederive these results, but rather to work out their consequences for Virasoro orbits and exhibit sharp transitions caused by otherwise smooth parameter variations, with a view towards fluid mechanics. In addition, we provide asymptotic formulas for orbit representatives of cnoidal waves in various limiting regions of their parameter space. Our hope is that these results pave the way for later applications such as the upcoming work \cite{OblakKozy}. We have also attempted to make the text self-contained and pedagogical, so that no prior familiarity with elliptic functions is assumed.\\

The paper is organized as follows. In section \ref{sec2} we introduce the KdV equation and briefly recall the classification of Virasoro coadjoint orbits. Section \ref{sec3} then describes the Virasoro orbits of cnoidal profiles, found thanks to the known solutions of the Lam\'e equation, leading to the bifurcation diagram of fig.\ \ref{s71b}. The asymptotic behaviour of orbital level curves in various corners of that diagram is studied in section \ref{sec36}, while the motion of profiles due to shoaling is described in section \ref{sec4}. We conclude and list some potential follow-ups in section \ref{sec5}. For completeness, appendix \ref{appA} reviews the Virasoro group; appendix \ref{appB} similarly reviews Weierstrass and Jacobi elliptic functions. Appendix \ref{appC} contains technical computations relevant for the asymptotic analysis of section \ref{sec36}.

\section{KdV and Virasoro coadjoint orbits}
\label{sec2}

In this section we set the stage by introducing the KdV equation and recalling how it is related to Virasoro coadjoint orbits. We also briefly review the classification of these orbits, which will be essential to find the orbits of cnoidal waves.

\subsection{The KdV equation}
\label{sec21}

The {\it KdV equation} describes the dynamics of a field $p(x,\tau)$ in $(1+1)$ space-time dimensions. We adopt the following normalization for this equation:
\be
\dot p
+3pp'
-\frac{c}{12}\,p'''
=
0,
\label{s5}
\ee
where $\dot p\equiv\der p/\der\tau$, $p'\equiv\der p/\der x$ and the Virasoro central charge $c\in\RR$ is a free parameter.\footnote{See appendix \ref{appA} for more on the Virasoro group. When $c=0$, eq.\ (\ref{s5}) is the inviscid Burgers equation.} In fluid mechanics, $p(x,\tau)$ describes a wave profile in a one-dimensional shallow water channel, with $x$ a comoving spatial coordinate and $\tau$ a slow time variable. (The quantities $p$, $c$, $x$ and $\tau$ are all dimensionless.) In that context the central charge is typically negative, but in order to be as close as possible to the literature on the Virasoro group, we will stick to the normalization (\ref{s5}) with arbitrary $c$.\\

In principle, the KdV equation (\ref{s5}) holds on an entire spatial real line with a coordinate $x\in\RR$, but we shall focus on solutions that have a fixed spatial period. Accordingly, we let KdV dynamics take place on a unit {\it circle}, whereupon any profile is $2\pi$-periodic:
\be
p(x+2\pi,\tau)
=
p(x,\tau)
\qquad\forall\,x,\tau\in\RR.
\label{s6}
\ee
This corresponds to focussing on spatially periodic profiles with wavelength $2\pi$. It entails no loss of generality within the description of periodic waves in fluid dynamics, where the derivation of KdV assumes the wavelength to be fixed anyway \cite{Ockendon} (see section \ref{sec4}).

\paragraph{Cnoidal waves.} With the above conventions, a {\it cnoidal wave} is a soliton solution of (\ref{s5}),
\be
p(x,\tau)
=
\frac{v}{3}
-\frac{cK(m)^2}{9\pi^2}(m+1)
+\frac{cK(m)^2}{3\pi^2}\,m\sn^2\biggl(\frac{K(m)}{\pi}(x-v\tau)\bigg|m\biggr),
\label{ss6}
\ee
where $\sn(\cdot|m)$ is the Jacobi elliptic sine with squared modulus $m$ and $K(m)$ is the complete elliptic integral of the first kind.\footnote{See section \ref{secB2} in appendix \ref{appB} for the definition of $\sn(x|m)$ and $K(m)$.} The other parameters specifying the solution (\ref{ss6}) are the dimensionless velocity $v$ and the central charge $c$.\footnote{Given eq.\ (\ref{ss6}), it is natural to call $v$ a `velocity'. But in fluid mechanics, the actual velocity of a wave (\ref{ss6}), as seen from a static laboratory frame, is $\propto(1+\epsilon v)$ with $\epsilon\ll1$. More on this in section \ref{sec4}.\label{FootNoteLab}} Thus, at fixed central charge and fixed wavelength, cnoidal waves are labelled by two parameters: the pointedness $m\in[0,1)$ and the velocity $v\in\RR$. Note that the precise combination of $(m,v)$-dependent coefficients in (\ref{ss6}) does not come out of the blue, as it is dictated by KdV dynamics: not every function of the form $A+B\sn^2[K(m)(x-v\tau)/\pi|m]$ solves eq.\ (\ref{s5}). Also note that, in fluid dynamics, the central charge $c$ is negative, so the profile (\ref{ss6}) is a wave with flat troughs and sharp crests, as is indeed observed in experiments.\\

For later convenience we trade the velocity $v$ for a rescaled version thereof,
\be
V
\equiv
\frac{2\pi^2\,v}{c\,K(m)^2}.
\label{t154}
\ee
The justification of this particular combination will become apparent in section \ref{sec31}. For the time being, just note that in these terms a cnoidal wave (\ref{ss6}) reads
\be
p(x,\tau)
=
\frac{cK(m)^2}{3\pi^2}
\Biggl[
\frac{V}{2}
-\frac{m+1}{3}
+m\sn^2\Biggl(\frac{K(m)}{\pi}\biggl(x-\frac{cK(m)^2}{2\pi^2}V\tau\biggr)\Bigg|m\Biggr)
\Biggr].
\label{CNOMEGA}
\ee
From now on we will label cnoidal waves by parameters $(m,V)$ instead of $(m,v)$. In fluid dynamics \cite{Ockendon}, cnoidal waves are often further restricted by requiring that their integral over one wavelength vanish; this turns $V$ into a function of $m$ (see eq.\ (\ref{ss31})), so the space of parameters loses one dimension. We refrain from imposing this restriction at the moment, but we will eventually use it in section \ref{sec4} to deal with shoaling.

\paragraph{KdV and Virasoro.} The KdV equation (\ref{s5}) is related to Virasoro symmetry as follows: if one thinks of $p(x,\tau)$ as the (left-moving component of the) stress tensor of a Lorentzian CFT in $(1+1)$ dimensions, then eq.\ (\ref{s5}) says that $p$ generates its own conformal transformations:
\be
\dot p
=
-pp'-2p'p+\frac{c}{12}p'''.
\label{s7}
\ee
Readers familiar with CFT may recognize in (\ref{s7}) an infinitesimal conformal transformation of $p$ generated by the vector field $\xi=p$. It is the coadjoint representation of the Virasoro algebra, in accordance with the fact that KdV is an Euler-Arnold equation for the Virasoro group \cite{Khesin}. This implies that any solution of (\ref{s5}) is confined to a single Virasoro coadjoint orbit. In other words, there always exists a time-independent profile $k(x)$ and a path $f_{\tau}(x)$ in the group $\Diff$ of diffeomorphisms of the circle such that
\be
p(x,\tau)
=
\bigl(f_{\tau}\cdot k\bigr)(x)
\label{s8}
\ee
where the dot denotes the coadjoint representation of the Virasoro group at central charge $c$. Explicitly, this means that
\be
p\bigl(f_{\tau}(x),\tau\bigr)
=
\frac{1}{\bigl(f_{\tau}'(x)\bigr)^2}
\Bigl[
k(x)
+
\frac{c}{12}\sfS[f_{\tau}](x)
\Bigr]
\equiv
\bigl(f_{\tau}\cdot k\bigr)(f_{\tau}(x))
\label{ss8}
\ee
where $\sfS[f]$ is the {\it Schwarzian derivative} of $f$,
\be
\sfS[f](x)
\equiv
\frac{f'''(x)}{f'(x)}-\frac{3}{2}\biggl(\frac{f''(x)}{f'(x)}\biggr)^2.
\label{shabadoo}
\ee
Indeed one can verify that (\ref{ss8}) solves (\ref{s5}) when $f_{\tau}(x)$ is the time-ordered exponential of the vector field $p(x,\tau)\partial_x$. This is related to Euler-Poincar\'e reconstruction \cite{Holm}, which we shall study in much greater detail in the upcoming publication \cite{OblakKozy}.\\

Note that in (\ref{s8}) the time-independent profile $k(x)$ need not (though it may) be the initial configuration $p(x,0)$. In particular, the simplest case occurs when $k(x)=k$ is independent of $x$, \ie {\it uniform} (or {\it constant}). Then eq.\ (\ref{s8}) says that $p(x,\tau)$ is, at any time $\tau$, a conformal transform of a fixed uniform field configuration. An appropriate analogy here is that the map $f_{\tau}$ sending $p(x,\tau)$ on the constant $k$ is a change of reference frames such that $k$ is the `rest-frame value' of $p(x,\tau)$. In general, it is not guaranteed that such a rest frame exists at all: some Virasoro orbits have no constant representative. The main goal of this paper is to understand whether cnoidal waves (\ref{CNOMEGA}) belong to orbits of constant profiles. To answer this, we now review the coadjoint orbits of the Virasoro group.

\subsection{Virasoro orbits and Hill's equation}
\label{sec22}

We have just seen that any solution (with wavelength $2\pi$) of the KdV equation (\ref{s5}) lies in a single Virasoro coadjoint orbit
\be
\cO_k
=
\big\{
f\cdot k
\big|
f\in\Diff
\big\}
\label{s10}
\ee
where the dot is the coadjoint representation (\ref{ss8}) and $k(x)$ is some time-independent profile. The question then is whether these orbits can be fully classified, and whether one can find an exhaustive, non-redundant set of `orbit representatives' such that each representative $k$ defines one orbit (\ref{s10}), and all orbits are accounted for. The complete answer to these questions \cite{Lazutkin} goes beyond the scope of this paper; see \cite{Balog:1997zz,Guieu,Khesin,Oblak:2016eij} for detailed reviews. In the remainder of this section we merely summarize the aspects of this classification that will be essential for orbits of cnoidal waves, and that also play a key role for inverse scattering through the `monodromy method' \cite{NovikovBook}.\\

The starting point is to find invariants that label different Virasoro orbits. To do this, pick a profile $p(x)$ and a central charge $c$; the corresponding {\it Hill's equation} \cite{Magnus,NovikovBook} is
\be
-\frac{c}{6}\psi''(x)+p(x)\psi(x)=0,
\label{s11}
\ee
where $\psi(x)$ is a function on $\RR$ that need not be periodic. If $\psi_1(x)$ and $\psi_2(x)$ are two linearly independent solutions of this equation, the periodicity (\ref{s6}) of $p(x)$ implies the existence of a {\it monodromy matrix} $\sfM\in\SL$ such that
\be
\bmm
\psi_1(x+2\pi) \\ \psi_2(x+2\pi) 
\emm
=
\sfM
\bmm
\psi_1(x) \\ \psi_2(x)
\emm
\qquad\forall\,x\in\RR.
\label{ss11}
\ee
One can then show that the conjugacy class of $\sfM$ is Virasoro-invariant; in other words, applying the coadjoint transformation (\ref{ss8}) to $p(x)$ while also transforming $\psi_1,\psi_2$ as densities\footnote{See the definition (\ref{sa4}) in appendix \ref{appA}.} with weight $-1/2$ changes the matrix $\sfM$ by at most conjugation in $\SL$. The same is true when changing the basis of solutions. As a result, the terminology of $\SL$ conjugacy classes carries over to Virasoro orbits: an orbit is
\begin{itemize}
\setlength\itemsep{0em}
\item {\it elliptic} if its monodromies are elliptic ($|\text{Tr}\,\sfM|<2$);
\item {\it hyperbolic} if its monodromies are hyperbolic ($|\text{Tr}\,\sfM|>2$);
\item {\it parabolic} if monodromies are parabolic ($|\text{Tr}\,\sfM|=2$ but non-zero off-diagonal terms);
\item {\it exceptional} if its monodromies are trivial ($\sfM=\pm\II$).
\end{itemize}
Note that the map which sends a Virasoro orbit on the conjugacy class of the corresponding monodromy matrices is not injective: different orbits may have identical monodromies. We will describe a less rough classification in section \ref{sec23}. As for the relation between this notion of monodromy and the Bloch ansatz of condensed matter physics, see section \ref{sec34}.\\

The conjugacy class of Hill's monodromy is thus an invariant parameter labelling Virasoro orbits. In particular, when $p(x)$ belongs to the orbit of a uniform profile $k$, the trace of $\sfM$ takes the same value both at $k$ and at $p(x)$. Now, when Hill's equation (\ref{s11}) has a constant potential $k$, its general solution is a sum of exponentials $\exp[\pm\sqrt{6k/c}\,x]$. If $k/c<0$, it is understood that the exponentials are oscillating, while if $k=0$ the solution collapses to $Cx+D$. In any event, the corresponding monodromy matrix has trace
\be
\text{Tr}\,\sfM
=
\begin{cases}
2\cosh\bigl(2\pi\sqrt{6k/c}\bigr) & \text{if }\;\,k/c\geq0,\\
2\cos\bigl(2\pi\sqrt{6|k/c|}\bigr) & \text{if }\;\,k/c\leq0.
\end{cases}
\label{s12}
\ee
Thus, if one solves Hill's equation with a profile $p(x)$ that admits a rest frame and computes its monodromy $\sfM$, the constant $k$ such that $p(x)=(f\cdot k)(x)$ satisfies (\ref{s12}). Of course, at this stage we do not yet know if $p(x)$ admits a rest frame to begin with; eq.\ (\ref{s12}) implies at least that a necessary condition for this frame to exist is $\text{Tr}\,\sfM\geq-2$.

\subsection{Classification and map of Virasoro orbits}
\label{sec23}

Conjugacy classes of monodromy matrices only provide a rough classification of Virasoro orbits. Different conjugacy classes certainly label different orbits, but the converse is untrue: different orbits may have identical monodromies. Thus, to obtain a more accurate classification, one must include a {\it second} invariant label of Virasoro orbits, namely the {\it winding number} $n$ of the ratio $\psi_2/\psi_1$. By definition, $n$ is the number of laps around a circle performed by the function $\psi_2(x)/\psi_1(x)$ as $x$ runs from $0$ to $2\pi$, when seen as a projective (stereographic) coordinate on the circle. One can show that this winding is invariant under $\Diff$: it takes the same value for any two profiles $p(x)$, $q(x)$ such that $p=f\cdot q$ for some $f\in\Diff$. Furthermore, along with monodromy matrices, it furnishes the complete classification of Virasoro orbits: two orbits are different if and only if they have different windings {\it and} their monodromies have different conjugacy classes \cite{Khesin}.\\

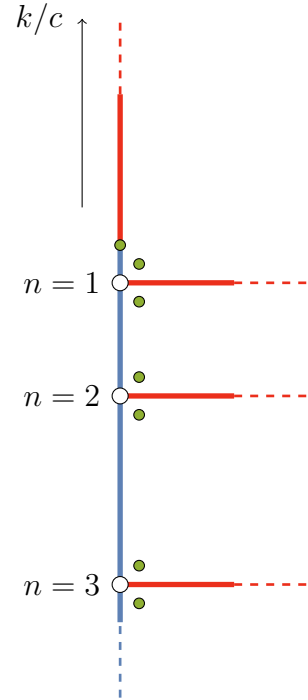
\begin{SCfigure}[2][t]
\centering
\begin{tikzpicture}
\draw [->] (-0.5,0.5) -- (-0.5,3) node[left] {$k/c\;$};
\draw[blueMath,dashed,line width=1pt] (0,-6)--(0,-5);
\draw[blueMath,line width=2pt] (0,-5)--(0,0);
\draw[redMath,line width=2pt] (0,0)--(0,2);
\draw[redMath,dashed,line width=1pt] (0,2)--(0,3);
\draw[redMath,line width=2pt] (0,-1/2)--(1.5,-1/2);
\draw[redMath,dashed,line width=1pt] (1.5,-1/2)--(2.5,-1/2);
\draw[redMath,line width=2pt] (0,-4/2)--(1.5,-4/2);
\draw[redMath,dashed,line width=1pt] (1.5,-4/2)--(2.5,-4/2);
\draw[redMath,line width=2pt] (0,-9/2)--(1.5,-9/2);
\draw[redMath,dashed,line width=1pt] (1.5,-9/2)--(2.5,-9/2);
\filldraw[fill=greenMath,draw=black] (0,0) circle (2pt);
\filldraw[fill=greenMath,draw=black] (1/4,-1/2+1/4) circle (2pt);
\filldraw[fill=greenMath,draw=black] (1/4,-1/2-1/4) circle (2pt);
\filldraw[fill=greenMath,draw=black] (1/4,-4/2+1/4) circle (2pt);
\filldraw[fill=greenMath,draw=black] (1/4,-4/2-1/4) circle (2pt);
\filldraw[fill=greenMath,draw=black] (1/4,-9/2+1/4) circle (2pt);
\filldraw[fill=greenMath,draw=black] (1/4,-9/2-1/4) circle (2pt);
\filldraw[fill=white,draw=black] (0,-1/2) circle (3pt) node[left] {$n=1\;$};
\filldraw[fill=white,draw=black] (0,-4/2) circle (3pt) node[left] {$n=2\;$};
\filldraw[fill=white,draw=black] (0,-9/2) circle (3pt) node[left] {$n=3\;$};
\end{tikzpicture}
\caption{A map of Virasoro orbits, adapted from \cite{Oblak:2016eij,Balog:1997zz}. Elliptic orbits are blue; hyperbolic ones are red; parabolic ones are green dots; exceptional ones are white dots (here visible at $n=1,2,3$). The value of $k/c$ increases along the vertical axis. In fig.\ \ref{s62} below we will see where cnoidal waves fit in this map.\label{s19}}
\end{SCfigure}

We now summarize the classification of orbits that stems from this seminal result \cite{Lazutkin} (for many more details, see \cite{Balog:1997zz,Guieu,Khesin,Oblak:2016eij}). First, not all orbits have constant representatives: at non-zero winding $n$, both hyperbolic and parabolic monodromies label orbits having no rest frame, respectively corresponding to `tachyonic' or `massless' Virasoro orbits.\footnote{This terminology is borrowed from \cite{Barnich:2014kra}, where Virasoro orbits are analogues of Poincar\'e mass shells.} In both cases, $n$ determines the sign of the trace of the monodromy matrix according to
\be
\text{Tr}\,\sfM
=
(-1)^n|\text{Tr}\,\sfM|,
\label{s17}
\ee
which will be important in section \ref{sec32}. By contrast, the other combinations of invariants all correspond to orbits that {\it do} have a rest frame:
\begin{itemize}
\setlength\itemsep{0em}
\item Winding $n\in\{0,1,2,...\}$ with elliptic monodromy labels the orbit of the constant
\be
k
=
-\frac{c}{6}\biggl(
\;\left\lfloor\frac{n+1}{2}\right\rfloor
+\frac{(-1)^n}{2\pi}\arccos(\text{Tr}\,\sfM/2)
\biggr)^2,
\label{t18}
\ee
which is indeed consistent with (\ref{s12}). We choose the branch of the $\arccos$ function such that $\arccos(-1)=\pi$ and $\arccos(1)=0$.
\item Winding $n=0$ with hyperbolic monodromy labels the orbit of the constant $k$ given by (\ref{s12}) with $\text{Tr}\,\sfM>2$:
\be
k
=
\frac{c}{6}
\biggl[
\frac{1}{2\pi}
\text{arcosh}\bigl(\text{Tr}\,\sfM/2\bigr)
\biggr]^2.
\label{s18}
\ee
\item Winding $n=0$ with parabolic monodromy labels the orbit of the constant $k=0$.
\item Winding $n\in\{1,2,3,...\}$ with $\sfM=(-1)^n\,\II$ labels the `exceptional orbit' of
\be
k=-\frac{n^2c}{24}.
\label{ss18}
\ee
(We shall return to the epithet `exceptional' shortly.)
\end{itemize}
This classification is summarized by the `map of orbits' of fig.\ \ref{s19}, in which every point corresponds to exactly one orbit:
\begin{itemize}
\item The points located on the vertical axis are orbits that admit a rest frame, with a constant representative given by eqs.\ (\ref{t18})-(\ref{ss18}). The value of $k/c$ increases as one moves upwards. On that axis, the points labelled by an integer $n$ are exceptional orbits with representatives (\ref{ss18}). All other points below $k=0$ are elliptic, with winding $n=\lfloor\sqrt{24|k/c|}\rfloor$. Above $k=0$, orbits are hyperbolic with winding $n=0$. 
\item All points that do {\it not} lie on the vertical axis in fig.\ \ref{s19} are orbits that have no rest frame, \ie {\it no} uniform representative. These points form continuous families of hyperbolic orbits (horizontal lines) and discrete pairs of parabolic orbits, with a winding number $n$ indicated on the left.
\end{itemize}

Importantly, fig.\ \ref{s19} is roughly consistent with the topology of the set of Virasoro orbits: points that can be continuously connected to one another in the map are indeed related by continuous deformations of orbit representatives.\footnote{The set of Virasoro orbits is an orbifold, so words such as `topology' and `continuity' must be handled with care. What we mean here is that functions on the circle that are close to one another (\eg in terms of the supremum norm) have Virasoro orbits whose points in fig.\ \ref{s19} are close to each other.} For example, the orbit of a constant $k+\varepsilon$ is at distance $\sim\varepsilon$ from that of the constant $k$ on the vertical line of the map. Similarly, hyperbolic orbits with non-zero winding $n$ (the $n^{\text{th}}$ horizontal line in fig.\ \ref{s19}) can be obtained as deformations of the constant (\ref{ss18}), which is why this family of orbits is connected to the vertical axis \cite{Balog:1997zz}. Note how exceptional orbits lie at the intersection of hyperbolic and elliptic orbits. This is closely related to the epithet `exceptional', which actually stresses the fact that constants of the form (\ref{ss18}) have an enhanced stabilizer. Indeed, the subgroup of $\Diff$ elements such that $f\cdot k=k$ (with $f\cdot k$ given by (\ref{ss8})) is three-dimensional only when the profile $k(x)$ belongs to the orbit of an exceptional constant (\ref{ss18}). For any other orbit, the stabilizer is only one-dimensional.\\

We now investigate where cnoidal waves fit in fig.\ \ref{s19}. As we shall see, some waves belong to the $n=1$ horizontal hyperbolic orbit, hence have no rest frame. They will turn out to be separated from other cnoidal waves (which do have a rest frame) by bifurcation lines where the uniform representative $k$ takes the exceptional value (\ref{ss18}) with $n=1$.

\newpage
\section{Orbits of cnoidal waves}
\label{sec3}

As reviewed in section \ref{sec21}, cnoidal waves (with wavelength $2\pi$) are travelling-wave solutions $p(x,\tau)=p(x-v\tau)$ of the KdV equation (\ref{s5}), with a time-independent shape described by the profile (\ref{CNOMEGA}) in terms of the rescaled velocity (\ref{t154}):
\be
p(x)
=
\frac{cK(m)^2}{3\pi^2}
\biggl[
\frac{V}{2}-\frac{m+1}{3}
+m\sn^2\biggl(\frac{K(m)}{\pi}x\bigg|m\biggr)
\biggr].
\label{s21}
\ee
In what follows we investigate the Virasoro orbits spanned by the conformal transforms of these functions as $m\in[0,1)$ and $V\in\RR$ vary, following the monodromy method outlined in section \ref{sec22}. The results are summarized in fig.\ \ref{s71b} and should look familiar to readers acquainted with the Lam\'e equation \cite{Whittaker,Arscott} or inverse scattering \cite{NovikovBook}. The link between Virasoro orbits and the Lam\'e band structure is established in section \ref{sec34}. The required background on elliptic functions can be found in appendix \ref{appB}.

\subsection{Lam\'e equation and cnoidal monodromies}
\label{sec31}

As explained in section \ref{sec22}, to find the Virasoro orbit of a profile $p(x)$, we first study its Hill's equation (\ref{s11}). Letting $y\equiv K(m)x/\pi$, for a cnoidal wave (\ref{s21}) this yields a {\it Lam\'e equation} with label $N=1$ \cite{Whittaker,Arscott}:
\be
-\frac{\dd^2\psi}{\dd y^2}(y)
+\biggl[
V-\frac{2m+2}{3}+2m\sn^2(y|m)
\biggr]
\psi(y)
=
0.
\label{s153}
\ee
The solutions of that equation are known explicitly; for comparison with the literature it is easiest to convert (\ref{s153}) to its Weierstrassian form, where the profile $\propto\sn^2$ is expressed in terms of a Weierstrass elliptic function $\wp$. Namely, one has\footnote{See the proof around eq.\ (\ref{b56}) in appendix \ref{appB}.}
\be
m\sn^2(y|m)
=
\wp\Bigl(y+iK(1{-}m),K(m),iK(1{-}m)\Bigr)
+\frac{m+1}{3},
\label{s57}
\ee
where $K(m)$ and $iK(1{-}m)$ are the half-periods of the $\wp$ function. Thus, introducing the complex variable $z\equiv y+iK(1{-}m)$ and writing $\psi(y)\equiv\phi(z)$, eq.\ (\ref{s153}) reads
\be
-\frac{\dd^2\phi}{\dd z^2}
+2\wp\bigl(
z,K(m),iK(1{-}m)
\bigr)\phi(z)
=
-V\,\phi(z).
\label{ss154}
\ee
This is the Weierstrassian form of the $N=1$ Lam\'e equation with eigenvalue $-V$. Solving it amounts to solving Hill's equation for cnoidal waves, which in turn will allow us to find their Virasoro orbits. By the way, we can now see why the redefinition (\ref{t154}) is useful: had we written everything in terms of the original velocity $v$, the eigenvalue on the right-hand side of (\ref{ss154}) would have been $\propto v/K(m)^2$ and would also have depended on the central charge $c$. For simplicity, from now on we write $\wp\bigl(\cdot,K(m),iK(1{-}m)\bigr)\equiv\wp(\cdot)$, neglecting to stress the half-periods.\\

\begin{figure}[t]
\centering
\includegraphics[width=0.49\textwidth]{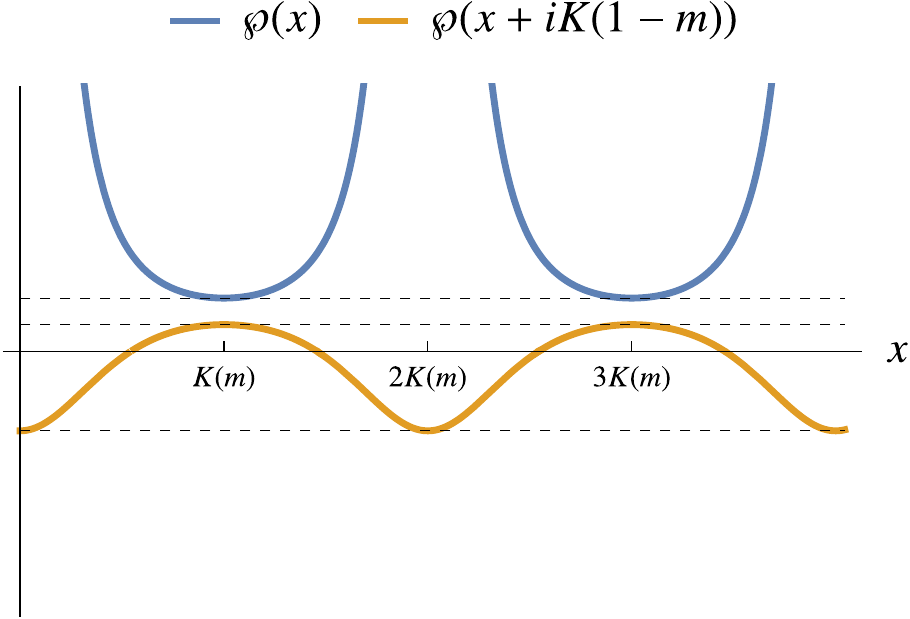}
\hfill
\includegraphics[width=0.49\textwidth]{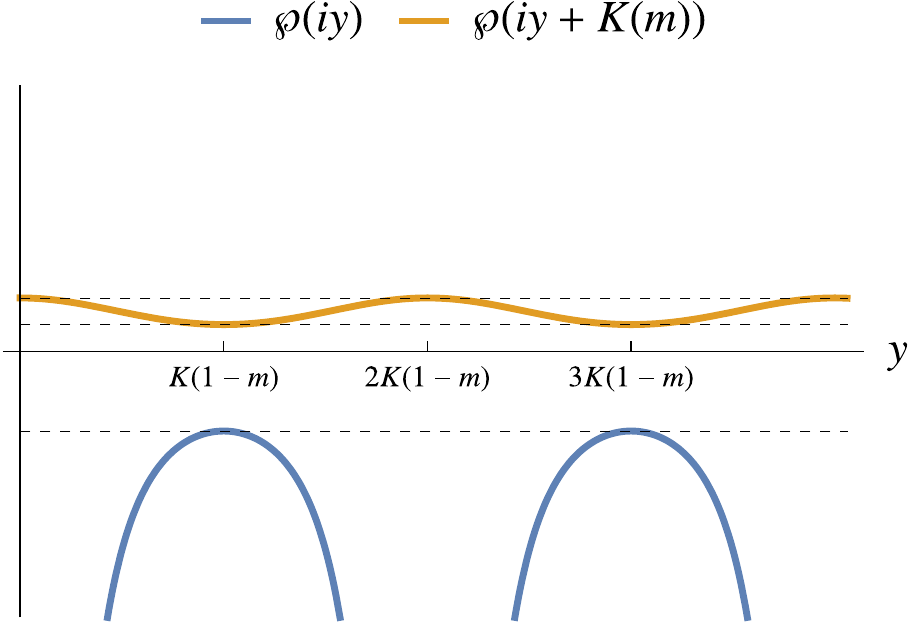}
\caption{The Weierstrass function $\wp(z)$ with $m=0.8$ along four lines in the complex $z$ plane where it takes real values. As functions of $x,y\in\RR$, $\wp(x)$ and $\wp(x+iK(1{-}m))$ have period $2K(m)$, while $\wp(iy)$ and $\wp(iy+K(m))$ have period $2K(1{-}m)$. The horizontal dashed lines indicate the extrema of the $\wp$ function: $\wp(iK(1{-}m))={-}\tfrac{m+1}{3}$, $\wp(K(m)+iK(1{-}m))=\tfrac{2m{-}1}{3}$ and $\wp(K(m))=\tfrac{2{-}m}{3}$. The plots make it manifest that the inverse function $\wp^{-1}$ in eq.\ (\ref{s156}) is multi-valued.\label{s157}}
\end{figure}

The solutions of (\ref{ss154}) can be written in closed form in terms of the Weierstrass zeta and sigma functions, $\zeta(z)$ and $\sigma(z)$, such that\footnote{See the more detailed definitions (\ref{b44}) and (\ref{bbb47}) in appendix \ref{appB}.} $\zeta'=-\wp$ and $\sigma'=\zeta\sigma$. Indeed, in these terms, two linearly independent solutions of the Lam\'e equation (\ref{ss154}) are
\be
\phi_{\pm}(z)
=
\frac{\pm\sigma\bigl(z\pm\wp^{-1}(V)\bigr)}{\sigma(z)\sigma\bigl(\wp^{-1}(V)\bigr)}
\,e^{\mp\zeta(\wp^{-1}(V))\,z},
\label{s156}
\ee
where $\wp^{-1}(V)\equiv a\in\CC$ is such that $\wp(a)=V$. To prove that (\ref{s156}) solves eq.\ (\ref{ss154}), one uses $\zeta'=-\wp$ and $\sigma'=\zeta\sigma$ to compute the second derivative of (\ref{s156}):
\be
\phi_{\pm}''(z)
=
\bigl(
-\wp(z\pm a)+\wp(z)
\bigr)\phi_{\pm(z)}
+\bigl(
\zeta(z\pm a)\mp\zeta(a)-\zeta(z)
\bigr)^2\phi_{\pm}(z),
\label{b53}
\ee
where $a\equiv\wp^{-1}(V)$. The addition formulas\footnote{These identities are derived at the end of appendix \ref{appBsec1}.}
\be
\wp(z+w)
=
-\wp(z)-\wp(w)+\biggl(\frac{\wp'(z)-\wp'(w)}{2\wp(z)-2\wp(w)}\biggr)^2,
\quad
\zeta(z+w)
=
\zeta(z)+\zeta(w)+\frac{\wp'(z)-\wp'(w)}{2\wp(z)-2\wp(w)}
\label{bb51}
\ee
then turn (\ref{b53}) into the Lam\'e equation (\ref{ss154}), as desired.\\

Note that the solutions (\ref{s156}) are (almost) uniquely specified by the rescaled velocity $V$ even though the inverse Weierstrass function $\wp^{-1}$ is multi-valued (see fig.\ \ref{s157}). Indeed, $\wp(z)$ is doubly periodic, so the value of $\wp^{-1}(V)$ is only fixed up to the addition of factors $2nK(m)+2in'K(1{-}m)$ with arbitrary integers $n,n'$. This ambiguity does not affect the solution (\ref{s156}) owing to the quasi-periodicity of $\zeta(z)$ and $\sigma(z)$.\footnote{See eqs.\ (\ref{bbb44}) and (\ref{bb47}) in appendix \ref{appB}.} In particular,
\be
\zeta(z+2K(m))
=
\zeta(z)+2K(m),
\qquad
\sigma(z+2K(m))
=
-e^{2\zeta(K(m))(z+K(m))}\sigma(z).
\label{t156}
\ee
Since the $\wp$ function is even, there is a further ambiguity in the sign of $\wp^{-1}(V)$; this exchanges the solutions $\phi_+$ and $\phi_-$ of eq.\ (\ref{s156}), but it does not affect the trace of the monodromy matrix (computed in eq.\ (\ref{s159}) below). This is enough to ensure that the solutions (\ref{s156}) uniquely determine a Virasoro orbit representative as a function of $(m,V)$.\\

\begin{figure}[t]
\centering
\includegraphics[width=0.49\textwidth]{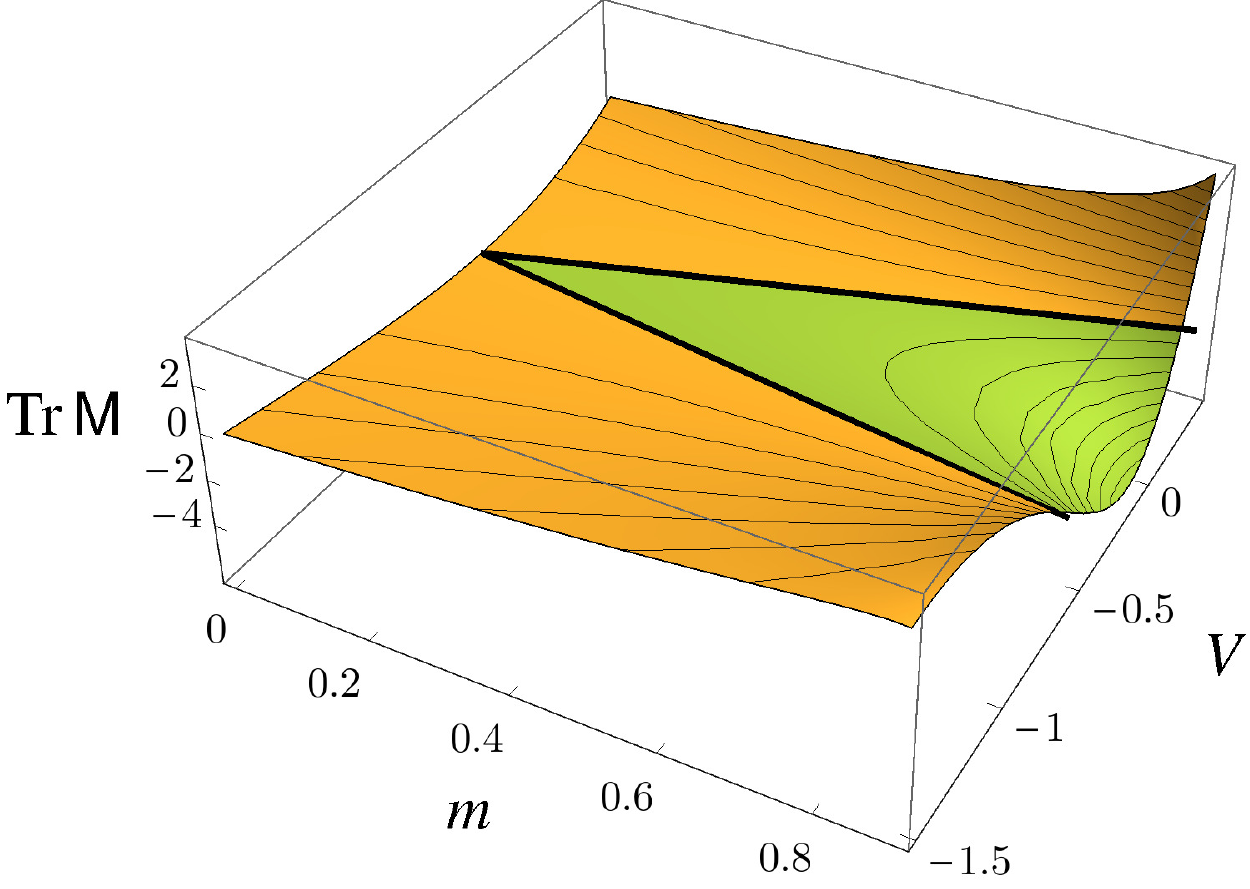}
\hfill
\includegraphics[width=0.49\textwidth]{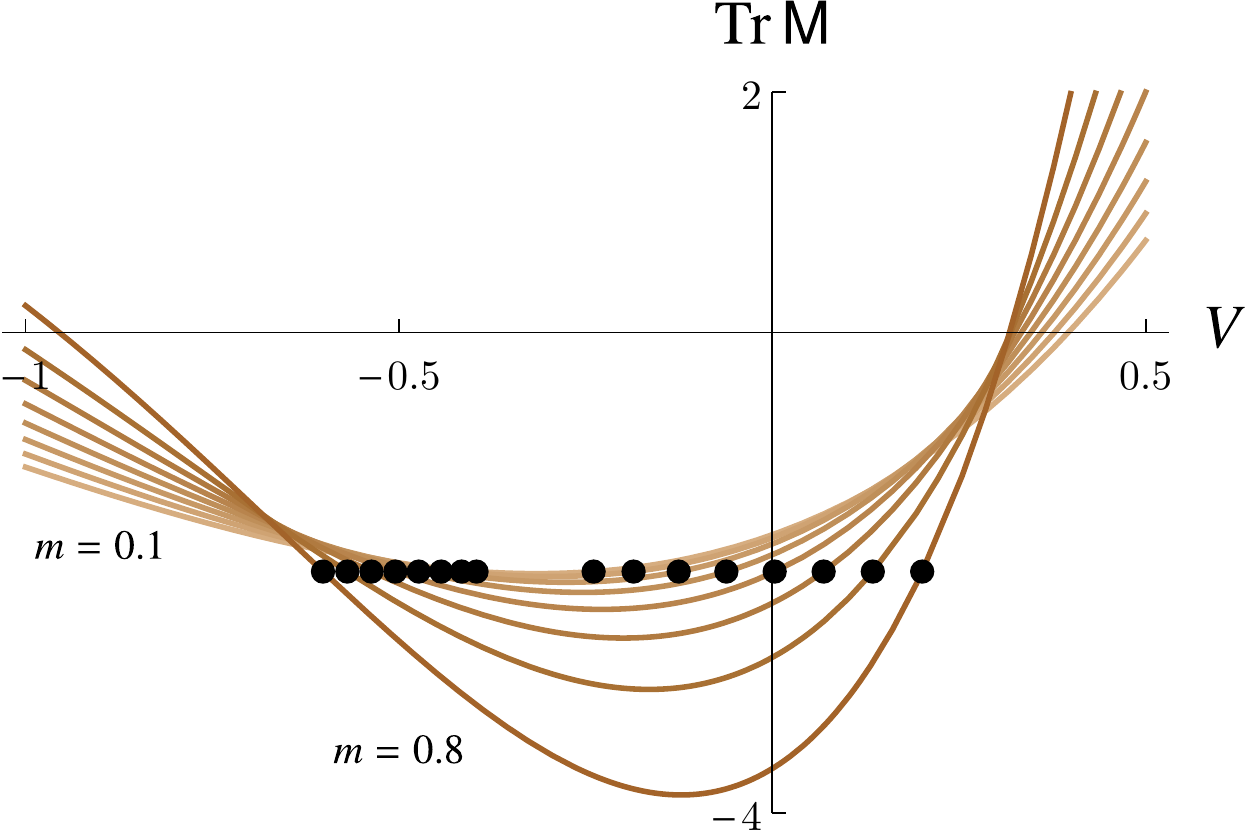}
\caption{The trace (\ref{s159}) of the monodromy of the solutions of the Lam\'e equation (\ref{ss154}) as a function of the parameters $m$ (pointedness) and $V$ (velocity and height) specifying a cnoidal wave (\ref{s21}). Note the wedge in $(m,V)$ space where $\text{Tr}\,\sfM<-2$, implying the lack of rest frame for profiles located in that region. {\it Left panel:} The wedge is highlighted in green and the black curves are level sets of $\text{Tr}\,\sfM$. We let $m$ run from $0$ to $0.9$ to avoid the problematic region $m\to1$ where $\text{Tr}\,\sfM$ blows up. {\it Right panel:} The trace is plotted as a function of $V$ at fixed $m$ ranging from $m=0.1$ to $m=0.8$; as $m$ increases, so does the size of the region where $\text{Tr}\,\sfM<-2$ (in fact it does so linearly with $m$).
\label{s59b}}
\end{figure}

To evaluate the monodromy matrix $\sfM$ of eq.\ (\ref{ss11}), recall that above (\ref{s153}) we rescaled $x\sim x+2\pi$ into $y\sim y+2K(m)$, so the monodromy that affects the solutions (\ref{s156}) occurs when $z=y+iK(1{-}m)$ increases by $2K(m)$. Using the periodicities (\ref{t156}) and the oddness of $\zeta$, we find a diagonal monodromy with trace
\be
\text{Tr}\,\sfM
=
2\cosh
\Bigl[
2K(m)\zeta\bigl(\wp^{-1}(V)\bigr)
-2\zeta\bigl(K(m)\bigr)\,\wp^{-1}(V)
\Bigr].
\label{s159}
\ee
This is real for all values of $m\in[0,1)$ and $V\in\RR$, despite the generally complex argument of the $\cosh$ function; in particular, $\text{Tr}\,\sfM$ need not be larger than $2$. We stress once more that the multi-valuedness of $\wp^{-1}$ does not entail any ambiguity in this expression. The trace is plotted in fig.\ \ref{s59b} as a function of $m,V$. In what follows we interpret this plot in terms of Virasoro orbits.

\subsection{Cnoidal orbits and bifurcations}
\label{sec32}

\begin{figure}[t]
\centering
\includegraphics[width=0.45\textwidth]{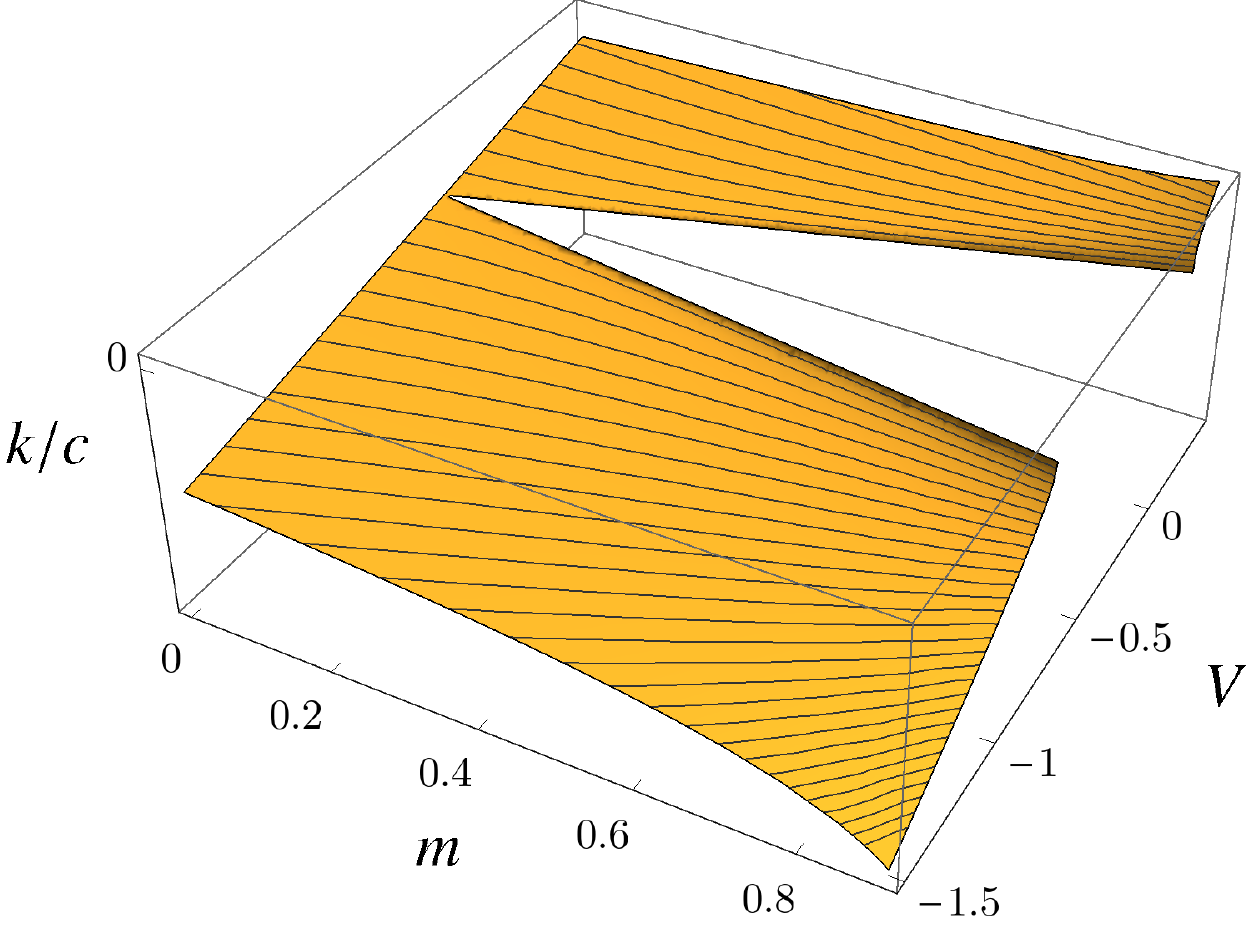}
\hfill
\includegraphics[width=0.45\textwidth]{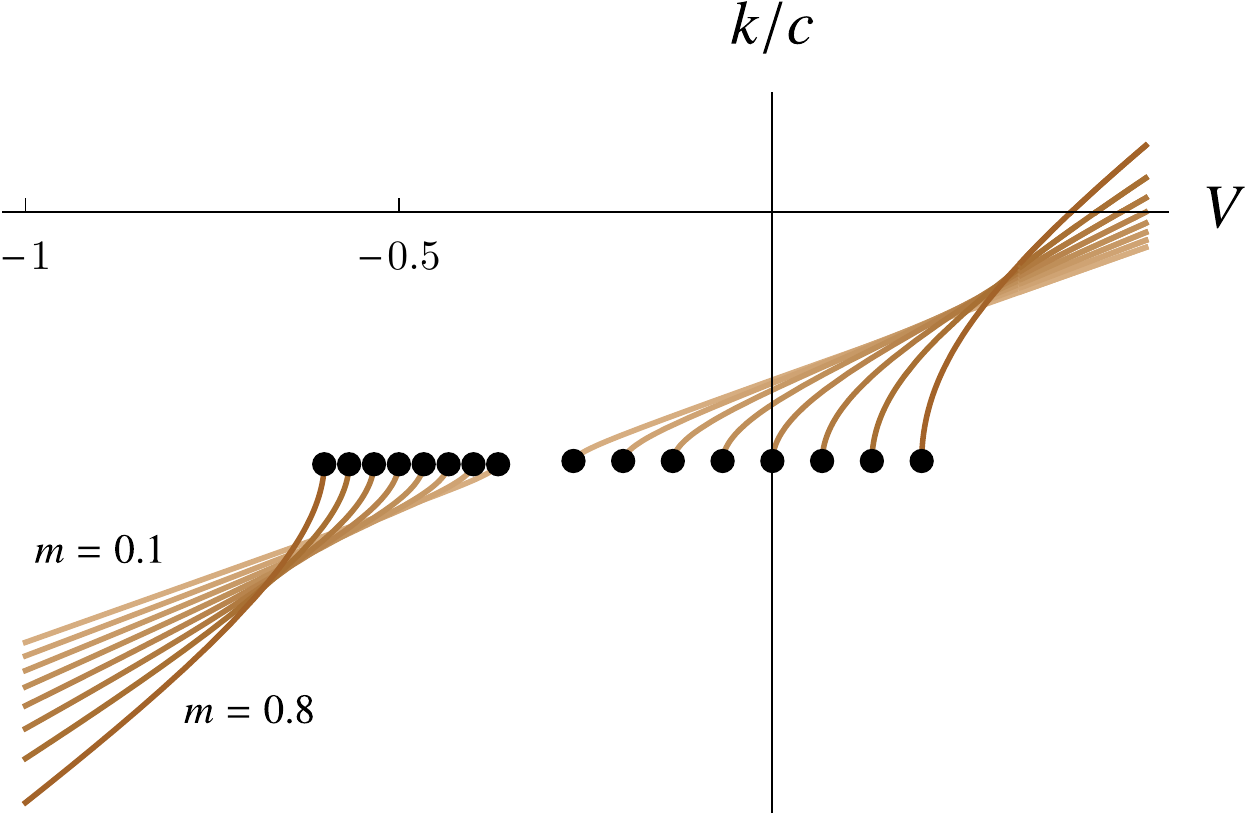}
\caption{The uniform orbit representative (\ref{ss159}) of a cnoidal profile (\ref{s21}) as a function of $m,V$. {\it Left panel:} The parameter $m$ ranges from $0$ to $0.9$ to avoid the singularity at $m\to1$. The black curves are level sets of $k$. Note the wedge where $k$ is not real; on its boundaries, $k=-c/24$. {\it Right panel:} The value of $k/c$ is plotted as a function of $V$ at fixed $m$ ranging from $m=0.1$ to $m=0.8$. Black dots indicate points where $k=-c/24$.\label{s60b}}
\end{figure}

The trace (\ref{s159}) yields the rest-frame representative $k$ of any cnoidal wave: using eq.\ (\ref{s18}),
\be
\boxed{
k
=
\frac{c}{6\pi^2}
\Bigl[
K(m)\zeta\bigl(\wp^{-1}(V)\bigr)
-
\zeta\bigl(K(m)\bigr)\wp^{-1}(V)
\Bigr]^2.
}
\label{ss159}
\ee
This is, in its most condensed form, the result we wanted to obtain: for any cnoidal profile that admits a rest frame, eq.\ (\ref{ss159}) gives the corresponding uniform configuration in terms of the wave's parameters $(m,V)$. It also tells us which profiles do admit a rest frame, and which ones do not: a wave with parameters $(m,V)$ has a rest frame if and only if the imaginary part of (\ref{ss159}) vanishes --- which occurs in the region where $\text{Tr}\,\sfM\geq-2$. By contrast, (\ref{ss159}) acquires an imaginary part when $\text{Tr}\,\sfM<-2$; the corresponding cnoidal waves have no rest frame. This is summarized in fig.\ \ref{s60b}.\\

Eq.\ (\ref{ss159}) is a well-known result in the literature on the Lam\'e equation, the KdV equation and inverse scattering \cite{Whittaker,Arscott,NovikovBook}, though to our knowledge it has never been applied to Virasoro orbits. In particular, the wedge where $\text{Tr}\,\sfM<-2$ in figs.\ \ref{s59b} and \ref{s60b} is really the Lam\'e band gap (see section \ref{sec34}). Our goal now is to classify the Virasoro orbits specified by (\ref{ss159}) as $m$ and $V$ vary. As we shall see, the result is as follows:
\begin{itemize}
\setlength\itemsep{0em}
\item For $V<{-}\tfrac{m+1}{3}$, the function (\ref{ss159}) is real and satisfies $k/c<-1/24$. The corresponding orbits have a rest frame and their winding number is $n\geq1$. They span the vertical line below the point $n=1$ in fig.\ \ref{s19}.
\item The lines $V=[(1\pm3)m-2]/6$ are the boundaries of the wedge in which $\text{Tr}\,\sfM\leq-2$. Along these lines, $k=-c/24$ (see eq.\ (\ref{KOMEGGA})). The corresponding Virasoro orbit is the exceptional one with winding $n=1$. It is the point $n=1$ in fig.\ \ref{s19}.
\item In the wedge where ${-}\tfrac{m+1}{3}<V<\tfrac{2m{-}1}{3}$, (\ref{ss159}) is complex. The corresponding Virasoro orbits have no rest frame. Since continuous orbit variations are accompanied by jumps of at most $\pm1$ in winding number, these orbits are in fact hyperbolic with winding $n=1$. They span the $n=1$ horizontal line in fig.\ \ref{s19}.
\item Above the wedge, when $V>\tfrac{2m{-}1}{3}$, the function (\ref{ss159}) is real and satisfies $k/c>-1/24$, with $k/c$ increasing monotonically with $V$, and $k=0$ when $V=\tfrac{2{-}m}{3}$. The corresponding orbits all have a rest frame and winding $n=0$, and span the vertical line above the point $n=1$ in fig.\ \ref{s19}. They are elliptic for $\tfrac{2m{-}1}{3}<V<\tfrac{2{-}m}{3}$, parabolic when $V=\tfrac{2{-}m}{3}$ and hyperbolic when $V>\tfrac{2{-}m}{3}$, but there is no bifurcation at $V=\tfrac{2{-}m}{3}$ in the sense that the function (\ref{ss159}) is regular there.
\end{itemize}
These statements are summarized in figs.\ \ref{s62} and \ref{s71b}. The remainder of this section is devoted to their proof, which relies on a detailed (and standard \cite{Whittaker,Arscott}) analysis of the functions $\wp^{-1}$ and $\zeta$ that appear in eq.\ (\ref{ss159}). Readers familiar with elliptic functions or the Lam\'e band structure may safely skip the proof and go straight to section \ref{sec36}.

\begin{SCfigure}[2][t]
\centering
\begin{tikzpicture}
\draw [->] (-0.5,0.5) -- (-0.5,3) node[left] {$k/c\;$};
\draw[redMath,dashed,line width=1pt] (0,-6)--(0,-5);
\draw[redMath,line width=2pt] (0,-5)--(0,0);
\draw[redMath,line width=2pt] (0,0)--(0,2);
\draw[redMath,dashed,line width=1pt] (0,2)--(0,3);
\draw[redMath,line width=2pt] (0,-1/2)--(1.5,-1/2);
\draw[redMath,dashed,line width=1pt] (1.5,-1/2)--(2.5,-1/2);
\draw[black,line width=2pt] (0,-4/2)--(1.5,-4/2);
\draw[black,dashed,line width=1pt] (1.5,-4/2)--(2.5,-4/2);
\draw[black,line width=2pt] (0,-9/2)--(1.5,-9/2);
\draw[black,dashed,line width=1pt] (1.5,-9/2)--(2.5,-9/2);
\filldraw[fill=redMath,draw=black] (0,0) circle (2pt);
\filldraw[fill=black,draw=black] (1/4,-1/2+1/4) circle (2pt);
\filldraw[fill=black,draw=black] (1/4,-1/2-1/4) circle (2pt);
\filldraw[fill=black,draw=black] (1/4,-4/2+1/4) circle (2pt);
\filldraw[fill=black,draw=black] (1/4,-4/2-1/4) circle (2pt);
\filldraw[fill=black,draw=black] (1/4,-9/2+1/4) circle (2pt);
\filldraw[fill=black,draw=black] (1/4,-9/2-1/4) circle (2pt);
\filldraw[fill=redMath,draw=black] (0,-1/2) circle (3pt) node[left] {$n=1\;$};
\filldraw[fill=redMath,draw=black] (0,-4/2) circle (3pt) node[left] {$n=2\;$};
\filldraw[fill=redMath,draw=black] (0,-9/2) circle (3pt) node[left] {$n=3\;$};
\end{tikzpicture}
\caption{The map of Virasoro orbits of fig.\ \ref{s19}, now stressing the fact that the orbits highlighted in red contain cnoidal waves. The vertical line above (respectively below) $n=1$ is spanned by the orbits located above (respectively below) the wedge in fig.\ \ref{s60b}. The horizontal line at $n=1$ is spanned by the orbits in the wedge of fig.\ \ref{s60b}. The wedge's boundaries consist of cnoidal waves that all belong to the exceptional $n=1$ orbit, here represented at the intersection of the vertical axis and the $n=1$ horizontal line.\label{s62}}
\end{SCfigure}

\begin{figure}[t]
\centering
\begin{tikzpicture}
\draw[->,line width=1.5pt] (0,0)--(0,1.5*0.89);
\draw[line width=1.5pt] (0,1.5*0.89)--(0,1.5*1.66);
\draw[->,line width=1.5pt] (0,1.5*1.66)--(1.5*1.129,1.5*1.66);
\draw[line width=1.5pt] (1.5*1.129,1.5*1.66)--(1.5*2.25,1.5*1.66);
\draw[->,line width=1.5pt] (1.5*2.25,1.5*1.66)--(1.5*2.25,1.5*0.79);
\draw[line width=1.5pt] (1.5*2.25,1.5*0.79)--(1.5*2.25,0);
\draw[->,line width=1.5pt] (1.5*2.25,0)--(1.5*1.125,0);
\draw[line width=1.5pt] (1.5*1.125,0)--(0,0);
\draw[->,line width=.5pt] (0,-.5)--(0,4) node[left] {Im};
\draw[->,line width=.5pt] (-.5,0)--(4.5,0) node[above] {Re$\;\;$};
\filldraw[fill=white,draw=black,thick] (0,0) circle (3pt) node[below,yshift=-.1cm,xshift=-.3cm] {$0$};
\filldraw (0,1.5*1.66) circle (3pt) node[left] {\small $iK(1{-}m)\,$};
\filldraw (1.5*2.25,1.5*1.66) circle (3pt) node[above] {\small $K(m)+iK(1{-}m)$};
\filldraw (1.5*2.25,0) circle (3pt) node[below,yshift=-.1cm] {\small $K(m)$};
\node at (0.8,1) {$\wp^{-1}(V)$};
\end{tikzpicture}
\hfill
\includegraphics[width=0.50\textwidth]{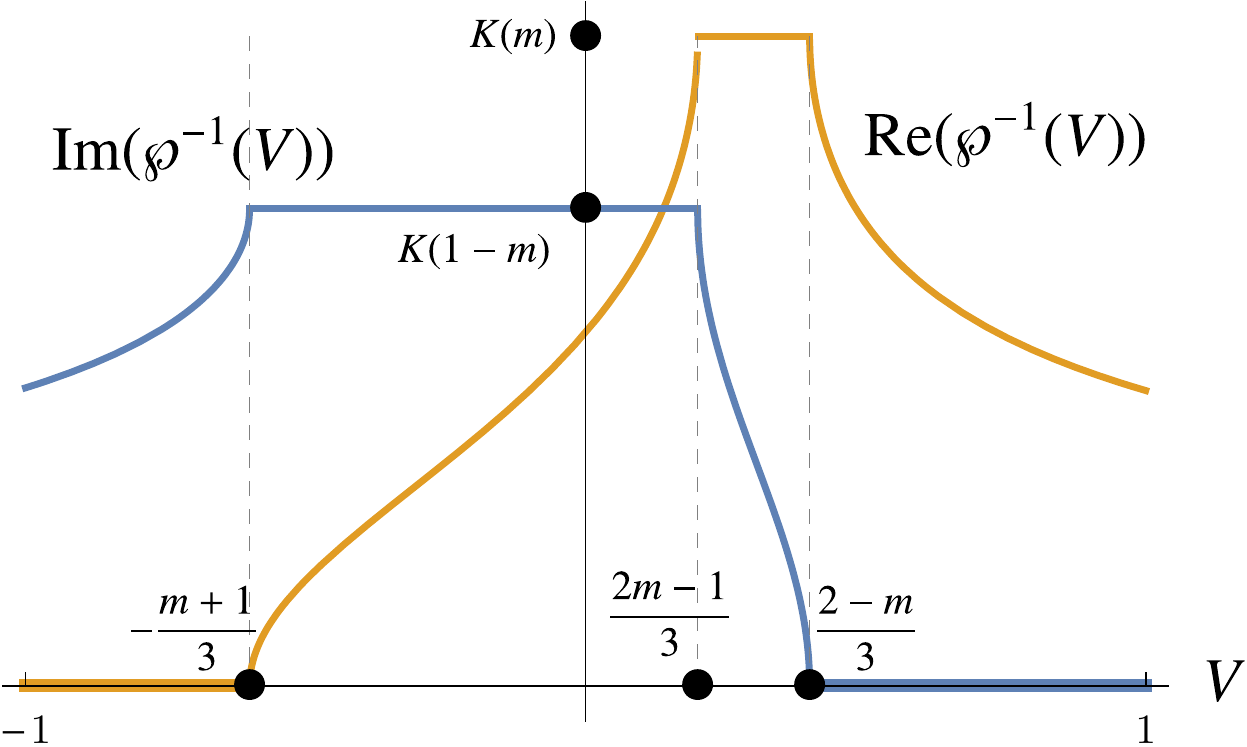}
\caption{{\it Left panel:} At fixed $m$ (here $m=0.8$), the function $\wp^{-1}(V)$ traces a closed path in the complex plane as $V\in\overline{\RR}$ ranges from ${-}\infty$ to ${+}\infty$. Up to standard ambiguities due to the multi-valuedness of $\wp^{-1}$, the path coincides with the rectangle (\ref{RACT}). {\it Right panel:} The real and imaginary parts of the function $\wp^{-1}(V)$. For $V\leq{-}\tfrac{m+1}{3}$, $\wp^{-1}(V)\in i\RR$; for ${-}\tfrac{m+1}{3}\leq V\leq\tfrac{2m-1}{3}$, $\wp^{-1}(V)\in\RR+iK(1{-}m)$; for $\tfrac{2m-1}{3}\leq V\leq\tfrac{2-m}{3}$, $\wp^{-1}(V)\in i\RR+K(m)$; finally, for $V\geq\tfrac{2-m}{3}$, $\wp^{-1}(V)\in\RR$. These four regions correspond to four families of Virasoro orbits --- respectively elliptic with non-zero winding, hyperbolic with unit winding, elliptic with zero winding, and hyperbolic with zero winding.\label{s67}}
\end{figure}
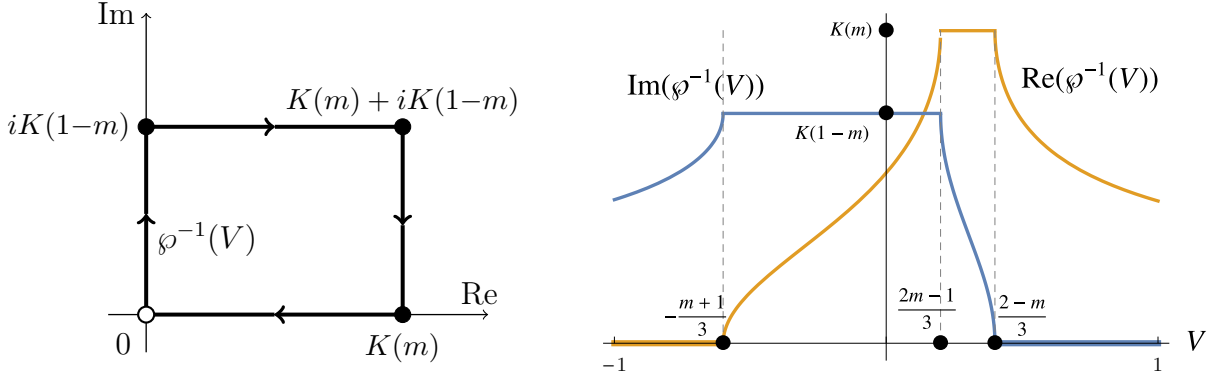

\paragraph{The inverse \texorpdfstring{$\boldsymbol{\wp}$}{$\wp$} function.} Let us investigate the uniform representative (\ref{ss159}) as a function of $V$, thinking of $m\in[0,1)$ as a fixed parameter. In particular, we need to study the values of the inverse Weierstrass function
\be
\wp^{-1}\bigl(V,K(m),iK(1{-}m)\bigr),
\qquad V\in\RR.
\label{s63}
\ee
Since $V$ is real, (\ref{s63}) can only take values $X$ for which $\wp(X)$ is real; this implies that $X$ belongs to the grid\footnote{See eq.\ (\ref{b43b}) in appendix \ref{appB}.}
\be
\bigl(\RR+iK(1{-}m)\ZZ\bigr)\cup\bigl(K(m)\ZZ+i\RR\bigr)
\subset\CC.
\label{grid}
\ee
We may therefore focus on the values of $\wp^{-1}(V)$ belonging to that grid when analysing (\ref{s159}) and (\ref{ss159}). The behaviour of the $\wp$ function on the grid is depicted in fig.\ \ref{s157}. The periods of the $\wp$ function further allow us to restrict attention to values of $\wp^{-1}(V)$ belonging to the following rectangle in the complex plane (see also fig.\ \ref{s67} and \cite{Arancibia:2014vua}):
\be
i[0,K(1{-}m)]
\cup
\bigl([0,K(m)]+iK(1{-}m)\bigr)
\cup
\bigl(K(m)+i[0,K(1{-}m)]\bigr)
\cup
[0,K(m)].
\label{RACT}
\ee
The relation (\ref{s57}) between the $\wp$ function and $\sn^2$ allows us to determine the values of $V$ at the corners of this rectangle: aside from $\wp(0)=\infty$ where $V\rightarrow\pm\infty$, they are\footnote{See eqs.\ (\ref{bb57})-(\ref{b58}) in appendix \ref{appB}.}
\be
\wp(iK(1{-}m))
=
-\frac{m+1}{3},
\quad
\wp(K(m)+iK(1{-}m))
=
\frac{2m-1}{3},
\quad
\wp(K(m))
=
\frac{2-m}{3}.
\ee
These corners will eventually correspond to bifurcations between sharply different types of Virasoro orbits. To see this, we now enumerate the four families of values of $\wp^{-1}(V)$ when $V$ increases from $-\infty$ to $+\infty$, neglecting sign and period ambiguities thanks to the restriction to (\ref{RACT}). These values are plotted in fig.\ \ref{s67}.
\begin{itemize}
\setlength\itemsep{0em}
\item For $V\leq-\tfrac{m+1}{3}$, the function $\wp^{-1}(V)\in i\RR$ has a monotonously increasing imaginary part. When $V=-\tfrac{m+1}{3}$, the derivative of $\wp^{-1}(V)$ diverges and $\wp^{-1}(-\tfrac{m+1}{3})=iK(1{-}m)$. In the opposite limit, $\lim_{V\to-\infty}\wp^{-1}(V)=i0$.
\item For $V\in[{-}\tfrac{m+1}{3},\tfrac{2m{-}1}{3}]$, $\wp^{-1}(V)=X+iK(1{-}m)$ has a monotonously increasing real part $X$. The derivative of $\wp^{-1}(V)$ diverges at the endpoints of the interval, with $\wp^{-1}({-}\tfrac{m+1}{3})=iK(1{-}m)$ and $\wp^{-1}(\tfrac{2m{-}1}{3})=K(m)+iK(1{-}m)$.
\item For $V\in[\tfrac{2m-1}{3},\tfrac{2-m}{3}]$, $\wp^{-1}(V)=K(m)+iY$ has a monotonously decreasing imaginary part $Y$. At the endpoints of the interval, the derivative of $\wp^{-1}(V)$ diverges, with $\wp^{-1}(\tfrac{2m{-}1}{3})=K(m)+iK(1{-}m)$ and $\wp^{-1}(\tfrac{2-m}{3})=K(m)$.
\item For $V\geq\tfrac{2-m}{3}$, $\wp^{-1}(V)\in\RR$ is monotonously decreasing. When $V=\tfrac{2-m}{3}$, its derivative diverges and $\wp^{-1}(\tfrac{2-m}{3})=K(m)$. In the opposite limit, $\lim_{V\to+\infty}\wp^{-1}(V)=0$.
\end{itemize}

\paragraph{Cnoidal orbits.} The different families of values of the inverse Weierstrass function (\ref{s63}), on different domains of $V\in\RR$, allow us to describe in detail the Virasoro orbits appearing in the bifurcation diagram of fig.\ \ref{s60b}. To see this, we list the values of the trace (\ref{s159}) and the constant representative (\ref{ss159}) when $V$ ranges from ${-}\infty$ to $+\infty$:
\begin{itemize}
\setlength\itemsep{0em}
\item For $V\leq{-}\tfrac{m+1}{3}$, the function $\wp^{-1}(V)$ is purely imaginary, as is the combination
\be
K(m)\zeta\big(\wp^{-1}(V)\big)
-\zeta\big(K(m)\big)\wp^{-1}(V).
\label{s68}
\ee
The trace (\ref{s159}) thus satisfies $|\text{Tr}\,\sfM|\leq2$, corresponding to elliptic Virasoro orbits. The constant (\ref{ss159}) is such that $k/c<0$, with $k/c$ reaching its upper bound at $V={-}\tfrac{m+1}{3}$, where $k=-c/24$. Indeed, at that point, the combination (\ref{s68}) becomes\footnote{See eq.\ (\ref{bb45}) in appendix \ref{appB}.}
\be
K(m)\zeta\big(iK(1{-}m)\big)-\zeta\big(K(m)\big)iK(1{-}m)
=
-i\pi/2.
\label{ss68}
\ee
Thus the boundary of the region $V\leq{-}\tfrac{m+1}{3}$ in the $(m,V)$ plane is the exceptional $n=1$ Virasoro orbit. In fig.\ \ref{s62}, the region $V<{-}\tfrac{m+1}{3}$ spans the vertical line below the point $n=1$. There is no lower bound on $k/c$.
\item In the region ${-}\tfrac{m+1}{3}\leq V\leq\tfrac{2m{-}1}{3}$, one has $\wp^{-1}(V)=X+iK(1{-}m)$ with real $X$, so the combination (\ref{s68}) takes the form
\be
K(m)\zeta\big(X+iK(1{-}m)\big)-\zeta\big(K(m)\big)(X+iK(1{-}m))
\stackrel{\text{(\ref{ss68})}}{=}
-\frac{i\pi}{2}+(\text{real}),
\nn
\ee
where $(\text{real})$ is some real function of $X$.\footnote{To be precise, $(\text{real})=K(m)\zeta(X)-\zeta(K(m))X+K(m)\wp'(X)/(2\wp(X)+(2m+2)/3)$.} To derive this we used the addition formula (\ref{bb51}) along with $\wp'(iK(1{-}m))=0$. Accordingly, the trace (\ref{s159}) is
\be
\text{Tr}\,\sfM
=
-2\cosh[2(\text{real})]
\leq-2
\nn
\ee
and the constant (\ref{ss159}) has a non-vanishing imaginary part. Thus the orbits of the bifurcation wedge ${-}\tfrac{m+1}{3}<V<\tfrac{2m{-}1}{3}$ are hyperbolic with an odd winding number; they have no rest frame, \ie no uniform representative. In fact, as argued in section \ref{sec31}, the winding number must be $n=1$, so the orbits in that wedge span the horizontal line at $n=1$ in fig.\ \ref{s62}. On the boundaries of the wedge, the constant (\ref{ss159}) takes the value $k=-c/24$. We have already shown this in (\ref{ss68}) for the lower boundary, while for the upper one we use\footnote{See eq.\ (\ref{APAKI}) in appendix \ref{appB}.}
\begin{align}
&
K(m)\zeta(\wp^{-1}(\tfrac{2m{-}1}{3}))-\zeta(K(m))\wp^{-1}(\tfrac{2m{-}1}{3})=\nn\\
&=
K(m)\zeta(K(m)+iK(1{-}m))-\zeta(K(m))(K(m)+iK(1{-}m))=-i\pi/2,\label{ipitu}
\end{align}
which yields $k=-c/24$ as announced. Thus, we have now confirmed that both boundaries of the forbidden wedge consist of a single exceptional orbit:
\be
\boxed{
k
=
-\frac{c}{24}
\qquad\text{for}\qquad
V
=
\frac{(1\pm3)m-2}{6}}
\qquad
\text{(wedge boundaries).}
\label{KOMEGGA}
\ee
\item In the region $\tfrac{2m{-}1}{3}\leq V\leq\tfrac{2{-}m}{3}$, one has $\wp^{-1}(V)\in i\RR+K(m)$, so the combination (\ref{s68}) is purely imaginary. (To prove this, use again the addition formula (\ref{bb51}).) Hence the trace (\ref{s159}) is such that $|\text{Tr}\,\sfM|\leq2$ and the orbits in that region are elliptic. The corresponding constant (\ref{ss159}) is such that $k/c\leq0$, with $k$ ranging monotonically from $k=-c/24$ at $V=\tfrac{2m{-}1}{3}$, to $k=0$ at $V=\tfrac{2{-}m}{3}$. Thus the orbits in the region $\tfrac{2m{-}1}{3}<V<\tfrac{2{-}m}{3}$ span the vertical interval between the exceptional point $n=1$ and the orbit of $k=0$ in the map of fig.\ \ref{s62}.
\item When $V\geq\tfrac{2{-}m}{3}$, the function $\wp^{-1}(V)$ and the combination (\ref{s68}) are both real, so the trace (\ref{s159}) is such that $\text{Tr}\,\sfM\geq2$, corresponding to hyperbolic Virasoro orbits with a constant representative such that $k/c\geq0$. The bound $k=0$ is reached at $V=\tfrac{2{-}m}{3}$, and there is no upper bound on $k/c$. Thus the orbits with $V>\tfrac{2{-}m}{3}$ span the vertical line above the point $k=0$ in fig.\ \ref{s62}.
\end{itemize}
This analysis covers the entire parameter space of cnoidal profiles and provides a detailed picture of the bifurcation diagram of fig.\ \ref{s60b}. The results are summarized in fig.\ \ref{s71b}.\\

\begin{figure}[h!]
\captionsetup{singlelinecheck=off}
\centering
\includegraphics[width=0.60\textwidth]{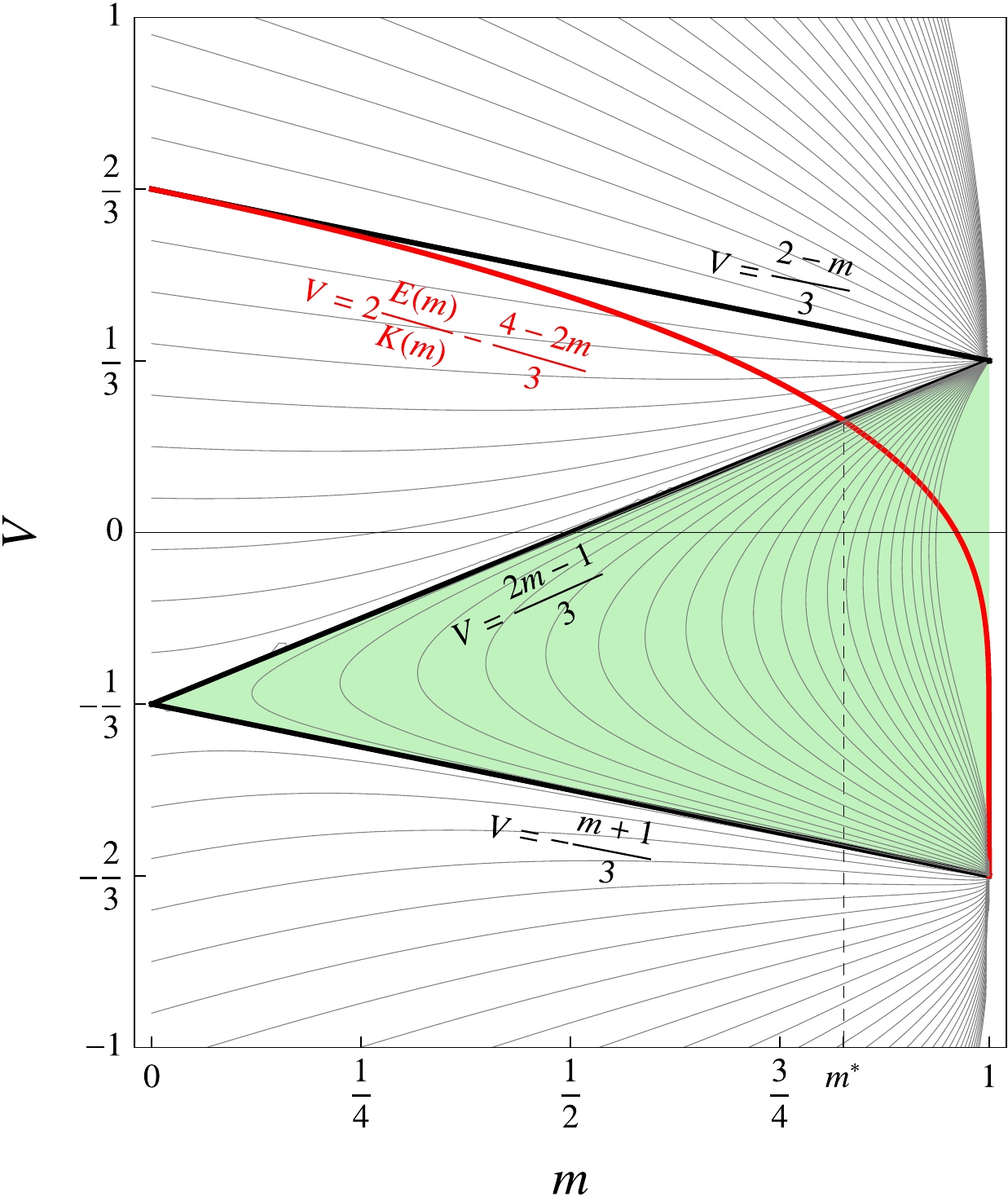}
\caption[foo bar]{The bifurcation diagram of fig.\ \ref{s60b}, now represented with a few level curves of $k$ and sharp boundaries between different classes of Virasoro orbits. As $V$ increases from $-\infty$ to $+\infty$ (at fixed $m\in[0,1)$), four distinct regions are apparent, corresponding to the four types of values of the inverse Weierstrass function (\ref{s63}):
\begin{itemize}
\setlength\itemsep{0em}
\item The region $V<{-}\tfrac{m+1}{3}$ comprises elliptic orbits with non-zero winding (the line below the exceptional point $n=1$ in fig.\ \ref{s62}).
\item The interval ${-}\tfrac{m+1}{3}<V<\tfrac{2m{-}1}{3}$ comprises hyperbolic orbits with unit winding (the horizontal line at $n=1$ in fig.\ \ref{s62}), which have no constant representative and span the green `forbidden wedge'; the curves in the wedge are level sets of $\text{Tr}\,\sfM<-2$, with $\text{Tr}\,\sfM$ decreasing monotonously from left to right.
\item The interval $\tfrac{2m{-}1}{3}<V<\tfrac{2{-}m}{3}$ consists of elliptic orbits with zero winding (the segment between $n=1$ and $k=0$ in fig.\ \ref{s62}).
\item $V>\tfrac{2{-}m}{3}$ comprises hyperbolic orbits with zero winding ($k/c>0$ in fig.\ \ref{s62}). 
\end{itemize}
There are three boundaries between these four regions:
\begin{itemize}
\setlength\itemsep{0em}
\item The (singular) boundary $V=[(1\pm3)m-2]/6$ of the `forbidden wedge' contains a single exceptional orbit, that of $k=-c/24$.
\item The (smooth) boundary $V=\tfrac{2{-}m}{3}$ contains the parabolic orbit of $k=0$.
\end{itemize}
For later reference, we also include the red curve that consists of cnoidal waves with vanishing average, given by eq.\ (\ref{ss31}). It is asymptotic to $V\sim\tfrac{2{-}m}{3}$ at small $m$ and lies in the forbidden wedge for sufficiently large $m$, eventually converging to $V=-2/3$ as $m\to1$. The crossing value $m^*=0.8261...$ is computed below eq.\ (\ref{s37b}).\label{s71b}}
\end{figure}
 \clearpage

We still need to address one last key feature of figs.\ \ref{s60b} and \ref{s71b}: as one can see \eg from the plot of $k$ in terms of $(m,V)$ in fig.\ \ref{s60b}, the derivative of $k$ with respect to $V$ diverges at the boundaries of the wedge, but there is no such divergence on the line $k=0$ separating elliptic and hyperbolic orbits with zero winding --- despite the fact that the derivative of $\wp^{-1}(V)$ does diverge there. To explain this, we compute the derivative of (\ref{ss159}),
\be
\frac{\der k}{\der V}
=
-\frac{c}{3\pi^2}
\Bigl(
K(m)\zeta\bigl(\wp^{-1}(V)\bigr)-\zeta(K(m))\wp^{-1}(V)
\Bigr)
\biggl(K(m)V+\frac{\zeta(K(m))}{\wp'\bigl(\wp^{-1}(V)\bigr)}\biggr),
\label{s72}
\ee
and investigate how it behaves at the three bifurcation lines in fig.\ \ref{s71b}. On the boundaries (\ref{KOMEGGA}) of the wedge, the first factor of (\ref{s72}) is non-zero (it equals $-i\pi/2$) while the second factor diverges due to the vanishing derivative $\wp'(\wp^{-1}(V))=\wp'(iK(1{-}m))=0$, resulting in the observed singularities
\be
\frac{1}{c}
\left.
\frac{\der k}{\der V}
\right|_{V=[(1\pm3)m-2]/6}
=
+\infty.
\nn
\ee
But the boundary between elliptic and hyperbolic orbits with zero winding, at $V=\tfrac{2{-}m}{3}$, is different: there the first factor of (\ref{s72}) vanishes, which cancels the divergence of $1/\wp'(K(m))=\infty$. In fact, using L'H\^opital's rule, one readily finds that the derivative of $k$ with respect to $V$ is finite on that line:\footnote{See eq.\ (\ref{KAPPOT}) in appendix \ref{appB} for the proof of $\wp''(K(m))=2(1-m)$.}
\be
\left.
\frac{\der k}{\der V}
\right|_{V=\tfrac{2{-}m}{3}}
=
\frac{c}{3\pi^2}\frac{\zeta\big(K(m)\big)^2}{\wp''\big(K(m)\big)}
=
\frac{c}{6\pi^2}\frac{\zeta\big(K(m)\big)^2}{1-m}.
\label{finidev}
\ee
This explains all the qualitative features of figs.\ \ref{s60b} and \ref{s71b}. As we now show, the finiteness of (\ref{finidev}) is a key qualitative difference between uniform representatives of cnoidal Virasoro orbits and crystal momenta in the Lam\'e band structure.

\subsection{Lam\'e band structure and multiple solitons}
\label{sec34}

\begin{figure}[t]
\centering
~~~~~~~~~~~~~
\includegraphics[width=0.3\textwidth]{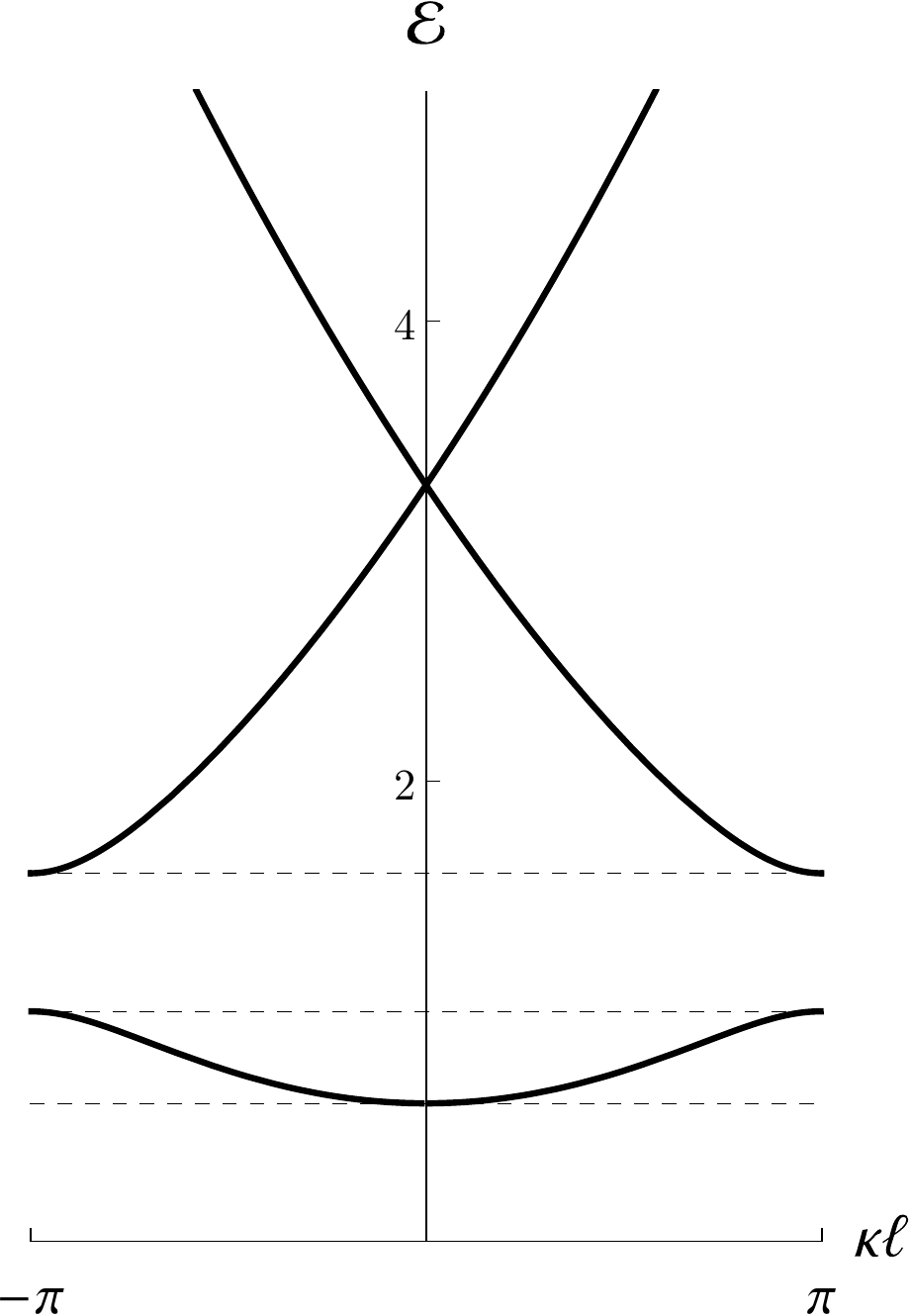}
\hfill
\includegraphics[width=0.3\textwidth]{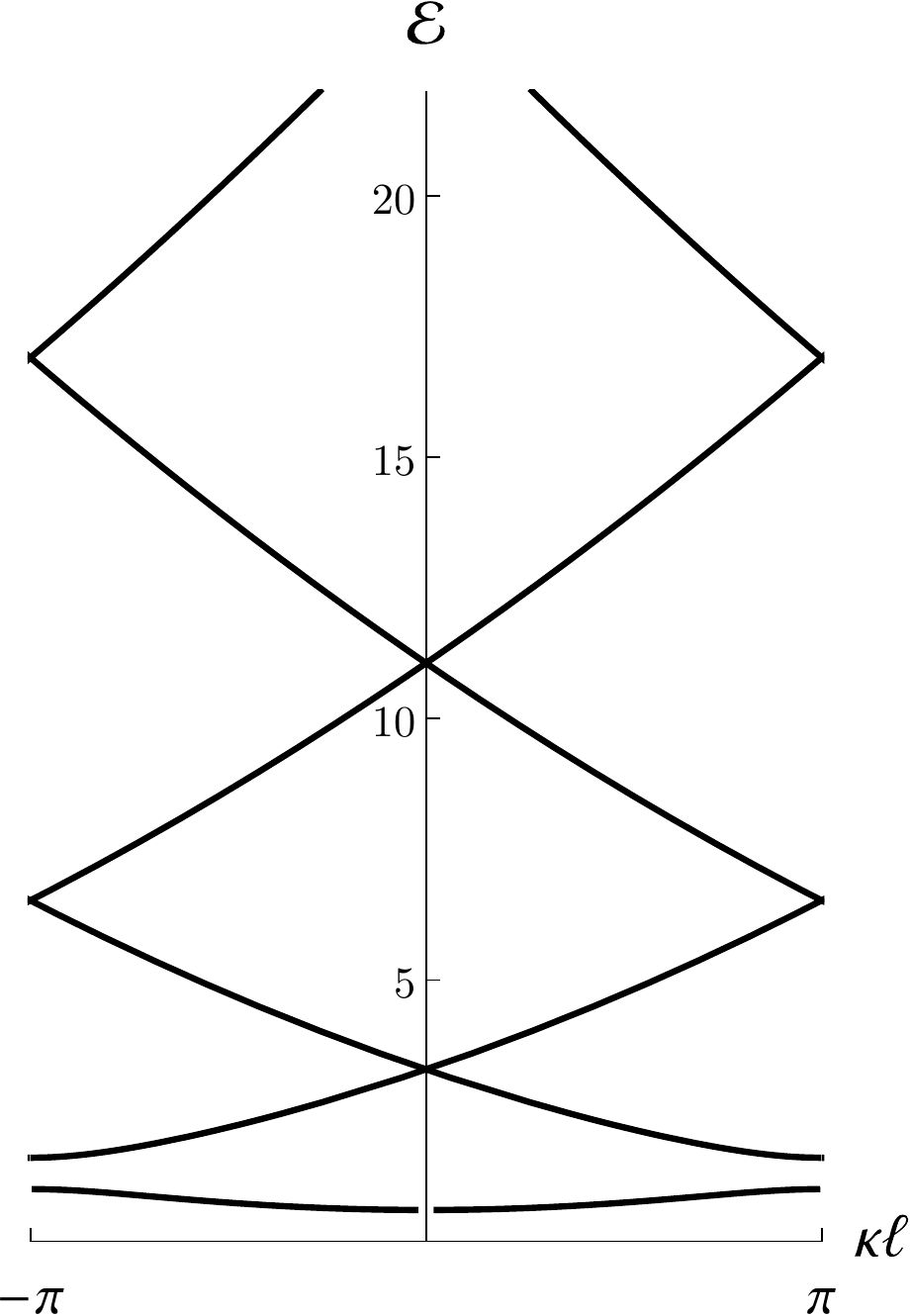}
~~~~~~~~~~~~~
\caption{The band structure (\ref{s20b}) of the $N=1$ Lam\'e potential at $m=0.6$, with crystal momenta restricted to the first Brillouin zone. There is one `valence band' where $\cE\in[m,1]$, and one semi-infinite `conduction band' where $\cE\in[m+1,+\infty)$. The two are separated by a gap of width $m$. The valence and conduction bands respectively correspond to elliptic Virasoro orbits with zero and non-zero windings, while the gap corresponds to the `forbidden wedge' of hyperbolic orbits with unit winding. In the {\it left panel}, band edges are indicated by dashed lines. The {\it right panel} displays the band structure at higher energies, where the details of the potential become irrelevant. Exceptional Virasoro orbits sit at points such that $\kappa\ell\in\{0,\pi\}$ (modulo $2\pi$), whose energies are asymptotic to (\ref{enerasym}).\label{ss19t}}
\end{figure}

The Lam\'e equation played an essential role in the derivation of fig.\ \ref{s71b}, so we now relate the classification of orbits of cnoidal waves to the $N=1$ Lam\'e band structure \cite{Whittaker,Arscott}. The starting point is a particle of mass $M$ in one dimension, subject to the potential
\be
U(X)
=
\frac{4\hbar^2}{M\ell^2}K(m)^2m\sn^2\Bigl(\frac{2K(m)}{\ell}X\Big|m\Bigr)
\label{b19b}
\ee
with `lattice spacing' $\ell$. In terms of $y\equiv 2K(m)X/\ell$, the corresponding time-independent Schr\"odinger equation with energy $E$ is the Lam\'e equation (\ref{s153}) with
\be
V
=
\frac{2m+2}{3}
-\frac{M\ell^2}{2\hbar^2K(m)^2}E
\equiv
\frac{2m+2}{3}-\cE.
\label{s19t}
\ee
Assuming now that the potential (\ref{b19b}) is that of a one-dimensional lattice, the wavefunction $\psi(y)$ is generally not $2\pi$-periodic, but has instead a certain monodromy determined by $m$ and $E$. This is customarily written in terms of a crystal momentum $\kappa\in\RR$:
\be
\psi(y+2K(m))
=
e^{i\kappa\ell}\psi(y).
\label{bb19b}
\ee
With the terminology of section \ref{sec22}, this yields an elliptic monodromy for $\psi$. Indeed, the function $\psi^*(y)$ satisfies the same Schr\"odinger equation as $\psi$, but with crystal momentum $-\kappa$ instead of $\kappa$. The trace of the monodromy of $\psi,\psi^*$ is thus $\text{Tr}\,\sfM=2\cos(\kappa\ell)$, which indeed corresponds to an elliptic matrix $\sfM$.\footnote{This is the generic case. In the exceptional case $\kappa\ell=n\pi$ mod $2\pi$, $\psi^*$ and $\psi$ have the same crystal momentum and may not be independent. The result $\text{Tr}\,\sfM=2\cos(\kappa\ell)$ is nevertheless correct even then, except that the linearly independent solutions used to compute $\text{Tr}\,\sfM$ may not be $\psi$ and $\psi^*$.} Since the monodromy is that of the Lam\'e equation, we can express $\text{Tr}\,\sfM$ as in eq.\ (\ref{s159}). Using eq.\ (\ref{s19t}), this implies
\be
\kappa\ell
=
\pm 2i\Bigl[
K(m)\zeta\bigl(\wp^{-1}(\tfrac{2m+2}{3}-\cE)\bigr)
-\zeta\bigl(K(m)\bigr)\wp^{-1}\bigl(\tfrac{2m+2}{3}-\cE\bigr)
\Bigr]+2\pi n.
\label{s20b}
\ee
At fixed $m$, this is a relation between the allowed energy $E=\tfrac{2\hbar^2K(m)^2}{M\ell^2}\cE$ and the corresponding crystal momentum $\kappa$. Values of energy for which $\kappa$ acquires an imaginary part are forbidden.\footnote{This need not be true when edge modes are included, but we do not consider this possibility here.} The resulting band structure is depicted in fig.\ \ref{ss19t}.\\

Eq.\ (\ref{s20b}) for the Lam\'e band structure is the square root of eq.\ (\ref{ss159}) for uniform representatives of cnoidal waves:
\be
\frac{6k}{c}
=
-\biggl(\frac{\kappa\ell}{2\pi}\biggr)^2,
\qquad
V=\frac{2m+2}{3}-\cE.
\label{squarek}
\ee
In particular, owing to the relation between the velocity $V$ and the energy $\cE$, the band edges $\cE\in\{m{+}1,1,m\}$ respectively correspond to $V\in\big\{{-}\tfrac{m+1}{3},\tfrac{2m-1}{3},\tfrac{2-m}{3}\big\}$; those are the boundaries of regions with elliptic orbits in fig.\ \ref{s71b}. More generally, in fig.\ \ref{s71b}, level curves of $k$ are also lines along which the crystal momentum $\kappa$ is constant. Non-elliptic orbits are excluded by the Bloch ansatz (\ref{bb19b}) --- at least as long as $\kappa$ is real. Thus, fig.\ \ref{s71b} is consistent (as expected) with the band structure of the $N=1$ Lam\'e equation. A key difference, however, occurs at $\cE=m$, the lower edge of the valence band, corresponding to $V=\tfrac{2{-}m}{3}$. Indeed, we saw in eq.\ (\ref{finidev}) that this line is not singular in terms of Virasoro orbits (despite its being a bifurcation between elliptic and hyperbolic orbits), and in particular that $\der k/\der V$ is finite there. By contrast, the derivative of the crystal momentum $\kappa$ with respect to $\cE$ does diverge at $\cE=m$, since the latter is a band edge. This is precisely because the uniform representative $k$ is the {\it square} (\ref{squarek}) of the crystal momentum $\kappa$: the square is responsible for the cancellation leading to the finite value (\ref{finidev}).\\

Note that the band structure of fig.\ \ref{ss19t} also exhibits points where $\kappa\ell=n\pi$ mod $2\pi$ with integer $n$, corresponding to a monodromy $\sfM=(-1)^n\II$. These points are exceptional Virasoro orbits. Their energy $E$ behaves at large $n$ as
\be
E\sim\frac{\hbar^2n^2\pi^2}{2M\ell^2}
\qquad
\text{as}
\qquad
n\rightarrow+\infty.
\label{enerasym}
\ee
For even $n$, those are the energy levels of a free particle on a circle with radius $\ell$: particles with high energy are insensitive to the potential, as expected.

\paragraph{Higher Lam\'e equations and multiple solitons.} For any periodic profile $p(x)$ depending on certain parameters, it is quite generally true that band edges correspond to exceptional Virasoro orbits. Indeed, band edges occur when the crystal momentum satisfies $\kappa\ell=n\pi$ (in the `extended zone scheme' \cite{Ashcroft}) and yields a uniform orbit representative $k=-n^2c/24$. This observation is especially helpful for the orbital analysis of the Lam\'e equation with a higher value of $N\in\mathbb{N}$,
\be
-\frac{\dd^2\Psi}{\dd y^2}
+N(N+1)m\sn^2(y|m)\Psi
=
\cE\Psi.
\label{LAMBADA}
\ee
Indeed, at fixed $m\in[0,1)$ the latter famously has exactly $N$ band gaps \cite{Ince,Ince2,Whittaker,Arscott}, \ie $N$ separate `valence bands' and one semi-infinite `conduction band'. The gaps vanish in the limit $m\to0$, while valence bands vanish in the opposite limit $m\to1$. In terms of Virasoro orbits, this corresponds to a bifurcation diagram in the $(m,\cE)$ plane that displays exactly $N$ `forbidden wedges' in which the profile
\be
p(x)
=
\frac{cK(m)^2}{6\pi^2}
\biggl[
N(N+1)m\sn^2\biggl(\frac{K(m)}{\pi}x\bigg|m\biggr)-\cE
\biggr]
\label{PAKA}
\ee
is not conformally equivalent to any uniform field configuration. In short, at fixed $m$, (\ref{PAKA}) has a Virasoro rest frame if and only if $\cE$ is not in a band gap.\\

It is easy to understand qualitatively the Virasoro orbits that are spanned by the profiles (\ref{PAKA}) as functions of $(m,\cE)$. Indeed, at $m=0$, the Lam\'e equation (\ref{LAMBADA}) is that of a free particle and (degenerate) band edges occur when $\kappa\ell=n\pi$, \ie when $k=-n^2c/24$, which in turn corresponds to $\cE=n^2$. At small $m$, gaps start opening around the values $\cE=1^2,2^2,3^2,...,N^2$. These are the $N$ `forbidden wedges'. On the boundaries of the $N^{\text{th}}$ wedge, all profiles have a rest frame with uniform representative $k=-n^2c/24$. The remainder of the resulting `bifurcation diagram' can be completed by continuity:
\begin{itemize}
\setlength\itemsep{0em}
\item Inside the $n^{\text{th}}$ wedge, wave profiles have no rest frame. They belong to hyperbolic Virasoro orbits with winding $n$. In figs.\ \ref{s19} and \ref{s62}, they span the $n^{\text{th}}$ horizontal line.
\item Between the $n^{\text{th}}$ and $(n{+}1)^{\text{th}}$ wedges, wave profiles have a uniform representative $k$ such that ${-}\tfrac{(n{+}1)^2}{24}<k/c<{-}\tfrac{n^2}{24}$. They belong to elliptic orbits with winding $n$. In figs.\ \ref{s19} and \ref{s62}, they span the interval between the points $n$ and $n{+}1$.
\item For $\cE$ above the highest ($N^{\text{th}}$) band gap, wave profiles have winding $\geq N$ and a uniform representative $k$ such that $k/c<{-}\tfrac{N^2}{24}$, corresponding to elliptic Virasoro orbits. In figs.\ \ref{s19} and \ref{s62}, they span the vertical line below the point $N$.
\item Conversely, for $\cE$ below the lowest ($n=1$) gap, wave profiles have winding zero and a uniform representative $k$ such that $k/c>{-}\tfrac{1}{24}$, with $k/c\leq0$ if $\cE$ belongs to the lowest allowed band and $k/c>0$ otherwise. This respectively corresponds to elliptic orbits when $\cE\in$ (band) and hyperbolic orbits when $\cE<$ (lowest band edge), with the transition at the band edge $k=0$ being a parabolic orbit. In figs.\ \ref{s19} and \ref{s62}, these profiles span the vertical line above the point $n=1$.
\end{itemize}
Thus, in terms of figs.\ \ref{s19} and \ref{s62}, the Lam\'e potentials (\ref{PAKA}) with fixed $N$ span the entire vertical line and the $N$ highest horizontal lines. The latter consist of orbits of profiles that lie in Lam\'e band gaps, while the former only contains cnoidal waves that lie in allowed energy bands. The value of $\cE$ increases monotonically as $k/c$ decreases. This is a natural generalization of the $N=1$ results summarized in fig.\ \ref{s62}.\\

Higher $N$ Lam\'e potentials, and finite-gap potentials more generally, are closely related to periodic multi-soliton solutions of the KdV equation. In fact, one can show that the solution of the KdV equation (\ref{s5}) with initial condition (\ref{PAKA}) is a sum of $N(N+1)/2$ cnoidal waves (\ref{CNOMEGA}) with specific {\it time-dependent} velocities $V_i(\tau)$, $i=1,...,N(N+1)/2$, such that $\sum_iV_i=\tfrac{m+1}{3}N(N+1)-\cE$ \cite{Airault}. In that sense, the classification of Virasoro orbits for $N$-gap Lam\'e potentials is also a classification of orbits of superpositions of $N(N+1)/2$ cnoidal waves. In particular, the values of $(m,\cE)$ that belong to one of the $N$ forbidden wedges specify wave profiles that lack a rest frame at any point along KdV time evolution. Superpositions of non-triangular numbers ($\neq N(N+1)/2$) of cnoidal waves can also be achieved using finite-gap potentials of a more general form than the Lam\'e expression (\ref{PAKA}) \cite{Verdier}. The Virasoro orbits corresponding to such profiles are, again, exceptional along band edges, elliptic in allowed bands, and hyperbolic with non-zero winding in band gaps. As usual, the latter have no rest frame.

\section{Asymptotics of level curves}
\label{sec36}

A striking aspect of fig.\ \ref{s60b} is the fact that, at large $|V|$, $k$ depends linearly on $V$. A similarly striking feature of fig.\ \ref{s71b} is the convergence of level curves of $k$ to the points $V=-2/3$ and $V=+1/3$ as $m\to1$. We now explain both of these phenomena, and others, by an asymptotic analysis of level curves in various regions of the $(m,V)$ plane: we first look at the regimes where $|V|$ goes to infinity, then turn to regions close to the bifurcation wedge, then analyse the singularities at $m\to1$, and finally investigate the region $m\to0$. This will result in various formulas yielding the approximate dependence of $V$ on $m$ (or vice-versa) along level curves; in particular, as we shall see, the singularity at $(m=1,V=+1/3)$ is qualitatively very different from the one at $(m=1,V=-2/3)$.\\

Note that this section merely serves to develop a detailed understanding of the level curves in fig.\ \ref{s71b}, but will otherwise be unimportant for the rest of the paper. In particular, readers interested in shoaling may want to go straight to section \ref{sec4}.\\

Other asymptotic analyses of the Lam\'e band structure have previously appeared in the literature, including its high-energy regime and various expansions of periodic wavefunctions at band edges \cite{Muller}, some of which have recently found applications in the context of resurgence \cite{Basar}. The limits $m\to0$ and $m\to1$ have also been studied \cite{Ince}, in particular owing to the `perturbative/non-perturbative duality' that relates them \cite{Dunne:1999zc}. Finally, the large $N$ limit of the general Lam\'e equation (\ref{LAMBADA}) is of interest as well, partly owing to its relation with the Mathieu equation \cite{Whittaker}. Aside from the latter limit which we do not consider at all, our results will indeed partially overlap the aforementioned references --- though, to our knowledge, some of the limits we obtain in section \ref{secNEAR} are new. As before we refer to appendix \ref{appB} for some background on elliptic functions, and to appendix \ref{appC} for their asymptotic behaviour.

\subsection{Large \texorpdfstring{$\boldsymbol{|V|}$}{$|V|$}}

Consider eq.\ (\ref{ss159}) in the limit where $|V|$ is large at fixed $m\in[0,1)$. When $V\to-\infty$, the inverse Weierstrass function (\ref{s63}) is purely imaginary (up to standard ambiguities). Since the $\wp$ function reduces to\footnote{See eq.\ (\ref{bb40}) in appendix \ref{appB}.} $1/z^2+\cO(z^2)$ for small $|z|$, one has $\wp^{-1}(V)\sim i|V|^{-1/2}(1+\cO(1/V^2))$ as $V\to-\infty$, while the Weierstrass $\zeta$ function is such that\footnote{See eq.\ (\ref{bb44}) in appendix \ref{appB}.} $\zeta(z)\sim1/z+\cO(z^3)$ as $z\to0$, so $\zeta\bigl(\wp^{-1}(V)\bigr)\sim-i|V|^{1/2}(1+\cO(1/V^2))$ as $V\to-\infty$. Similarly, when $V\to+\infty$, the inverse Weierstrass function (\ref{s63}) is purely real (again up to standard ambiguities). Then $\wp^{-1}(V)\sim V^{-1/2}(1+\cO(1/V^2))$ and $\zeta(\wp^{-1}(V))\sim V^{1/2}(1+\cO(1/V^2))$ as $V\to+\infty$. Thus, in both limits, eq.\ (\ref{ss159}) reduces to\footnote{We use the symbol $\sim$ to denote standard asymptotic equalities and asymptotic series. In particular, $f(x)\sim g(x)$ as $x\to a$ iff $\lim_{x\to a}f(x)/g(x)=1$.}
\begin{align}
\frac{6\pi^2 k}{c}
&\sim
\pm
\biggl[
K(m)\sqrt{|V|}
\mp
\frac{\zeta\big(K(m)\big)}{\sqrt{|V|}}
\biggr]^2
\big(1+\cO(1/V^2)\big)
\quad\text{as}\quad V\to\pm\infty
\label{b18b}\\
&\sim
K(m)^2V-2K(m)\zeta\big(K(m)\big)+\cO(1/|V|)
\quad\text{as}\quad V\to\pm\infty.
\label{KARAPA}
\end{align}
This is a linear dependence on $V$, as expected from fig.\ \ref{s60b}, with the same slope and the same additive constant in both limits $V\to\pm\infty$. (In fact, the asymptotic series of $k$ in powers of $1/V$ takes the exact same form in both regions $V\to\pm\infty$.) The slope increases with $m$, eventually becoming infinite as $m\to1$; we return to the latter limit in much greater detail below.\\

As a by-product, we can investigate the asymptotics of the curves where the uniform representative (\ref{ss159}) takes exceptional values $k=-n^2c/24$ with $n\geq2$. In principle, such curves are bifurcations from the point of view of Virasoro orbits and their symplectic structure, just as much as the lines $V=[(1\pm3)m-2]/6$ where $k=-c/24$. However, in contrast to the latter case, these curves are completely regular. In any case, they are all contained in the region $V<{-}\tfrac{m+1}{3}$. At fixed $m$ and large $n$, $V$ goes to $-\infty$ and eq.\ (\ref{b18b}) yields a result analogous to eq.\ (\ref{enerasym}):
\be
k=-\frac{n^2c}{24}
\qquad
\Leftrightarrow
\qquad
V\sim-\frac{\pi^2n^2}{4K(m)^2}
\quad\text{as}\quad
\frac{n}{K(m)}\to+\infty.
\label{XAMAN}
\ee
Note that here we write a limiting condition on $n/K(m)$ rather than $n$ alone. This is because the key fact in deriving (\ref{XAMAN}) was that $\lim_{n\to+\infty}V=-\infty$ at fixed $m$, as the simplification (\ref{b18b}) only works when $|V|$ is very large, which in turn occurs only when $n\gg K(m)$. In other words, the approximation (\ref{XAMAN}) is valid only in the lower left corner of the $(m,V)$ plane, where the increasing function $K(m)$ is small enough. By contrast, when $m$ gets close to $1$, the function $K(m)$ grows without bound and the asymptotic formula (\ref{XAMAN}) breaks down, no matter how large $n$ is.

\subsection{Near the bifurcations}

We now investigate the behaviour of level curves when $V$ approaches its values $[(1\pm3)m-2]/6$ on the boundaries of the bifurcation wedge, along which $k=-c/24$. We assume $m\in[0,1)$ to be fixed throughout. The computation works in the same way both when $V\to{-}\tfrac{m+1}{3}$ from below and when $V\to\tfrac{2m{-}1}{3}$ from above, so we only describe the former. Namely, using\footnote{See eq.\ (\ref{KAPPOT}) in appendix \ref{appB}.} $\wp''(iK(1{-}m))=2m$ one finds
\be
\wp^{-1}(-\tfrac{m+1}{3}-\nu)
=
iK(1{-}m)-i\sqrt{\frac{\nu}{m}}
+\cO(\nu).
\label{LAOS}
\ee
Using $\zeta'=-\wp$, it also follows that
\be
\zeta\big(\wp^{-1}(-\tfrac{m+1}{3}-\nu)\big)
\,\sim\,
\zeta\big(iK(1{-}m)\big)-\frac{m+1}{3}\,i\,\sqrt{\frac{\nu}{m}}+\cO(\nu)
\qquad\text{as }\nu\to0^+,
\nn
\ee
and this can be plugged into (\ref{ss159}) to obtain an asymptotic relation between $k/c$, $m$ and $V$ near the lower boundary of the bifurcation wedge. A similar argument applies near the upper boundary, resulting in a general asymptotic relation that can be written as
\be
\frac{k}{c}
\sim
-\frac{1}{24}
+
\frac{1}{6\pi}
\Big(\zeta\big(K(m)\big)+\tfrac{(1\pm3)m-2}{6}K(m)\Big)
\sqrt{\frac{\big|2V-\tfrac{(1\pm3)m-2}{3}\big|}{m(2-(1\pm1)m)}},
\qquad V\to\tfrac{(1\pm3)m-2}{6}.
\label{XAM}
\ee
This approximate square root behaviour is apparent in the right panel of fig.\ \ref{s60b}. Note again that the formula breaks down in the limit $m\to1$, to which we now turn.

\subsection{Near \texorpdfstring{$\boldsymbol{m=1}$}{$m=1$}}
\label{secNEAR}

As we have seen, the limit $m\to1$ is problematic: $K(m)$ blows up and, for instance, eqs.\ (\ref{XAMAN}) and (\ref{XAM}) do not apply. This limit is singular, and indeed we were careful not to take $m$ too large in figs.\ \ref{s59b} and \ref{s60b}. In fig.\ \ref{s71b}, the singularity translates into the convergence of level curves of $k$ to the points $V=\{{-}2/3,1/3\}$ as $m\to1$. To prove this convergence we must return to eq.\ (\ref{ss159}) for $k$ and find its asymptotic behaviour in the limit $m\to1$. As we shall see, the result depends crucially on whether $k>-c/24$ or $k<-c/24$.\\

When $m\to1$, the complete elliptic integral of the first kind satisfies \cite[sec.\ 19.12]{DLMF}
\be
K(m)
\,\stackrel{m\to1}{\sim}\,
-\log\sqrt{1-m}+\log 4+\cO\bigl((1{-}m)\log(1{-}m)\bigr).
\label{KAPADI}
\ee
Accordingly, the half-period $K(m)$ in the inverse Weierstrass function $\wp^{-1}$ and the zeta function $\zeta$ of eq.\ (\ref{ss159}) goes to infinity. As for the imaginary period, it remains finite:
\be
K(1{-}m)
\,\stackrel{m\to1}{\sim}\,
\frac{\pi}{2}\Bigl(1+\frac{1-m}{4}\Bigr)
+\cO\bigl((1{-}m)^2\bigr).
\label{KAPAFI}
\ee
Thus, in order to understand the asymptotics of (\ref{ss159}) as $m\to1$, we need to know how Weierstrass functions behave when one of their periods diverges. This follows from standard asymptotic formulas that can be found \eg in \cite[sec.\ 23.12]{DLMF}, and it is addressed in greater detail in appendix \ref{appC}. In particular, one has
\be
\zeta(K(m))
\,\stackrel{m\to1}{\sim}\,
\frac{\pi}{2K(1{-}m)}
-\frac{\pi^2K(m)}{4K(1{-}m)^2}
\biggl[
\frac{1}{3}
+\cO\bigl((1{-}m)^2\bigr)
\biggr]
\label{n20t}
\ee
and more generally, at fixed $z$,
\be
\zeta(z)
\,\stackrel{m\to1}{\sim}\,
-\frac{\pi^2}{4K(1{-}m)^2}
\biggl[
\frac{z}{3}
-\frac{2K(1{-}m)}{\pi}\cotanh\Bigl(\frac{\pi z}{2K(1{-}m)}\Bigr)
+\cO\bigl((1{-}m)^2\bigr)
\biggr].
\label{NEWT}
\ee
Starting from this, it is immediate to prove that any level curve (\ref{ss159}) converges to $V=-2/3$ or $V=+1/3$ as $m\to1$. Indeed, in that limit the curve reduces to
\be
\pm\sqrt{\frac{6\pi^2k}{c}}
\;\stackrel{m\to1}{\sim}\;
\frac{\pi K(m)}{2K(1{-}m)}\cotanh\Bigl(\frac{\pi}{2K(1{-}m)}\wp^{-1}(V)\Bigr)
-\frac{\pi}{2K(1{-}m)}\wp^{-1}(V).
\label{ASAMAN}
\ee
Here the left-hand side is a finite (generally non-zero) constant, while the $K(m)$ in the first term on the right-hand side diverges when $m\to1$. In order for the left- and right-hand sides to have the same order, this divergence must be cancelled --- either by forcing the $\cotanh$ function to go to zero ($\wp^{-1}(V)$ is generally complex, so this is possible), or by having $\wp^{-1}(V)$ diverge as well so that the difference $\wp^{-1}-K\cotanh(\wp^{-1})$ remains finite. Which mechanism occurs depends on whether $k/c$ is below or above $-1/24$, and this, in turn, determines whether $V\to-2/3$ or $V\to+1/3$ in the limit:
\begin{itemize}
\setlength\itemsep{0em}
\item Suppose first that $k/c<-1/24$. Then $V<{-}\tfrac{m+1}{3}$ and $\wp^{-1}(V)$ is purely imaginary, ranging monotonically from $\wp^{-1}(-\infty)=0$ to $\wp^{-1}({-}\tfrac{m+1}{3})=iK(1{-}m)$. Since the latter value is finite when $m\to1$ (cfr eq.\ (\ref{KAPAFI})), $\wp^{-1}(V)$ cannot diverge in that region when $m\to1$. Thus the only way to cancel the divergence of $K(m)$ on the right-hand side of eq.\ (\ref{ASAMAN}) is to force the $\cotanh$ function to go to zero. This, in turn, only occurs when $\wp^{-1}(V)=iK(1{-}m)$, \ie when $V={-}\tfrac{m+1}{3}$. The latter converges to $V=-2/3$ in the limit $m\to1$.
\item Now suppose $k/c>-1/24$; in fact, let us even assume $k/c>0$ for definiteness. Then $V>\tfrac{2{-}m}{3}$ and $\wp^{-1}(V)$ is purely real, ranging monotonically from $\wp^{-1}(\tfrac{2{-}m}{3})=K(m)$ to $\wp^{-1}(+\infty)=0$. Since $K(m)$ diverges when $m\to1$ (cfr eq.\ (\ref{KAPADI})), $\wp^{-1}(V)$ also diverges near $V=\tfrac{2{-}m}{3}$, and this can be used to cancel the divergence of $K(m)$. Indeed, when $\wp^{-1}(V)$ is real and very large, the leading term of the right-hand side of (\ref{ASAMAN}) is $K-\wp^{-1}(V)$, which vanishes for $V=\tfrac{2{-}m}{3}$. The latter converges to $V=+1/3$ in the limit $m\to1$. A similar argument yields $V\to+1/3$ as $m\to1$ for level curves such that $-1/24<k/c<0$.
\end{itemize}
We have thus shown that all level curves of $k$ converge to $V=-2/3$ or $V=+1/3$ when $m\to1$, with the former limit occurring when $k/c<-1/24$ and the latter when $k/c>-1/24$. This confirms the pattern visible in fig.\ \ref{s71b}.\\

We stress that the mechanism leading to the limits $V\to\{{-}2/3,+1/3\}$ is radically different depending on whether $k/c<-1/24$ or $k/c>-1/24$. This difference translates into correspondingly different asymptotic `shapes' of level curves near the line $m=1$. Indeed, in appendix \ref{appC} we show that these level curves satisfy the following asymptotic behaviour (now refining the established result as $V\to\{{-}2/3,+1/3\}$):
\begin{itemize}
\setlength\itemsep{0em}
\item If $k/c<-1/24$, the level curve of $k$ converges to the point $(m=1,V=-2/3)$ and is tangent to the line $m=1$ in the $(m,V)$ plane. In fact, the level curve sticks to that line `with infinite velocity', as small changes in $V$ lead to non-perturbative changes in $m$ in that region:
\be
1-m
\:\stackrel{V\to-\tfrac{2}{3}}{\sim}\;
\frac{16}{e^2}
\exp\biggl[
-\pi\frac{\sqrt{|24k/c|}-1}{\sqrt{|V+2/3|}}
\biggr].
\label{NOPOTO}
\ee
Note how the difference $\sqrt{|24k/c|}-1$ controls the `strength' of the non-perturbative dependence of $1-m$ on $|V+2/3|$. For exceptional orbits (\ref{ss18}) having $n\geq2$, this difference reduces to $n-1$.
\item If $k/c>-1/24$, the level curve of $k$ converges to the point $(m=1,V=+1/3)$ with a finite, non-zero slope that depends on $k/c$ as follows:
\be
V
\;\stackrel{m\to1}{\sim}\;
\frac{1}{3}+
\biggl[
\cosh^2\Bigl(\sqrt{6\pi^2 k/c}\Bigr)-\frac{2}{3}
\biggr]
(1-m)+\cO\bigl((1{-}m)^2\bigr),
\label{POTO}
\ee
where the $\cosh$ becomes a $\cos$ when $-1/24\geq k/c\geq0$. The slope of $V$ as a function of $1-m$ thus increases with $k/c$, while it reaches its minimum when $k/c=-1/24$.
\end{itemize}
Both of these results are consistent with the patterns visible in fig.\ \ref{s71b}, and they can be confirmed by more detailed plots (which we omit). Note also, again from fig.\ \ref{s71b}, that the orbits in the forbidden wedge have level curves that interpolate between a perturbative behaviour at $(m=1,V=1/3)$ and a non-perturbative one at $(m=1,V=-2/3)$.

\subsection{Near \texorpdfstring{$\boldsymbol{m=0}$}{$m=0$}}
\label{sec444}

The last limiting region that remains to be studied in the $(m,V)$ plane is the one where $m\to0$. In fact, on the line $m=0$, any cnoidal profile (\ref{s21}) is just a constant, resulting in a simple linear dependence of $V$ on $k$:
\be
\lim_{m\to0}V
=\frac{24k}{c}+\frac{2}{3}
\quad\text{along a level curve.}
\label{VAKO}
\ee
This is apparent in fig.\ \ref{s60b} and it is consistent with the asymptotic relation (\ref{KARAPA}). A much less obvious question is how (\ref{VAKO}) gets modified when $m$ is small, but not quite zero: from fig.\ \ref{s71b} we expect a linear relation between $V$ and $m$ as $m\to0$, but the slope certainly depends on $k/c$. To evaluate this slope, we start once more from eq.\ (\ref{ss159}) and take the limit $m\to0$, where $K(m)\sim\tfrac{\pi}{2}(1+\tfrac{m}{4})$ is finite while\footnote{See eq.\ (\ref{KLAZ}) in appendix \ref{appC}.} \cite[sec.\ 23.12]{DLMF}
\be
\zeta(z)
\,\stackrel{m\to0}{\sim}\,
\frac{\pi^2}{4K(m)^2}\bigg(\frac{z}{3}+\frac{2K(m)}{\pi}\cot\Big(\frac{\pi z}{2K(m)}\Big)+\cO\big(m^2\big)\bigg).
\label{PAKATTO}
\ee
Using this and assuming $k/c>0$ for definiteness, the level curve (\ref{ss159}) simplifies into
\be
\sqrt{\frac{24k}{c}}
\sim
\cot\Big[\frac{\pi\wp^{-1}(V)}{2K(m)}\Big]+\cO(m^2)
\quad\Rightarrow\quad
V\sim\wp\bigg(\frac{2K(m)}{\pi}\arccot\bigg[\sqrt{\frac{24k}{c}}\bigg]\bigg)
+\cO(m^2).
\label{VAPAT}
\ee
Finally, using the asymptotic form of the Weierstrass $\wp$ function at large $K(1{-}m)$ that follows from the derivative of (\ref{PAKATTO}),\footnote{See eq.\ (\ref{KLA}) in appendix \ref{appC}.} we find that (\ref{VAPAT}) simplifies into
\be
V\sim\Big(1-\frac{m}{2}\Big)\bigg(\frac{24k}{c}+\frac{2}{3}\bigg)+\cO(m^2),
\qquad m\to0.
\label{VANT}
\ee
This linear behaviour is indeed consistent with the level curves visible in fig.\ \ref{s71b}, up to one apparent caveat: when $k=-c/24$, (\ref{VANT}) reduces to $V\sim-1/3+m/6$, which is {\it not} the correct behaviour $V=[(1\pm3)m-2]/6$. This seemingly suggests that the result (\ref{VANT}) is incorrect (at least for small $|24k/c-1|$), but the actual situation is more subtle: formula (\ref{VANT}) is {\it correct}, but the omitted coefficient of $m^2$ diverges when $k=-c/24$ (and only then).\footnote{In fact, that coefficient is proportional to $|V+1/3|$ to a negative power.} As a result, (\ref{VANT}) really holds for all values of $k$ except $k=-c/24$, where the expansion of $V$ in powers of $m$ becomes ill-defined. This divergence of the second-order coefficient can, again, be confirmed by detailed plots which we omit here.

\section{Wave shoaling}
\label{sec4}

So far our analysis has been entirely mathematical: we investigated the coadjoint orbits of cnoidal waves from a group-theoretic viewpoint, irrespective of any phenomenological realization. It is tempting, however, to think of fig.\ \ref{s71b} as a bifurcation diagram for cnoidal waves, where slow variations of $(m,V)$ would lead to visible changes in the waves' behaviour. This is especially true of transitions into the `forbidden wedge', since they correspond to dramatic changes in the symmetry properties of Virasoro orbits: outside the wedge, there always exists a frame in which the wave profile is conformally equivalent to a constant $k$, whereas such a frame {\it never} exists inside the wedge. It is natural to wonder if such abrupt changes have correspondingly abrupt observable consequences.\\

\begin{figure}[t]
\centering
\includegraphics[width=.6\textwidth]{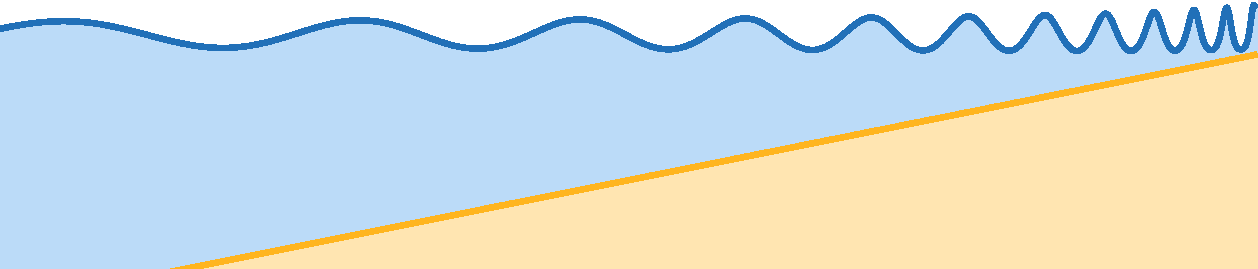}
\caption{A shoal linking a deep sea region on the left to a shallow region on the right. As waves propagate on the free water surface from left to right, they get deformed --- higher crests, shorter wavelengths. This behaviour is referred to as {\it shoaling}.\label{s23}}
\end{figure}

As an attempt to answer this question, we now investigate wave shoaling in shallow water --- a phenomenon extensively studied in fluid dynamics and coastal engineering \cite{Dean}. Shoals, or sandbanks, are near-coastal regions where the seabed rises from some large depth to nearly water level (see fig.\ \ref{s23}). Waves propagating on the water/air interface then undergo a series of deformations: they get refracted and reflected, change shape and velocity, and lose energy through friction. If the waves keep a cnoidal profile throughout this process, then they effectively trace a path in the bifurcation diagram of fig.\ \ref{s71b}. Our task is to find that path. In fact, under suitable assumptions, the solution is elementary: the path will consist of profiles with zero average, given by eq.\ (\ref{ss31}) below. The latter formula only involves the parameters $m$ and $V$, without reference to the wavelength $\lambda$ or the average depth $h$. Its most striking aspect is that it crosses the (upper) boundary of the bifurcation wedge of fig.\ \ref{s71b}, at a critical value of $m$ that we compute below eq.\ (\ref{s37b}). By relating $m$ to the wavelength $\lambda$ and the depth $h$, one obtains their corresponding critical values and concludes that shoaling allows cnoidal waves to cross into the forbidden wedge when the fluid becomes sufficiently shallow. Tempting as it is to interpret this crossing as an observable effect --- related for instance to wave breaking \cite{Dean,Svendsen02,Brun} ---, we will unfortunately find no such implication of our analysis. (It is conceivable that a different conclusion would hold upon modelling shallow water dynamics through, say, the Camassa-Holm equation \cite{Camassa} that notoriously supports breaking waves which KdV lacks \cite{Constantin98}, but this modification goes beyond the scope of the present work; we shall return to it briefly in section \ref{sec5}.)\\

This section is organized as follows. We first briefly review the KdV description of shallow water dynamics, so as to make contact between the abstract concepts of section \ref{sec2} and observable fluid-mechanical effects. Then we turn to wave shoaling applied to a train of cnoidal waves propagating in a fluid with a gently sloping bottom, following the derivation of \cite{Svendsen}. Finally, we relate shoaling to a path of profiles with zero average in the bifurcation diagram of fig.\ \ref{s71b}.

\subsection{The KdV equation in shallow water dynamics}
\label{sec41}

We consider a two-dimensional perfect fluid with uniform density $\rho$ embedded in a constant gravitational field $\bg$. We let $X,Y$ be spatial (laboratory) coordinates, with $Y$ pointing along the vertical axis so that $\bg=-g{\bf 1}_{Y}$. Assuming the domain where the fluid lives to be simply connected, and provided vorticity vanishes, the fluid's velocity field can be written as a gradient $\bu=\nabla\phi$. The velocity potential then satisfies
\be
\nabla^2\phi=0
\qquad\text{and}\qquad
\frac{\der\phi}{\der t}+gY+\frac{p}{\rho}+\demi(\nabla\phi)^2=0
\label{s25}
\ee
at any point $(X,Y)$ and any time $t$. We denote by $p$ the pressure inside the fluid; for shallow water, it is essentially constant throughout the fluid, so we may neglect it by redefining $\phi\mapsto\phi-pt/\rho$ (which does not affect the velocity $\nabla\phi$). Finally, we assume that the fluid has a free upper surface with a depth profile $D(X,t)$ (unknown at this stage) and a fixed lower boundary at $Y=0$. This entails the boundary conditions
\be
\left.\frac{\der\phi}{\der Y}\right|_{Y=0}=0
\qquad\text{and}\qquad
\left.\frac{\der\phi}{\der Y}\right|_{Y=D}
=
\frac{\der D}{\der t}+\left.\frac{\der\phi}{\der X}\right|_{Y=D}\frac{\der D}{\der X},
\label{s26}
\ee
while the second equation in (\ref{s25}), evaluated at $Y=D(X,t)$, becomes
\be
\left.\frac{\der\phi}{\der t}\right|_{Y=D}+gD+\left.\demi(\nabla\phi)^2\right|_{Y=D}=0.
\label{ss26}
\ee
Along with the conservation equation $\nabla^2\phi=0$ (which holds throughout the fluid), the three boundary conditions (\ref{s26})-(\ref{ss26}) are the tools used to find the unknown functions $\phi(X,Y,t)$ and $D(X,t)$.\\

KdV dynamics emerges in the shallow water limit of eqs.\ (\ref{s25})-(\ref{s26})-(\ref{ss26}), when describing the (slow) time evolution of one wave. The limit is designed as follows. First, let us assume to be describing a train of identical waves with wavelength $\lambda$. We also assume that the average depth of each wave is a fixed parameter $h$, given by
\be
\frac{1}{\lambda}\int_0^{\lambda}\dd X\,D(X,t)
=
h
\qquad\forall\, t.
\label{s27}
\ee
Using the two dimensionful parameters $(\lambda,h)$ that describe the setup, we define dimensionless spatial coordinates $(X',Y')$ and a dimensionless time $t'$ by
\be
X'\equiv\frac{X}{\lambda},
\qquad
Y'\equiv\frac{Y}{h},
\qquad
t'\equiv\frac{\sqrt{gh}}{\lambda}t.
\nn
\ee
Now, a key assumption for shallow water dynamics is that the average depth is much smaller than the wavelength, \ie $h^2/\lambda^2\equiv\epsilon\ll 1$. Starting from this fact, one solves eqs.\ (\ref{s25})-(\ref{s26})-(\ref{ss26}) perturbatively, as power series in $\epsilon$. The depth field is then written as
\be
D(X',t')
\sim
h\Bigl[1+\epsilon\,\frac{p(X',t')}{2\pi}+\cO(\epsilon^2)\Bigr]
\qquad\text{in the limit }\epsilon\to0.
\label{t28}
\ee
At lowest order, the height deviation $p$ satisfies a (linear) wave equation, so $p$ is in general a superposition of left- and right-moving perturbations. To derive KdV as a slow-time evolution equation for the wave profile, we pick one of the two possibilities --- say a right-moving wave $p(X'-t')$ --- and take care of subleading corrections thanks to a (dimensionless) slow time coordinate $\tau\equiv\epsilon\,t'/2$. In terms of the comoving coordinate $x\equiv2\pi(X'-t')$, then, the first subleading correction to the wave equation yields the KdV equation (\ref{s5}) for $p$ with a negative central charge $c=-32\pi^3$:
\be
\frac{\der p}{\der\tau}+3p\frac{\der p}{\der x}+\frac{(2\pi)^3}{3}\frac{\der^3p}{\der x^3}=0.
\label{ss29}
\ee
Recall that we assumed, at the outset, that $p$ describes a wave with fixed wavelength $\lambda$; in terms of the dimensionless comoving coordinate $x$, this makes $p$ $2\pi$-periodic as in eq.\ (\ref{s6}). As for the constraint (\ref{s27}), it means that the profile $p$ has vanishing average:
\be
\int_0^{2\pi}\dd x\,p(x)=0.
\label{s29}
\ee

\paragraph{Cnoidal waves.} In terms of dimensionless quantities, a cnoidal wave solution of (\ref{ss29}) with wavelength $2\pi$ is given by eq.\ (\ref{ss6}) with $c=-32\pi^3$. In principle it contains two free dimensionless parameters ($m$ and $V$), but the zero-average condition (\ref{s29}) forces them to be related. Indeed, using the identity\footnote{See eq.\ (\ref{NADA}) in appendix \ref{appB}.}
\be
\label{EMINEM}
\int_0^{2K(m)}\dd x\,\dn^2(x|m)
=
2E(m)
\ee
where $E(m)$ is the complete elliptic integral of the second kind, eq.\ (\ref{s29}) yields
\be
\boxed{%
V
=
\frac{2E(m)}{K(m)}-\frac{4-2m}{3}}
\qquad
\text{(zero-average curve).}
\label{ss31}
\ee
Accordingly, a cnoidal wave (\ref{CNOMEGA}) with vanishing average reads
\be
p(x,\tau)
=
-\frac{c}{3\pi^2}K(m)^2
\biggl[
\dn^2\biggl(\frac{K(m)}{\pi}\Bigl(x+32\pi K(m)^2\bigl[\tfrac{E(m)}{K(m)}-\tfrac{2-m}{3}\bigr]\tau\Bigr)\bigg|m\biggr)
-\frac{E(m)}{K(m)}
\biggr].
\label{s31}
\ee
In terms of dimensionful variables, we need to plug (\ref{s31}) into eq.\ (\ref{t28}) for the depth of fluid at $X,t$. The result, now also including $c=-32\pi^3$, is
\be
D(X,t)
\sim
h+\frac{16}{3}\frac{h^3}{\lambda^2}K(m)^2
\biggl[
\dn^2\biggl(
\frac{2K(m)}{\lambda}
\big(X-\sqrt{gh}\,t\big)
\bigg|m\biggr)
-\frac{E(m)}{K(m)}
+\cO(\epsilon)
\biggr],
\label{s32}
\ee
which holds in the limit $h^2/\lambda^2=\epsilon\to0$. This is how a train of cnoidal waves looks like as seen from a laboratory with a static spatial coordinate $X$ and time $t$. It is uniquely determined by one dimensionless parameter $m$ and two dimensionful parameters $h,\lambda$ --- respectively the pointedness, the average depth and the wavelength.

\subsection{Wave shoaling}
\label{sec43}

We consider the standard setup \cite{Svendsen,Dean} of a train of cnoidal waves (\ref{s32}) incoming from the left and encountering a shoal, that is, a region where the seabed gently slopes towards a beach. We assume that the slope is so weak that (i) at any point, the actual wave profile is accurately described by a train of cnoidal waves with a horizontal seabed, and (ii) there is no reflection of waves. As a result, the wave profile at any point takes the form (\ref{s32}), except that the parameters $(h,\lambda,m)$ are $X$-dependent. Since the $X$-dependence of $h$ is assumed to be known --- it is a monotonously increasing function with very weak slope ---, the game is to find relations between $h$ and the remaining parameters $\lambda,m$ so as to predict the wave's shape throughout the shoal.\\

\begin{figure}[t]
\centering
\includegraphics[width=0.50\textwidth]{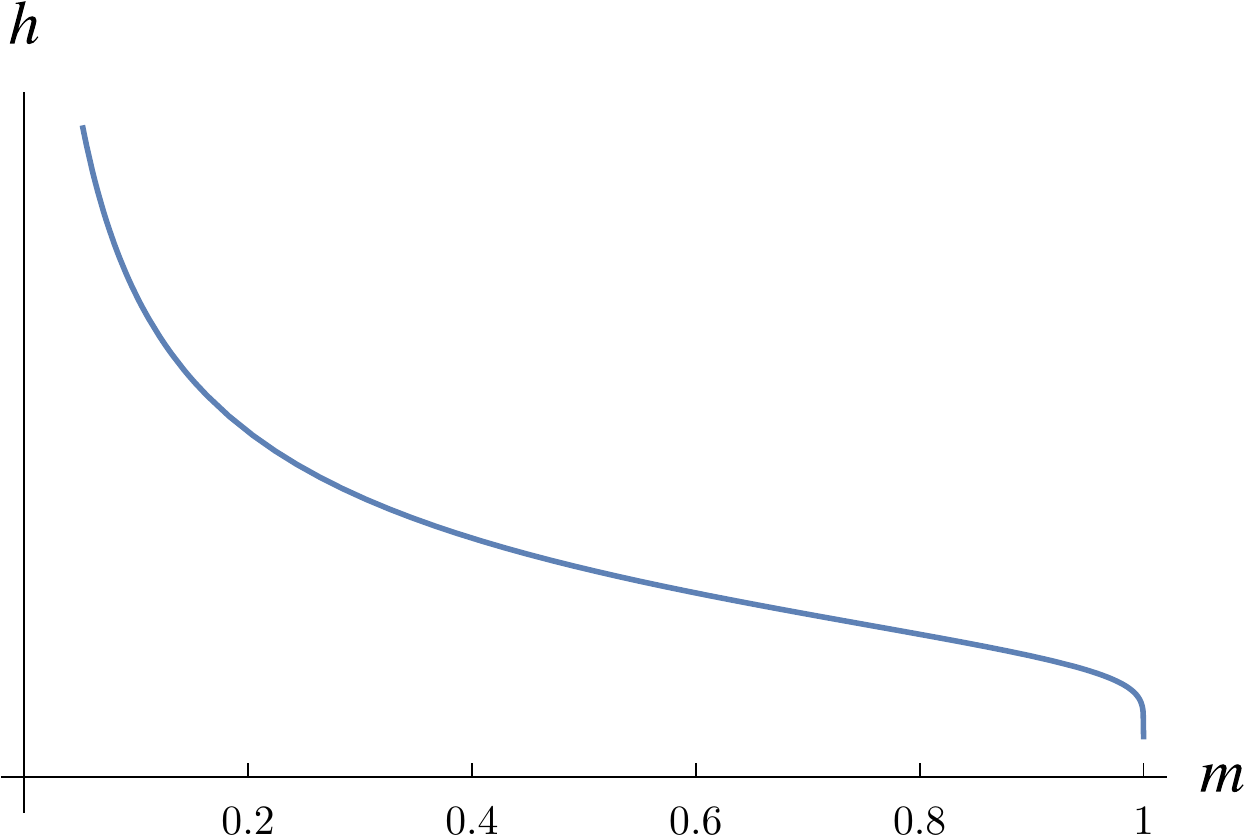}
\caption{The curve (\ref{s36}). As anticipated in fig.\ \ref{s23}, $m$ is small when $h$ is large, and vice-versa. The curve is injective, so eq.\ (\ref{s36}) can be inverted to express $m$ in terms of $h$.\label{FIGs36}}
\end{figure}

Since our purpose here is merely to illustrate Virasoro orbit transitions as they occur in shoaling (as opposed to offering a new derivation of shoaling {\it per se}), we follow the argument of \cite{Svendsen}. In that approach, two parameters are assumed to be kept fixed throughout the shoaling process:
\begin{itemize}
\item[(i)] The period $T$ of the wave (at any position $X$). Owing to the explicit formula (\ref{s32}) and the fact that $\dn^2(\cdot|m)$ has period $2K(m)$, one has
\be
T
=
\frac{\lambda}{\sqrt{gh}}[1+\cO(\epsilon)]
\qquad\Rightarrow\qquad
\lambda\sim\sqrt{gh}\,T.
\label{s34}
\ee
Thus, to leading order in $\epsilon=h^2/\lambda^2$ and at fixed parameters $g,T$, the wavelength is proportional to $\sqrt{h}$. Shallower waters yield shorter wavelengths.
\item[(ii)] The energy transport $\cF$ of the wave train, \ie the integral, over one period $T$, of the wave's energy density multiplied by its velocity\footnote{Thus $\cF$ has dimensions of force, \ie mass $\times$ length $\times$ time$^{-2}$.} \cite{Svendsen}:
\be
\cF
\equiv
\rho g\sqrt{gh}
\int_0^T\dd t\,\frac{h^6}{\lambda^4}p^2(X,t).
\label{s35}
\ee
Using eq.\ (\ref{s31}) for $p$ and the identity (\ref{EMINEM}) together with\footnote{See eq.\ (\ref{TADAM}) in appendix \ref{appB}.}
\be
\int_0^{2K(m)}\dd x\,\dn^4(x|m)
=
\frac{2m-2}{3}K(m)+\frac{8-4m}{3}E(m),
\label{MYSTERY}
\ee
the energy transport (\ref{s35}) becomes
\be
\cF
=
\frac{256}{9}\rho g\frac{h^6}{\lambda^3}
\Bigl[
\frac{m-1}{3}K(m)^4+\frac{4-2m}{3}E(m)K(m)^3-E(m)^2K(m)^2
\Bigr]
\bigl(1+\cO(\epsilon)\bigr).
\nn
\ee
The requirement $\cF=$ cst is the second conservation equation that can be used to describe wave shoaling. It yields a relation between $h$, $\lambda$ and $m$. Using eq.\  (\ref{s34}) to express $\lambda$ as a function of $h$, one obtains the following relation between $h$ and $m$:
\be
h^{9/2}
=
\frac{27}{256}
\frac{\sqrt{g}}{\rho}
T^3\cF
\Bigl[
(m-1)K(m)^4+2(2-m)E(m)K(m)^3-3E(m)^2K(m)^2
\Bigr]^{-1}.
\label{s36}
\ee
The resulting graph $h(m)$ is plotted in fig.\ \ref{FIGs36}. As expected, smaller $h$ corresponds to higher $m$, and vice-versa. Note that the dependence of $h$ on the combination $\sqrt{g}T^3\cF/\rho$ could not have been guessed on the basis of dimensional analysis alone!
\end{itemize}
Eqs.\ (\ref{s34})-(\ref{s36}) are two constraints on the three parameters $(\lambda,h,m)$ that determine the profile (\ref{s32}). Thus, shoaling results in a path traced by the wave profile in $(\lambda,h,m)$ space. The average depth can be taken as the parameter along the path, with $\lambda$ and $m$ expressed as functions of $h$ through eqs.\ (\ref{s34}) and (\ref{s36}). These formulas ensure that, when $h$ decreases slowly, $\lambda$ decreases while $m$ increases, as shown in fig.\ \ref{s23}.

\subsection{Shoaling and Virasoro orbits}
\label{sec44}

We now relate wave shoaling to a path in the bifurcation diagram of fig.\ \ref{s71b}, \ie we think of shoaling as a curve in the $(m,V)$ plane of dimensionless cnoidal parameters. Given the constraint (\ref{s29}) that the average of the profile vanishes, the curve is simply given by eq.\ (\ref{ss31}). It is utterly insensitive to the behaviour of $h$ and $\lambda$ determined by eqs.\ (\ref{s34})-(\ref{s36}).\\

The zero-average curve (\ref{ss31}) is depicted in fig.\ \ref{s71b}. It has two striking features: first, at small $m$ it is asymptotic to the straight line $V=\tfrac{2{-}m}{3}$ that separates elliptic and hyperbolic orbits with zero winding. Secondly, as $m$ increases, the curve first deviates towards the inside of the region consisting of elliptic orbits with zero winding, but then reaches a critical value $m^*$ where it crosses the bifurcation line $V=\tfrac{2m{-}1}{3}$. For $m$ larger than $m^*$, the curve lies in the `forbidden wedge' of hyperbolic orbits with unit winding. Analytically, the critical value of $m$ is obtained by matching eq.\ (\ref{ss31}) with eq.\ (\ref{KOMEGGA}) for the boundaries of the wedge (with a plus sign in (\ref{KOMEGGA})):
\be
V
=
\frac{2E(m^*)}{K(m^*)}-\frac{4-2m^*}{3}
\stackrel{!}{=}
\frac{2m^*-1}{3}
\qquad\Rightarrow\qquad
\frac{E(m^*)}{K(m^*)}
=
\demi.
\label{s37b}
\ee
The solution can be found numerically\footnote{To do this, iterate the map $m\mapsto(2E(m)/K(m))-1+m$, which has a unique stable fixed point.} to be $m^*=0.8261147659849698...$. Note that there is no value of $m\in[0,1)$ for which the curve (\ref{ss31}) crosses the lower boundary of the forbidden wedge, since the solution of $(2E/K)-(4-2m)/3={-}\tfrac{m+1}{3}$ is $m=1$.\\

Eq.\ (\ref{s36}) says that a wave train propagating into shallower waters (smaller $h$) gradually becomes more steeply peaked (larger $m$). Accordingly, waves incoming from deep water and moving towards a beach follow the zero-average curve of fig.\ \ref{s71b} from left to right. If the water becomes shallow enough, there is a point at which $m$ crosses the critical value (\ref{s37b}) where the path enters the `forbidden wedge'. At that point, the average depth is
\be
h^*
=
\biggl(\frac{gT^6\cF^2}{\rho^2}\biggr)^{1/9}
\biggl(\frac{3}{4K(m^*)^{4/3}}\biggr)^{2/3}
\simeq
0.3905
\biggl(\frac{gT^6\cF^2}{\rho^2}\biggr)^{1/9}
\nn
\ee
where we used (\ref{s37b}) to simplify the right-hand side of (\ref{s36}). Thus, shoaling definitely produces transitions between different Virasoro orbits; some of these transitions are sharp bifurcations between orbits of radically different nature --- for instance orbits that have a rest frame and orbits that do not. Once a wave enters in the forbidden wedge, there exists no one-dimensional diffeomorphism (conformal transformation) that turns it into a uniform profile. Of course, the fluid is in fact two-dimensional and there exists a {\it two-dimensional} diffeomorphism that `flattens' the profile; the point is that, in the forbidden wedge, this transformation {\it must} shuffle the fluid's particles along a vertical direction, since a purely horizontal diffeomorphism would not be able to flatten the wave.\\

Unfortunately, at this stage it is unclear whether these abrupt geometric bifurcation have observable consequences. A natural candidate for such phenomena would be wave breaking \cite{Dean,Svendsen02}, \eg the convective breaking \cite{Brun} that occurs when the velocity of fluid particles at the top of a wave's crest becomes larger than the wave's velocity. One might hope that the transition  (\ref{s37b}) into the forbidden wedge somehow reflects that phenomenon. However, this is not the case for an elementary reason: the breaking point of a wave is an intricate function of its parameters $(h,\lambda,m)$ \cite{Brun}, and is not given by as simple a criterion as eq.\ (\ref{s37b}). Nevertheless, it is still conceivable that (\ref{s37b}) does describe an observable effect; for instance, we will show in \cite{OblakKozy} that it corresponds to a bifurcation that dramatically affects the Euler-Poincar\'e reconstruction \cite{Holm} of the KdV equation. Whether this, or similar effects, can be seen in a lab, is left for further study.

\section{Conclusion and outlook}
\label{sec5}

The purpose of this paper has been to initiate a group-theoretic study, based on Virasoro symmetry, of shallow water dynamics as described by the KdV equation. We addressed the problem of finding Virasoro orbits of cnoidal waves, asking in particular which of these waves can be turned into uniform profiles thanks to diffeomorphisms. The result of this investigation is summarized by fig.\ \ref{s71b}, whose structure mimics the sequence of orbits that appear as $V$ increases at fixed $m$:
\begin{table}[h]
\centering
\begin{tabular}{|c|l|}
\hline
Range of $V$ & Type of Virasoro orbit\Tstrut\Bstrut\\
\hline
$({-}\infty,{-}\tfrac{m+1}{3})$ & Elliptic with winding $n\geq1$\Tstrut\\[.1cm]
$\{{-}\tfrac{m{+}1}{3}\}$ & Exceptional orbit $n=1$\\[.1cm]
$({-}\tfrac{m{+}1}{3},\tfrac{2m{-}1}{3})$ & Hyperbolic with winding $n=1$ (`forbidden wedge')\\[.1cm]
$\{\tfrac{2m{-}1}{3}\}$ & Exceptional orbit $n=1$\\[.1cm]
$(\tfrac{2m{-}1}{3},\tfrac{2{-}m}{3})$ & Elliptic with winding $n=0$\\[.1cm]
$\{\tfrac{2{-}m}{3}\}$ & Parabolic with winding $n=0$\\[.1cm]
$(\tfrac{2{-}m}{3},+\infty)$ & Hyperbolic with winding $n=0$\Bstrut\\
\hline
\end{tabular}
\end{table}

\noindent Except for the wedge ${-}\tfrac{m+1}{3}<V<\tfrac{2m{-}1}{3}$, all these orbits have uniform representatives given by eq.\ (\ref{ss159}). This classification is closely related to the band structure of the $N=1$ Lam\'e equation \cite{Ince,Ince2,Whittaker,Arscott}, and was thus implicitly known thanks to the inverse scattering approach to the KdV equation \cite{Gardner,Novikov,NovikovBook}. To our knowledge, however, an analysis of the resulting classification that is as detailed as the one of this paper, including the asymptotics of section \ref{sec36}, has never appeared in the literature.\\

We have also seen how wave shoaling produces transitions between different Virasoro orbits by following the curve (\ref{ss31}) in the $(m,V)$ plane. As we saw, this curve enters in the `forbidden wedge' at the critical value of $m$ given by (\ref{s37b}). Beyond that critical value, there exists no reference frame in which the wave profile is uniform --- at least when restricting oneself to one-dimensional diffeomorphisms. It is unclear if this bifurcation entails any observable effect, but a first hint that the answer might be positive will be provided in the follow-up \cite{OblakKozy}, where we will show how the bifurcation diagram of fig.\ \ref{s71b} affects certain Berry phases carried by KdV solitons. In particular, the bifurcation lines at the boundaries of the forbidden wedge then correspond to noticeable, abrupt changes in the late-time behaviour of Euler-Poincar\'e reconstruction \cite{Holm} --- a close analogue of the equation of motion for fluid particles.\\

From a broader perspective, the general theme of this work is relevant not only for KdV, but also for other nonlinear wave systems in $1+1$ dimensions. For instance, the Camassa-Holm equation \cite{Camassa} models shallow water dynamics with slightly different scalings than the ones giving rise to KdV, but it is also (crucially!) an Euler-Arnold equation for the Virasoro group. This makes it amenable to geometric arguments essentially identical to those of this paper, up to the fact that the solitons to be considered would not be cnoidal waves, but `coshoidal waves' (\ie periodic `peakons') or their smooth version \cite{Boyd}. Beyond solitons, a further intriguing aspect of this system is the presence of wave-breaking profiles: in contrast to KdV, not all solutions of Camassa-Holm are smooth for arbitrarily late times --- they develop singularities \cite{Constantin98}. Any such blow-up takes place, by construction, on a single coadjoint orbit of the Virasoro group, and it would be interesting to see what aspect of Virasoro geometry (if any) is sensitive to it.
Similar questions can be raised more generally for any Virasoro-based Euler-Arnold equation, notoriously including the Hunter-Saxton system \cite{Hunter} that describes waves in liquid crystals. We leave such considerations for future work.\\

Another plausible application of the classification presented here occurs in three-dimensional gravity. Indeed, it was shown in \cite{Perez:2016vqo} that suitable fall-off conditions on the dreibein and spin connection of Anti-de Sitter gravity lead to KdV dynamics on the space-time boundary. In that context, the Virasoro orbit of a KdV soliton is interpreted as a family of gravitational backgrounds. For example, from that perspective, the forbidden wedge of fig.\ \ref{s71b} consists of tachyons, while its boundaries are cnoidal profiles that all belong to the orbit of the pure AdS$_3$ background under Brown-Henneaux diffeomorphisms \cite{Brown:1986nw}. Similarly, cnoidal waves in the region $\tfrac{2m{-}1}{3}<V<\tfrac{2{-}m}{3}$ belong to orbits of conical deficits (massive particles in AdS$_3$), while waves for which $V>\tfrac{2{-}m}{3}$ belong to orbits of BTZ black holes \cite{Banados:1992wn}. It would be interesting to see whether coupling gravity to some bulk degrees of freedom can lead to adiabatic motion that mimics wave shoaling, producing transitions between `massive' configurations that have a rest frame, and `tachyonic' ones that do not.

\section*{Acknowledgements}

I am grateful to G.\ Kozyreff for collaboration on the related project \cite{OblakKozy}. This work is supported by the Swiss National Science Foundation, and also partly by the NCCR SwissMAP.

\appendix

\section{The Virasoro group}
\label{appA}

In this appendix we review basic definitions regarding the Virasoro group, its algebra, its coadjoint representation, and the notion of densities on the circle. We refer \eg to \cite{Guieu} or \cite[chap.\ 6]{Oblak:2016eij} for much more thorough presentations.

\paragraph{Virasoro group.} The Virasoro group is the universal central extension of $\Diff$, so we first describe the latter. Given a coordinate $x\in\RR$, a (lift of a) diffeomorphism of the circle is a smooth map $f:\RR\rightarrow\RR:x\mapsto f(x)$ such that
\be
f'(x)>0,
\qquad
f(x+2\pi)
=
f(x)+2\pi.
\label{ta1}
\ee
The set $\Diff$ of such maps forms a group under composition. To introduce the Virasoro group, define the {\it Bott cocycle} $\sfC(f,g)$, for all $f,g\in\Diff$, as
\be
\sfC(f,g)
\equiv
-\frac{1}{48\pi}\int_0^{2\pi}\!\dd x\,\log(f'\circ g)\frac{g''}{g'}.
\label{sa1}
\ee
The {\it Virasoro group} then consists of pairs $(f,\alpha)\in\Diff\times\RR$, with multiplication
\be
(f,\alpha)\cdot(g,\beta)
=
\bigl(f\circ g,\alpha+\beta+\sfC(f,g)\bigr).
\label{ssa1}
\ee
One can verify that the definition of the Bott cocycle (\ref{sa1}) makes this product associative. The identity is $(\II,0)$ with $\II(x)=x$, and the inverse of $(f,\alpha)$ is $(f,\alpha)^{-1}=(f^{-1},-\alpha)$, where $f^{-1}$ is the (unique) diffeomorphism such that $f(f^{-1}(x))=f^{-1}(f(x))=x$.

\paragraph{Adjoint representation and Virasoro algebra.} The Lie algebra of $\Diff$ consists of infinitesimal diffeomorphisms, that is, vector fields $\xi(x)\der_x$ on the circle. Accordingly, the algebra of the Virasoro group (which extends $\Diff$ by $\RR$) consists of pairs $(\xi,\alpha)\in\Vect\oplus\RR$. As for any Lie group, the adjoint representation is defined by
\be
\Ad_{(f,\alpha)}(\xi,\beta)
\equiv
\left.\frac{\dd}{\dd t}\right|_{t=0}\Bigl[(f,\alpha)\cdot(e^{t\xi},t\beta)\cdot(f,\alpha)^{-1}\Bigr]
\nn
\ee
where $e^{t\xi}$ is the exponential of $t\xi$; it is the flow of the vector field $\xi$ at time $t$, so $e^{t\xi}(x)=x+t\xi(x)+\cO(t^2)$. Using the group operation (\ref{ssa1}), one finds the explicit formula
\be
\Ad_{(f,\alpha)}(\xi,\beta)
=
\Bigl(\Ad_f\xi,\beta-\frac{1}{24\pi}\int_0^{2\pi}\sfS[f](x)\xi(x)\Bigr)
\label{sa2}
\ee
where $\sfS$ is the Schwarzian derivative (\ref{shabadoo}) and $\Ad_f\xi$ denotes the standard transformation law of vector fields under diffeomorphisms:
\be
\bigl(\Ad_f\xi\bigr)(f(x))
=
f'(x)\xi(x).
\label{ssa2}
\ee
By differentiating the adjoint representation (\ref{sa2}) with respect to its argument $f$, one can show that the Lie bracket reproduces the standard Virasoro algebra of two-dimensional CFTs. See \cite[sec.\ 6.4]{Oblak:2016eij} for details.

\paragraph{Coadjoint representation.} The space of KdV wave profiles is, by definition, the (smooth) dual space of the Virasoro algebra. It consists of pairs $(p,c)$ where $p(x)$ is a $2\pi$-periodic function\footnote{Technically, $p=p(x)\dd x^2$ is a quadratic density rather than a function.} while $c\in\RR$ is a number --- a {\it central charge}. Each such dual element $(p,c)$ pairs with the Virasoro algebra as
\be
\left<
(p,c),(\xi,\alpha)
\right>
\equiv\frac{1}{2\pi}\int_0^{2\pi}\!\dd x\,p(x)\xi(x)+c\alpha.
\nn
\ee
As for any Lie group, the {\it coadjoint representation} of the Virasoro group is defined by
\be
\left<\Ad^*_{(f,\alpha)}(p,c),(\xi,\beta)\right>
\equiv
\left<(p,c),\Ad_{(f,\alpha)^{-1}}(\xi,\beta)\right>.
\nn
\ee
Using eq.\ (\ref{sa2}) along with suitable integrations by parts and changes of variables, a cumbersome but straightforward calculation yields $\Ad^*_{(f,\alpha)}(p,c)=(f\cdot p,c)$, where $f\cdot p$ is the transformation law of a (chiral) CFT stress tensor $p$ under a conformal transformation $f$:
\be
(f\cdot p)(f(x))
=
\frac{1}{(f'(x))^2}\Bigl[p(x)+\frac{c}{12}\sfS[f](x)\Bigr].
\label{sa3}
\ee
This is the formula that we used in (\ref{ss8}). Furthermore, one can define an infinitesimal coadjoint representation by differentiating (\ref{sa3}): $\delta_{\xi}p\equiv\left.\der_t\right|_{t=0}
\bigl(e^{t\xi}\cdot p\bigr)$. The result is
\be
\delta_{\xi}p
=
-\xi p'-2p'\xi+\frac{c}{12}\xi'''
\label{ssa3}
\ee
where the third derivative $\xi'''$ is a remnant of $f'''$ in the Schwarzian derivative (\ref{shabadoo}). Upon setting $\xi=p$, the right-hand side of (\ref{ssa3}) coincides with that of the KdV equation (\ref{s7}).

\paragraph{Densities.} Spaces of densities form important representations of the Virasoro group. By definition, a {\it density with weight $h\in\RR$} is a function $\psi(x)$ (not necessarily $2\pi$-periodic) that transforms under $\Diff$ as
\be
(f\cdot\psi)(f(x))
\equiv
(f'(x))^{-h}\psi(x).
\label{sa4}
\ee
In CFT language, one would say that $\psi(x)$ is a (chiral) primary field with weight $h$. For example, eq.\ (\ref{ssa2}) says that vector fields are densities of weight $h=-1$; and the coadjoint representations (\ref{sa3}) says that wave profiles would be densities with weight $h=+2$ if it weren't for the central charge and the Schwarzian derivative. Finally, the solutions of Hill's equation (\ref{s11}) are densities of weight $-1/2$.

\section{Weierstrass and Jacobi elliptic functions}
\label{appB}

In this appendix we review elementary properties of the elliptic functions that play a key role for cnoidal waves and their orbits. We start by introducing Weierstrass elliptic functions and derive some of the important identities used in section \ref{sec3}. Then we define Jacobi elliptic functions as generalizations of trigonometric functions, and show that their squares coincide (up to additive and multiplicative constants) with the Weierstrass $\wp$ function. The standard reference on these matters is the book \cite{Whittaker}, especially chapters {\sc xx} and {\sc xxii}. Certain recent references also contain accessible reviews on this subject --- see \eg \cite{Pastras:2017wot} for Weierstrass functions and \cite{Singh} for Jacobi functions.

\subsection{Weierstrass elliptic functions}
\label{appBsec1}

In section \ref{sec3}, the Weierstrass elliptic function $\wp(z)$ (and its cousins $\zeta$ and $\sigma$) played a key role. Accordingly, we now introduce these functions and review some of their properties. For more details we refer to \cite[chap.\ {\sc xx}]{Whittaker} or \cite{Pastras:2017wot}.

\paragraph{Elliptic functions.} An {\it elliptic function} is a complex function $F$, defined on some subset of the complex plane, which is meromorphic and doubly periodic. The latter condition means that there exist non-zero complex numbers $\omega_1$, $\omega_2$ such that $\text{Im}(\omega_2/\omega_1)\neq0$ and
\be
F(z+2n\omega_1+2n'\omega_2)
=
F(z)
\qquad\forall\,n,n'\in\ZZ
\nn
\ee
for all $z\in\CC$ that belong to the domain of $F$. The numbers $\omega_1$ and $\omega_2$ are called the {\it half-periods} of $F$, and the periods $\{2\omega_1,2\omega_2\}$ generate a {\it lattice}
\be
\Gamma
\equiv
\big\{
2n\omega_1+2n'\omega_2
\big|
n,n'\in\ZZ
\big\}
\subset
\CC.
\label{bb37}
\ee
A {\it unit cell} of that lattice is a parallelogram with corners $\{z,z+2\omega_1,z+2\omega_1+2\omega_2,z+2\omega_2\}$ for any fixed $z\in\CC$.  Two numbers $z,w\in\CC$ are said to be {\it congruent} if they differ by a point in $\Gamma$, \ie if there exist integers $n,n'$ such that $w=z+2n\omega_1+2n'\omega_2$.\\

The interplay between meromorphicity and double periodicity is responsible for the richness of elliptic functions. For instance, given an elliptic function $F(z)$, the following key properties follow from the residue theorem applied to the integral of $F(z)$ along the boundary of a unit cell:
\begin{enumerate}
\setlength\itemsep{0em}
\item[(i)]\label{MastaVanish} The sum of residues of $F$ over its poles in a unit cell vanishes.
\item[(ii)]\label{MastaPole} If $F$ has no poles, then it is constant. The number of poles of $F$ in a unit cell, weighted by their order, is called the {\it order} of $F$.
\item[(iii)] There exists no elliptic function of order one. (This is a corollary of (i) and (ii).)
\end{enumerate}
Further important properties can be derived as follows. Let $H(z)$ be an analytic function, and let $G(z)$ be a meromorphic function with $N$ poles $u_i$ and $N'$ roots $v_j$ (of respectives orders $\mu_i$, $\nu_j$) enclosed by a counterclockwise contour $C$. Then the residue theorem yields
\be
\frac{1}{2\pi i}\oint_C\dd z\,H(z)\frac{G'(z)}{G(z)}
=
\sum_{i=1}^N\mu_i\,G(u_i)
-
\sum_{j=1}^{N'}\nu_j\,G(v_j).
\label{b38}
\ee
Thus, if $F(z)$ is an elliptic function, one finds that:
\begin{enumerate}
\setlength\itemsep{0em}
\item[(iv)] For any $z_0\in\CC$, the number of roots (weighted by their order) of the equation $F(z)=z_0$ in a unit cell equals the order of $F$. (The proof follows from (\ref{b38}) with $C$ the boundary of a unit cell, $H(z)=1$ and $G(z)=F(z)-z_0$.)
\item[(v)] Let $u_1,...,u_N$ be the poles of $F$ in a unit cell, with orders $\mu_1,...,\mu_N$ respectively; let $v_1,...,v_{N'}$ be the roots of $F$ in a unit cell, with orders $\nu_1,...,\nu_{N'}$ respectively. Then the weighted sum of locations of poles is congruent to the weighted sum of locations of roots, \ie there exist integers $n,n'$ such that
\be
\sum_{i=1}^N
\mu_iu_i
=
\sum_{j=1}^{N'}\nu_jv_j
+2n\omega_1+2n'\omega_2.
\label{b39}
\ee
(The proof follows from eq.\ (\ref{b38}) with $H(z)=z$ and $G(z)=F(z)$.)
\end{enumerate}

\paragraph{Weierstrass elliptic functions.} Let $\omega_1,\omega_2\in\CC^*$ be such that $\text{Im}(\omega_2/\omega_1)\neq0$. Then the {\it Weierstrass elliptic function} with half-periods $\omega_1,\omega_2$ is
\be
\wp(z,\omega_1,\omega_2)
\equiv
\frac{1}{z^2}
+\sum_{(m,n)\in\ZZ^2\backslash\{(0,0)\}}
\biggl[
\frac{1}{(z+2m\omega_1+2n\omega_2)^2}
-\frac{1}{(2m\omega_1+2n\omega_2)^2}
\biggr].
\label{bb40}
\ee
It is an even function whose Laurent series at the origin reads
\be
\wp(z,\omega_1,\omega_2)
=
\frac{1}{z^2}
+\frac{g_2(\omega_1,\omega_2)}{20}z^2
+\frac{g_3(\omega_1,\omega_2)}{28}z^4+\cO(z^6)
\label{b40}
\ee
in terms of the {\it Weierstrass invariants}\footnote{Up to normalization, these invariants are Eisenstein series of weights $4$ and $6$ \cite{Apostol}.}
\be
\begin{split}
g_2(\omega_1,\omega_2)
&\equiv
60\sum_{(m,n)\in\ZZ^2\backslash\{(0,0)\}}
\frac{1}{(2m\omega_1+2n\omega_2)^4},\\
g_3(\omega_1,\omega_2)
&\equiv
140\sum_{(m,n)\in\ZZ^2\backslash\{(0,0)\}}
\frac{1}{(2m\omega_1+2n\omega_2)^6}.
\end{split}
\label{bbb40}
\ee
Note that there is no constant ($z$-independent) term in the expansion (\ref{b40}). We will use this at the end of section \ref{secB2} to relate $\wp$ to $\sn^2$ (see eq.\ (\ref{b56})).\\

As $\wp$ is elliptic with half-periods $\omega_1,\omega_2$, so is any integer power of $\wp$ or of its derivatives. This is true, in particular, of the combination $\wp'^2-4\wp^3+g_2\wp+g_3$. Using the series (\ref{b40}), one verifies that this combination vanishes at $z=0$. Since it has the same periods as the $\wp$ function, its regularity at $z=0$ implies that it actually has no poles at all, so by property (iii) it vanishes identically. It follows that, in fact,
\be
\wp'^2
=
4\wp^3-g_2\wp-g_3.
\label{b41}
\ee
This differential equation is a defining property of the $\wp$ function: its general solution is $F(z)=\wp(z-z_0)$ for some constant $z_0$. Note that eq.\ (\ref{b41}) is an integrated form of the KdV equation (\ref{s5}) for travelling waves, with $g_1,g_2$ as constants of integration.\footnote{If $p(x,\tau)=p(x-v\tau)$ solves (\ref{s5}), then $F(x)\equiv p(x)-v/3$ solves $3FF'-\frac{c}{12}F'''=0$, which can be integrated into $\frac{3}{2}F^2-\frac{c}{12}F''=g_2/8$, then into $\frac{c}{3}F'^2=4F^3-g_2F-g_3$. This is eq.\ (\ref{b41}) upon rescaling $x$.} Thus, any KdV soliton is related to the Weierstrass $\wp$ function \cite{NovikovBook}. Eq.\ (\ref{s57}) for cnoidal waves was a particular instance of that fact; we shall prove that equation below (see eq.\ (\ref{b56})).\\

The values of the $\wp$ function at its half-periods played an important role in section \ref{sec32}, as they coincide with the boundaries $V\in\{{-}\tfrac{m+1}{3},\tfrac{2m-1}{3},\tfrac{2-m}{3}\}$ of the four regions shown in fig.\ \ref{s71b} (more on that in section \ref{secB2}). These values are conventionally written as
\be
\wp(\omega_1)
\equiv e_1,
\quad
\wp(\omega_2)
\equiv e_2,
\quad
\wp(\omega_3)
\equiv e_3,
\qquad
\text{with}
\quad\omega_3\equiv-\omega_1-\omega_2.
\label{b43}
\ee
The periodicities and the evenness of the $\wp$ function imply that $\wp'(\omega_1)=\wp'(\omega_2)=\wp'(\omega_3)=0$, which in turn implies that the numbers (\ref{b43}) are the three roots of the polynomial on the right-hand side of eq.\ (\ref{b41}): $4e_i^3-g_2e_i-g_3=0$. Thus, writing the polynomial as $4t^3-g_2t-g_3=4(t-e_1)(t-e_2)(t-e_3)$, one reads off the relations
\be
e_1+e_2+e_3=0,
\qquad
e_1e_2+e_2e_3+e_3e_1
=
-\frac{g_2}{4},
\qquad
e_1e_2e_3
=
\frac{g_3}{4}.
\label{bb43}
\ee
A special case, relevant for Jacobi elliptic functions, occurs when $\omega_1$ is purely real and $\omega_2$ is purely imaginary (for definiteness, take $\omega_1>0$ and $\text{Im}(\omega_2)>0$). Then the definition (\ref{bb40}) ensures that $\wp(z)$ is real if and only if $z$ belongs to the grid
\be
\bigl(\RR+\ZZ\omega_2\bigr)
\cup
\bigl(\ZZ\omega_1+i\RR\bigr),
\label{b43b}
\ee
of which eq.\ (\ref{grid}) is a special case. In that situation, the Weierstrass invariants (\ref{bbb40}) and the roots (\ref{b43}) are all real.

\paragraph{Zeta functions.} The {\it Weierstrass zeta function}\footnote{Not to be confused with the unrelated Riemann zeta function.} with half-periods $\omega_1,\omega_2$ is the unique function $\zeta(z,\omega_1,\omega_2)$ such that (keeping periods implicit)
\be
\zeta'(z)
=
-\wp(z)
\qquad\text{and}\qquad
\lim_{z\to0}
\Bigl(
\zeta(z)-\frac{1}{z}
\Bigr)
=
0.
\label{b44}
\ee
Since $\wp$ is even, $\zeta$ is odd. Owing to the defining series (\ref{bb40}) of the $\wp$ function, the series representation of the zeta function is
\be
\zeta(z,\omega_1,\omega_2)
=
\frac{1}{z}
+
\sum_{(m,n)\in\ZZ^2\backslash\{(0,0)\}}
\biggl(
\frac{1}{z+2m\omega_1+2n\omega_2}
-\frac{1}{2m\omega_1+2n\omega_2}
+\frac{z}{(2m\omega_1+2n\omega_2)^2}
\biggr).
\label{bb44}
\ee
Importantly, the zeta function is {\it not} periodic: it is only quasi-periodic,
\be
\zeta(z+2\omega_i)
=
\zeta(z)+2\zeta(\omega_i)
\qquad
\forall\,i=1,2,
\label{bbb44}
\ee
as follows from the oddness of $\zeta$ and the integral of the periodicity relation $\wp(z+2\omega_i)=\wp(z)$. In particular, $\zeta(z)$ is {\it not} an elliptic function. For instance, the first equation in (\ref{t156}) is a special case of (\ref{bbb44}).\\

The values of $\zeta(z)$ at half-periods satisfy an important identity that we used in eqs.\ (\ref{ss68}) and (\ref{ipitu}). Namely, note that in any unit cell of the $\wp$ function, $\zeta(z)$ has only one pole of order one. The residue theorem then yields
\be
\oint_{\der\,\text{cell}}
\dd z\,\zeta(z)
=
2\pi i,
\label{b45}
\ee
where the cell's contour is traced in a counterclockwise direction. Now, assuming that $\text{Im}(\omega_2/\omega_1)>0$ and using the quasi-periodicity (\ref{bbb44}), the left-hand side of (\ref{b45}) reduces to $-4\omega_1\zeta(\omega_2)+4\omega_2\zeta(\omega_1)$. It follows that
\be
\omega_1\zeta(\omega_2)
-\omega_2\zeta(\omega_1)
=
-i\pi/2,
\label{bb45}
\ee
of which eq.\ (\ref{ss68}) is a special case.\footnote{One sometimes writes $\zeta(\omega_i)\equiv\eta_i$, whereupon (\ref{bb45}) reads $\omega_1\eta_2-\omega_2\eta_1=-i\pi/2$.} With the opposite orientation for the periods ($\text{Im}(\omega_2/\omega_1)<0$), the right-hand side of (\ref{bb45}) would be $+i\pi/2$.

\paragraph{Sigma functions.} The {\it Weierstrass sigma function} with half-periods $\omega_1,\omega_2$ is the unique function $\sigma(z,\omega_1,\omega_2)$ such that (keeping periods implicit)
\be
\frac{\sigma'(z)}{\sigma(z)}
=
\zeta(z)
\qquad\text{and}\qquad
\sigma(z)\stackrel{z\to0}{\sim}
z.
\label{bbb47}
\ee
Using the series representation (\ref{bb44}) of the zeta function, it follows that
\be
\sigma(z,\omega_1,\omega_2)
=
z
\prod_{(m,n)\in\ZZ^2\backslash\{(0,0)\}}
\biggl(
\frac{z+2m\omega_1+2n\omega_2}{2m\omega_1+2n\omega_2}
e^{-\frac{z}{2m\omega_1+2n\omega_2}+\frac{z^2}{2(2m\omega_1+2n\omega_2)^2}}
\biggr).
\nn
\ee
It is thus an analytic, odd function, with first order roots located at the poles of $\wp(z)$. It is also quasi-periodic in the sense that
\be
\sigma(z+2\omega_i)
=
-e^{2\zeta(\omega_i)(z+\omega_i)}
\sigma(z)
\qquad\forall\,i=1,2,
\label{bb47}
\ee
which follows from the oddness of $\sigma$ and the integral of the relation (\ref{bbb44}). In particular, it is {\it not} an elliptic function. The second equation in (\ref{t156}) is a special case of eq.\ (\ref{bb47}).

\paragraph{Universality of Weierstrass functions.} We now derive an important fact that will lead to addition formulas, and that we will eventually use in section \ref{secB2}. Namely, the Weierstrass $\wp$ function is universal in the sense that any elliptic function with half-periods $\omega_1,\omega_2$ can be expressed in terms of $\wp(z,\omega_1,\omega_2)$ and $\wp'(z,\omega_1,\omega_2)$. Indeed, let $F(z)$ be an elliptic function and write it as
\be
F(z)
=
\demi\bigl(F(z)+F(-z)\bigr)
+\demi\frac{F(z)-F(-z)}{\wp'(z)}\wp'(z),
\label{b47}
\ee
so we may assume that $F$ is even without loss of generality. Now let $u_1,...,u_N$ be the poles of $F$ in a unit cell, with orders $\mu_1,...,\mu_N$, and let $v_1,...,v_{N'}$ be its roots in a unit cell, with orders $\nu_1,...,\nu_{N'}$. Then the function
\be
F(z)
\cdot
\frac{\prod_{i=1}^N(\wp(z)-\wp(u_i))^{\mu_i}}{\prod_{j=1}^{N'}(\wp(z)-\wp(v_j))^{\nu_j}}
\nn
\ee
is elliptic but has no poles, so in fact it is a constant. It follows that
\be
F(z)
=
C
\frac{\prod_{j=1}^{N'}(\wp(z)-\wp(v_j))^{\nu_j}}{\prod_{i=1}^N(\wp(z)-\wp(u_i))^{\mu_i}}
\label{b48}
\ee
for some constant $C$. Thus, any elliptic function (\ref{b47}) is the sum of a rational function of $\wp(z)$ and the product of $\wp'(z)$ with another rational function of $\wp(z)$. In particular, there exists an algebraic relation between any elliptic function and (i) its derivative, (ii) the same function with shifted arguments. This implies that elliptic functions satisfy {\it algebraic addition theorems}: $F(z+w)$ can be written as a rational function of $F(z)$ and $F(w)$. We shall return to this shortly.\\

In practice, a helpful way to write any elliptic function is in terms of the sigma function (though the latter is not elliptic). To see this, let $F(z)$ be an elliptic function with poles $u_1,...,u_N$ (with orders $\mu_1,...,\mu_N$) and roots $v_1,...,v_{N'}$ (with orders $\nu_1,...,\nu_{N'}$) in a unit cell. Eq.\ (\ref{b39}) allows us to choose this cell in such a way that $\sum_i\mu_i u_i=\sum_j\nu_j v_j$; we assume such a choice has been made. Then consider the function
\be
F(z)
\cdot
\frac{\prod_{i=1}^N(\sigma(z-u_i))^{\mu_i}}{\prod_{j=1}^{N'}(\sigma(z-v_j))^{\nu_j}}.
\label{b49}
\ee
This function has neither roots nor poles in a cell. Furthermore, using the quasi-periodicity (\ref{bb47}) of $\sigma$ and the choice $\sum\mu u=\sum\nu v$, one sees that (\ref{b49}) is an elliptic function, so in fact it is a constant. It follows that
\be
F(z)
=
C
\frac{\prod_{j=1}^{N'}(\sigma(z-v_j))^{\nu_j}}{\prod_{i=1}^N(\sigma(z-u_i))^{\mu_i}}
\label{bb49}
\ee
for some constant $C$. In particular, any elliptic function is completely determined, up to a multiplicative constant, by its poles and roots (and their degrees). We will use this at the end of section \ref{secB2} to express the square of a Jacobi sine as a Weierstrass $\wp$ function.

\paragraph{Addition formulas for Weierstrass functions.} We saw above that elliptic functions satisfy algebraic addition theorems. We now work this out for the $\wp$ function, as it is a key tool in proving that (\ref{s156}) solves the Lam\'e equation (\ref{ss154}). Let $z,w\in\CC\backslash\Gamma$ (where $\Gamma$ is the lattice (\ref{bb37})) and define 
\be
c(z,w)
\equiv
\frac{\wp'(z)-\wp'(w)}{\wp(z)-\wp(w)},
\qquad
d(z,w)
\equiv
\frac{\wp(z)\wp'(w)-\wp(w)\wp'(z)}{\wp(z)-\wp(w)}.
\label{b50}
\ee
Thinking of $z,w$ as fixed parameters, we define an elliptic function $F(x)\equiv\wp'(x)-c(z,w)\wp(x)-d(z,w)$. By construction, in any unit cell, the latter has only one pole of order three (at $x=0$ and congruent points). Furthermore, it has two simple roots: one at $x=z$, the other at $x=w$. Property (iv) above then implies the existence of a third root, while property (v) implies that this third root is congruent to $x=-z-w$. Accordingly, $F(-z-w)=0$. Now consider the function
\be
G(x)
\equiv
\wp'(x)^2
-
(c\wp(x)+d)^2
=
F(x)(\wp'(x)+c\wp(x)+d),
\nn
\ee
which can also be written as follows thanks to the Weierstrass differential equation (\ref{b41}):
\be
G
=
4\wp^3
-
c^2
\wp^2
-
(2cd+g_2)\wp
-(d^2+g_3).
\label{b51}
\ee
This function vanishes whenever $F$ vanishes, hence $G(z)=G(w)=G(-z-w)=0$. Thus the polynomial on the right-hand side of (\ref{b51}) can be written as
\be
4t^3-c^2t^2-(2cd+g_2)t-(d^2+g_3)
=
4(t-\wp(z))(t-\wp(w))(t-\wp(-z-w)).
\nn
\ee
Expanding the right-hand side of this expression and comparing with the left-hand side, then using the definition (\ref{b50}) of $c(z,w)$, one finds the first equation in (\ref{bb51}): the expected addition theorem for the $\wp$ function. (The right-hand side of that formula contains $\wp'$, but the latter can be expressed in terms of $\wp$ thanks to the Weierstrass equation (\ref{b41}), which yields an algebraic relation between $\wp(z+w)$ and $\wp(z),\wp(w)$.)\\

An addition theorem can similarly be derived for the sigma function, as follows (though $\sigma(z)$ is not elliptic). Fix $w\in\CC\backslash\Gamma$ and think of $\wp(z)-\wp(w)$ as a function of $z$. In each cell, it has a second-order pole congruent to $z=0$, and simple roots congruent to $z=\pm w$. According to eq.\ (\ref{bb49}), this implies there exists a constant $C$ such that $\wp(z)-\wp(w)=C(\sigma(z-w)\sigma(z+w))/\sigma(z)^2$. The value of $C$ is found by expanding this equation around $z=0$, which results in
\be
\wp(z)-\wp(w)
=
-
\frac{\sigma(z-w)\sigma(z+w)}{\sigma(z)^2\sigma(w)^2}.
\nn
\ee
This is a pseudo-addition theorem for the $\sigma$ function. Differentiating with respect to $z$ and $w$ and using eq.\ (\ref{bbb47}), one also obtains the pseudo-addition theorem for $\zeta$ announced in eq.\ (\ref{bb51}). The latter identity implies in particular that $\zeta(\omega_1+\omega_2)=\zeta(\omega_1)+\zeta(\omega_2)$, since $\wp'(\omega_i)=0$. Along with eq.\ (\ref{bb45}), this implies that 
\be
\omega_1\zeta(\omega_1+\omega_2)
-\zeta(\omega_1)(\omega_1+\omega_2)
=
-i\pi/2,
\label{APAKI}
\ee
which we used in eq.\ (\ref{ipitu}).

\subsection{Jacobi elliptic functions}
\label{secB2}

To define Jacobi elliptic functions, we shall first think of them as generalizations of trigonometric functions. Then we shall review their addition identities, which will allow us to extend them to the complex plane and show that they are, in fact, elliptic. We will then relate them to the $\wp$ function. For more details we refer to \cite[chap.\ {\sc xxii}]{Whittaker} or \cite{Singh}.

\paragraph{Generalizing trigonometry.} Consider a plane $\RR^2$ with Cartesian coordinates $(x,y)$. Let
\be
x^2+(1-m)y^2
=
1
\label{ss146}
\ee
be the equation of an ellipse with eccentricity (or `modulus') $\sqrt{m}$.\footnote{One also writes $\sqrt{m}=k$, but we stick to $m$ to avoid confusion with the parameter $k$ of eq.\ (\ref{ss159}).} To introduce Jacobi elliptic functions we need to relate two coordinates on the ellipse. First, any point can be labelled by the angle $\theta$ such that $x=r\cos\theta$ and $y=r\sin\theta$; in these terms, eq.\ (\ref{ss146}) specifying the ellipse reads
\be
r
=
\frac{1}{\sqrt{1-m\sin^2\theta}}.
\label{s146}
\ee
Alternatively one may label a point on the ellipse by a coordinate $u$ such that $\dd u=r(\theta)\dd\theta$, with $u=0$ at the 'east pole' $\theta=0$. It is the purely angular part of the arc length, \ie the {\it incomplete elliptic integral of the first kind} (sometimes denoted $F(\theta|m)$):
\be
u(\theta)
=
\int_0^{\theta}
\frac{\dd\theta'}{\sqrt{1-m\sin^2\theta'}}.
\label{JAMMIN}
\ee
Then the {\it Jacobi sine and cosine functions} with parameter $m$ are defined as the unique functions $\sn(u|m)$ and $\cn(u|m)$ such that, for any $\theta$,
\be
\sn\bigl(u(\theta)\big|m\bigr)
\equiv\sin\theta,
\qquad
\cn\bigl(u(\theta)\big|m\bigr)
\equiv\cos\theta.
\label{t147}
\ee
One also defines the {\it delta amplitude}
\be
\dn(u|m)
\equiv
\sqrt{1-m\sn(u|m)^2}.
\label{s148}
\ee
The definition (\ref{t147}) readily implies that $\cn(u|m)^2+\sn(u|m)^2=1$ (for real $u$), and also that Jacobi elliptic functions are periodic: their period in $u$ is conventionally written as $4K(m)$ in terms of the {\it complete elliptic integral of the first kind}
\be
K(m)
\equiv
\int_0^{\pi/2}
\frac{\dd\theta}{\sqrt{1-m\sin^2\theta}}.
\label{ss147}
\ee
Similarly, the delta amplitude (\ref{s148}) has period $2K(m)$. Note that $K(m)$ is monotonously increasing as a function of $m$, and $K(m)\sim-\log\sqrt{1-m}$ as $m\to1$. In section \ref{sec4} we also used the {\it complete elliptic integral of the second kind}
\be
E(m)
\equiv
\int_0^{\pi/2}\dd\theta\sqrt{1-m\sin^2\theta}
=
\int_0^{K(m)}\dd u\dn(u|m)^2,
\label{NADA}
\ee
which is monotonously decreasing as a function of $m$ and satisfies $E(1)=1$. When $m=0$, the ellipse becomes a circle and Jacobi elliptic functions reduce to standard trigonometric functions $\sn(u|0)=\sin u$, $\cn(u|0)=\cos u$, while $\dn(u|0)=1$. The complete elliptic integrals then reduce to $K(0)=E(0)=\pi/2$.

\paragraph{Derivatives.} To see that cnoidal waves (\ref{ss6}) solve the KdV equation (\ref{s5}), one needs to differentiate Jacobi elliptic functions. For example, using the chain rule, one finds
\begin{align}
\cn(u(\theta)|m)
&=
\cos\theta
=
\frac{\dd}{\dd\theta}\sin\theta
=
\frac{\dd}{\dd\theta}\sn(u(\theta)|m)
=
\frac{\dd u}{\dd\theta}\frac{\dd}{\dd u}\sn(u|m)\big|_{u(\theta)}\nn\\
&\stackrel{\text{(\ref{JAMMIN})}}{=}
\frac{1}{\sqrt{1-m\sin^2\theta}}\frac{\dd}{\dd u}\sn(u|m)\big|_{u(\theta)}
\stackrel{\text{(\ref{s148})}}{=}
\frac{1}{\dn(u(\theta)|m)}\frac{\dd}{\dd u}\sn(u|m)\big|_{u(\theta)}.\nn
\end{align}
Similar arguments can be used to show all three differentiation identities
\begin{align}
\frac{\dd}{\dd u}\sn(u|m)
&=
\dn(u|m)\cn(u|m),\nn\\
\frac{\dd}{\dd u}\cn(u|m)
&=
-\dn(u|m)\sn(u|m),\label{sb35}\\
\frac{\dd}{\dd u}\dn(u|m)
&=
-m\sn(u|m)\cn(u|m),\nn
\end{align}
which generalize derivatives of trigonometric functions. In particular, it readily follows that (\ref{ss6}) is indeed a solution of the KdV equation (\ref{s5}).\\

Note that the derivatives in (\ref{sb35}) can be used to prove eq.\ (\ref{MYSTERY}). Indeed, keeping the $m$ dependence implicit, one has \cite[pp.\ 87-88]{Lawden}
\be
\frac{\dd^2}{\dd u^2}\Big(\dn(u)^N\Big)
=
-N(N+1)\dn(u)^{N+2}+(2-m)N^2\dn(u)^N-N(N-1)(1-m)\dn(u)^{N-2}.
\nn
\ee
Thus, if we let $I_N\equiv\int_0^{2K(m)}\dd u\dn(u)^N$, we find the recursion relation
\be
(N+1)I_{N+2}
=
(2-m)NI_N-(N-1)(1-m)I_{N-2}.
\label{TADAM}
\ee
For $N=2$, using $I_2=2E(m)$ by eq.\ (\ref{NADA}), it follows that $3I_4=(2m-2)K(m)+(8-4m)E(m)$. This proves eq.\ (\ref{MYSTERY}).\\

The remainder of this appendix is devoted to the proof of the relation between squares of Jacobi functions and the Weierstrass $\wp$ function. This relation was stated in eq.\ (\ref{s57}) and was essential to relate Hill's equation for cnoidal waves to the Lam\'e equation. Aside from this application, however, most of the rest of this appendix may safely be skipped.

\paragraph{Addition formulas.} The derivatives (\ref{sb35}) can be used to obtain addition formulas for Jacobi elliptic functions. Indeed, consider two variables $u,v$ such that $u+v\equiv\alpha$  is a constant, and write $\sn\,u\equiv s_1$, $\sn\,v\equiv s_2$, $\cn\,u\equiv c_1$, $\cn\,v\equiv c_2$, $\dn\,u\equiv d_1$ and $\dn\,v\equiv d_2$. (From now on we omit to stress the parametric $m$ dependence.) Then, since $u+v$ is constant, the derivatives (\ref{sb35}) along with elementary identities yield
\be
\frac{\dd}{\dd u}(s_1c_2+s_2c_1)
=
(d_1-d_2)(c_1c_2-s_1s_2),
\qquad
\frac{\dd}{\dd u}(d_1+d_2)
=
-m(c_1c_2-s_1s_2)(s_1c_2-s_2c_1).
\nn
\ee
Noting that $-m(s_1^2c_2^2-s_2^2c_1^2)=d_1^2-d_2^2$, this gives $\der_u\log[s_1c_2+s_2c_1]=\der_u\log[d_1+d_2]$. It follows that the ratio $(s_1c_2+s_2c_1)/(d_1+d_2)$ is constant when $u+v$ is constant; the value of the ratio is found by setting $u=0,v=\alpha$. Along with a similar line of thought for $s_1c_2-s_2c_1$ and $d_1-d_2$, one finds
\be
\frac{d_1+d_2}{s_1c_2+s_2c_1}
=
\frac{\dn(u+v)+1}{\sn(u+v)},
\qquad
\frac{d_1-d_2}{s_1c_2-s_2c_1}
=
\frac{\dn(u+v)-1}{\sn(u+v)}.
\nn
\ee
The difference between these equations gives an addition formula for $\sn(u)$; their sum yields an addition formula for $\dn(u)$. Including also the addition formula for $\cn(u)$, which can be derived in a similar manner, one finally obtains
\be
\sn(u+v)
=
\frac{s_1^2-s_2^2}{s_1c_2d_2-s_2c_1d_1},
\quad
\cn(u+v)
=
\frac{s_1c_1d_2-s_2c_2d_1}{s_1c_2d_2-s_2c_1d_1},
\quad
\dn(u+v)
=
\frac{s_1c_2d_1-s_2c_1d_2}{s_1c_2d_2-s_2c_1d_1}.
\label{b32b}
\ee
It is not obvious that these equations reduce to the standard addition formulas for trigonometric functions when $m=0$. To see that they do, one can actually rewrite (\ref{b32b}) in a more suggestive form that follows \eg from the extension of Jacobi elliptic functions to the complex plane. We omit this argument here for brevity and refer instead to \cite{Singh}.

\paragraph{Jacobi elliptic functions in the complex plane.} The addition formulas (\ref{b32b}) allow us to extend Jacobi elliptic functions to the complex plane by analytic continuation. Namely, note first that by defining a variable $t\in[0,1]$ such that $\sin\theta=2t/(1+t^2)$ for $\theta\in[0,\pi/2]$, we can rewrite the incomplete elliptic integral (\ref{JAMMIN}) as
\be
u(t)
=
2\int_0^t\frac{\dd s}{\sqrt{1+2(1-2m)s^2+s^4}}.
\label{sb37}
\ee
With this parametrization, the Jacobi sine and cosine read
\be
\sn(u(t)|m)
=
\frac{2t}{1+t^2},
\qquad
\cn(u(t)|m)
=
\frac{1-t^2}{1+t^2}.
\label{tb38}
\ee
From this perspective it is now straightforward to define Jacobi elliptic functions with a purely imaginary argument: multiplying (\ref{sb37}) by $i$ we find
\be
iu(t)\big|_m
=
\int_0^t\frac{i\dd s}{\sqrt{1+2(1-2m)s^2+s^4}}
=
\int_0^{it}\frac{\dd z}{\sqrt{1+2(1-2(1-m))z^2+z^4}}
=
u(it)\big|_{1-m}.
\nn
\ee
Using this, we can define for instance the Jacobi elliptic sine on the imaginary axis:
\be
\sn(iu(t)|m)
=
\sn(u(it)|1-m)
\stackrel{\text{(\ref{tb38})}}{=}
=
\frac{2it}{1-t^2}
=
\frac{i\sn(u(t)|1-m)}{\cn(u(t)|1-m)}.
\nn
\ee
More generally, one has {\it Jacobi's imaginary transformations}
\be
\sn(iu|m)
=
\frac{i\sn(u|1-m)}{\cn(u|1-m)},
\quad
\cn(iu|m)
=
\frac{1}{\cn(u|1-m)},
\quad
\dn(iu|m)
=
\frac{\dn(u|1-m)}{\cn(u|1-m)}.
\label{ttb38}
\ee
Given the known elliptic functions on the real line, these transformations define the continuation of these functions to the imaginary axis, minus the points where $\cn(u|1-m)$ vanishes, \ie points of the form $(2n+1)iK(1{-}m)$ with $n\in\ZZ$. Note that these definitions readily imply that $\cn(z|m)$ and $\dn(z|m)$ have period $4iK(1{-}m)$ as functions of $z\in i\RR$, while $\sn(z|m)$ has period $2iK(1{-}m)$.\\

Having defined Jacobi elliptic functions both on the real line and on the imaginary axis, one can extend them on the complex plane using addition formulas. For instance, one defines $\sn(z)=\sn(x+iy)$ by eq.\ (\ref{b32b}) with $u=x$ and $v=iy$, and with Jacobi functions of purely imaginary arguments defined by (\ref{ttb38}). The resulting functions $\sn(z|m)$, $\cn(z|m)$ and $\dn(z|m)$ turn out to be meromorphic (they satisfy the Cauchy-Riemann equations), with respective real periods $4K(m)$, $4K(m)$, $2K(m)$ and imaginary periods $2iK(1{-}m)$, $4iK(1{-}m)$, $4iK(1{-}m)$. Accordingly, they are by definition elliptic functions in the sense of section \ref{appBsec1}. Due to the $\cn$ in the denominator of (\ref{ttb38}), these functions diverge at all points of the form $2nK(m)+i(2n'+1)K(1{-}m)$ with integers $n,n'$. Finally, the derivative identities (\ref{sb35}) still hold when $u$ is a complex variable, as do the addition formulas (\ref{b32b}). See again \cite[chap.\ {\sc xxii}]{Whittaker} or \cite{Singh} for details.

\paragraph{Relation to Weierstrass functions.} Since Jacobi elliptic functions are elliptic, and since any elliptic function can be written in terms of the Weierstrass $\wp$ function, it is possible to express the complex functions $\sn(z)$, $\cn(z)$ and $\dn(z)$ in terms of $\wp(z)$ with suitable periods. We will not review this rewriting here, as it is not essential for our purposes. However, the one relation that does play a key role for us is eq.\ (\ref{s57}), which relates $\sn^2$ and $\wp$. Accordingly, we now derive this equation.\\

To start, note that $\sn^2(z|m)$ is an even elliptic function with periods $2K(m)$ and $2iK(1{-}m)$. Owing to eq.\ (\ref{b48}), it is thus possible to write $\sn^2(z|m)$ as a rational function of $\wp(z,K(m),iK(1{-}m))$. Now recall our conclusion below eq.\ (\ref{bb49}) that any elliptic function is fully determined, up to normalization, by its poles and roots. In the case at hand, $\sn^2(z|m)$ has exactly one pole of order 2, and exactly one root of order 2, in each unit cell.\footnote{The root's locations are $2nK(m)+2in'K(1{-}m)$, while the pole's locations are $2nK(m)+i(2n'+1)K(1{-}m)$. The multiplicity follows from the fact that $\sn(z|m)$ only has first-order roots.} Since the $\wp$ function is the unique elliptic function having exactly one pole of order 2 in each cell, there must exist constants $A,B\in\CC$ such that
\be
\sn^2(z|m)
=
A\wp\bigl(z+iK(1{-}m),K(m),iK(1{-}m)\bigr)+B
\label{b56}
\ee
where the shift in the argument of $\wp$ ensures that the locations of poles and roots match on the two sides of the equation.\\

Our only remaining task is to find $A$ and $B$ in eq.\ (\ref{b56}). To do this, we expand the left-hand side around $z=iK(1{-}m)$; for simplicity we also restrict attention to purely imaginary arguments, where we can use eqs.\ (\ref{ttb38}) directly (the extension to complex arguments then follows form meromorphicity). Thus we find
\be
\sn^2(iy+iK(1{-}m)|m)
\stackrel{y\to0}{\sim}
-\frac{\Bigl[\sn\bigl(K(1{-}m)\big|1{-}m\bigr)+\frac{y^2}{2}\sn''\bigl(K(1{-}m)\big|1{-}m\bigr)\Bigr]^2}{\Bigl[y\cn'\bigl(K(1{-}m)\big|1{-}m\bigr)+\frac{y^3}{6}\cn'''\bigl(K(1{-}m)\big|1{-}m\bigr)\Bigr]^2}
\nn
\ee
where we used $\sn'(K(1{-}m)|1{-}m)=\cn(K(1{-}m)|1{-}m)=\cn''(K(1{-}m)|1{-}m)=0$. Then using the derivatives (\ref{sb35}) along with $\sn'(K(1{-}m)|1{-}m))=0$, we get
\be
m\,\sn^2\bigl(z+iK(1{-}m)\big|m\bigr)
\stackrel{z\to0}{\sim}
\frac{1}{z^2}+\frac{m+1}{3}.
\nn
\ee
Comparing this with the expansion (\ref{bb40}) of the $\wp$ function, we deduce that $A=1/m$ and $B=(m+1)/(3m)$ in eq.\ (\ref{b56}), which yields the expected result (\ref{s57}).

\paragraph{Special values.} The values of $V$ at the bifurcations $V\in\{{-}\tfrac{m+1}{3},\tfrac{2m-1}{3},\tfrac{2-m}{3}\}$ follow from the relation (\ref{s57}) between the $\wp$ function and $\sn^2$. Indeed, we have
\be
\begin{split}
e_2
&\equiv
\wp(\omega_2)=\wp(iK(1{-}m))=m\,\sn^2(0)-\frac{m+1}{3}=-\frac{m+1}{3},\\
e_3
&\equiv
\wp(\omega_3)=\wp(K(m)+iK(1{-}m))=m\,\sn^2(K(m))-\frac{m+1}{3}=\frac{2m-1}{3},
\end{split}
\label{bb57}
\ee
and the third value can be deduced from eq.\ (\ref{bb43}), namely
\be
e_1
=
-e_2-e_3
=
\frac{m+1}{3}+\frac{1-2m}{3}
=
\frac{2-m}{3}.
\label{b58}
\ee
Since $e_2$, $e_3$ and $e_1$ are the values of $V$ at the corners of the rectangle depicted in fig.\ \ref{s67}, this explains why $V\in\{{-}\tfrac{m+1}{3},\tfrac{2m-1}{3},\tfrac{2-m}{3}\}$ are the boundaries of regions with sharply different orbits in fig.\ \ref{s71b}.\\

It also follows from (\ref{bb57})-(\ref{b58}) and eqs.\ (\ref{bb43}) that, for $\omega_1=K(m)$ and $\omega_2=iK(1{-}m)$, the Weierstrass invariants are $g_2=\tfrac{4}{3}(m^2-m+1)$, $g_3=\tfrac{4}{27}(2m^3-3m^2-3m+2)$. Accordingly, using the Weierstrass equation (\ref{b41}), we find the second derivatives
\be
\wp''(iK(1{-}m))=2m,
\qquad
\wp''(K(m)+iK(1{-}m))=2m(m-1),
\qquad
\wp''(K(m))=2(1-m).
\label{KAPPOT}
\ee
We used the latter of these relations in (\ref{finidev}), while the first one appeared above (\ref{LAOS}).

\section{Asymptotics of elliptic functions}
\label{appC}

Here we derive the asymptotic behaviour of level curves of $k$ as displayed in section \ref{sec36} (see in particular eqs.\ (\ref{NOPOTO}) and (\ref{POTO})). Specifically, we analyse Weierstrass functions in the region where one of their periods goes to infinity. In terms of eq.\ (\ref{ss159}) and fig.\ \ref{s71b}, this corresponds to the regions $m\to1$ and $m\to0$, where either $K(m)$ or $K(1{-}m)$ blows up.\\

We start by listing some standard asymptotic identities. First, when one of the periods of the $\wp$ function goes to $\infty$ while the other one is finite (say $|\omega_2/\omega_1|\to\infty$ with $\text{Im}(\omega_2/\omega_1)>0$), one has \cite[eq.\ 23.12.1]{DLMF}
\be
\wp(z,\omega_1,\omega_2)
\sim
\frac{\pi^2}{4\omega_1^2}
\biggl[
-\frac{1}{3}
+\frac{1}{\sin^2(\pi z/(2\omega_1))}
+8\bigl(1-\cos(\pi z/\omega_1)\bigr)e^{2\pi i\omega_2/\omega_1}
+\cO(e^{4\pi i \omega_2/\omega_1})
\biggr]
\label{KLA}
\ee
We used this in section \ref{sec444} to simplify (\ref{VAPAT}) into (\ref{VANT}). As regards the Weierstrass zeta function, we use the defining equation $\zeta'=-\wp$ to integrate (\ref{KLA}) into
\be
\zeta(z,\omega_1,\omega_2)
\stackrel{\omega_2\to\infty}{\sim}
\frac{\pi^2}{4\omega_1^2}
\biggl[
\frac{z}{3}
+\frac{2\omega_1}{\pi}\cot\Bigl(\frac{\pi z}{2\omega_1}\Bigr)
-8\Bigl(z-\frac{\omega_1}{\pi}\sin(\pi z/\omega_1)\Bigr)e^{2\pi i\omega_2/\omega_1}
+\cO(e^{4\pi i \omega_2/\omega_1})
\biggr],
\label{KLAZ}
\ee
where the integration constant is set to zero by requiring $\lim_{z\to0}(\zeta(z)-1/z)=0$. Eq.\ (\ref{PAKATTO}) is a special case of this result, as is eq.\ (\ref{NEWT}) when $\omega_1$ goes to $+\infty$ while $\omega_2$ remains finite. In particular, we deduce from (\ref{NEWT}) that, when $\omega_1=K(m)$ and $\omega_2=iK(1{-}m)$,
\be
\zeta(iK(1{-}m))
\;\stackrel{m\to1}{\sim}\;
-\frac{i\pi^2}{12K(1{-}m)}+\cO\bigl((1{-}m)^2\bigr).
\nn
\ee
Eq.\ (\ref{n20t}) then follows upon using identity (\ref{bb45}). This justifies all the arguments of section \ref{sec36}, except for eqs.\ (\ref{NOPOTO}) and (\ref{POTO}), to which we now turn.

\paragraph{Proof of (\ref{NOPOTO}).} We consider the limit $m\to1$ of a level curve for which $k/c<-1/24$. Anticipating the non-perturbative behaviour in eq.\ (\ref{NOPOTO}), we shall think of the curve as a function $m=f_k(V)$ and ask how it behaves as $V\to-\tfrac{m+1}{3}$. With this in mind, we let $V\equiv-\tfrac{m+1}{3}-\nu$ and expand
\begin{align}
\wp^{-1}(V)
&\sim
iK(1{-}m)-i\sqrt{\frac{\nu}{m}}+\cO(\nu),\nn\\
\zeta(\wp^{-1}(V))
&\sim
\zeta(iK(1{-}m))-i\sqrt{\frac{\nu}{m}}\frac{m+1}{3}+\cO(\nu).\nn
\nn
\end{align}
Plugging this in the level curve (\ref{ss159}), we find
\be
K(m)\frac{m+1}{3\sqrt{m}}-\frac{\zeta(K(m))}{\sqrt{m}}
\;\stackrel{\nu\to0}{\sim}\;
\frac{1}{\sqrt{\nu}}
\biggl[
\sqrt{\left|\frac{6\pi^2 k}{c}\right|}-\frac{\pi}{2}\biggr].
\label{c47b}
\ee
Here, as a function of $m$, the left-hand side is a monotonously increasing function that blows up as $m\to1$. Since the right-hand side of (\ref{c47b}) is large when $\nu\to0$, we conclude that $m$ is close to $1$ and $K(m)$ is well approximated by the expansion (\ref{KAPADI}). Using also eq.\ (\ref{n20t}), the relation (\ref{c47b}) reduces to
\be
-\frac{1}{2}\log(1-m)+\cO(1)
\;\stackrel{\nu\to0}{\sim}\;
\frac{1}{\sqrt{\nu}}
\biggl[
\sqrt{\left|\frac{6\pi^2 k}{c}\right|}-\frac{\pi}{2}\biggr].
\nn
\ee
This indicates that the $1/\sqrt{\nu}$ divergence of $V$ translates into a logarithmic divergence of $1-m$. It motivates the definition
\be
1-m
=
1-f_k(V)
\equiv
\exp\biggl[-\frac{\pi}{\sqrt{|V+2/3|}}\bigl(\sqrt{|24k/c|}-1\bigr)\biggr]
\cdot
F_k(V)
\nn
\ee
where $F_k(V)$ is some function, finite at $V=-2/3$. Plugging this ansatz back into the above expansions, we find $F_k(-2/3)=16/e^2$, proving the announced result (\ref{NOPOTO}).

\paragraph{Proof of (\ref{POTO}).} Consider now the limit $m\to1$ of a level curve for which $k/c>-1/24$. For definiteness we take $k/c>0$, but the same derivation works for $-1/24<k/c<0$. Then, as $m\to1$ at fixed $k$, $V$ becomes close to $\tfrac{2-m}{3}$ and $\wp^{-1}(V)$ therefore gets close to $K(m)$. The latter diverges at $m=1$, so $\wp^{-1}(V)$ is very large. As explained in section \ref{secNEAR}, this divergence of $\wp^{-1}(V)$ balances out that of $K(m)$, ensuring that the asymptotic relation (\ref{ASAMAN}) holds:
\be
\sqrt{\frac{6\pi^2 k}{c}}
\sim
K(m)-\wp^{-1}(V)+\cO(1{-}m).
\nn
\ee
(In writing this we expanded $K(1{-}m)\sim\tfrac{\pi}{2}+\cO(1{-}m)$ and chose $\wp^{-1}(V)$ to range from $0$ to $K(m)$ as $V$ runs from $+\infty$ to $\tfrac{2-m}{3}$.) It follows that
\be
\wp^{-1}(V)\sim
K(m)-\sqrt{\frac{6\pi^2 k}{c}}+\cO(1{-}m)
\qquad\text{as }m\to1.
\label{c45t}
\ee
As we shall see, the $\cO(1{-}m)$ correction will yield a subleading, quadratic correction to the leading linear dependence of $V$ on $1{-}m$ near $m=1$.\\

The relation (\ref{c45t}) suggests that we define a function $\nu_k(m)$ by
\be
V
\equiv
\wp\biggl(K(m)-\sqrt{\frac{6\pi^2 k}{c}}\biggr)+\nu_k(m).
\label{c46t}
\ee
Plugging this ansatz back into (\ref{c45t}) and Taylor-expanding $\wp^{-1}$ around $\wp(K(m)-\sqrt{6\pi^2k/c})$, we find that $\nu_k(m)$ must be of order $(1-m)^2$ when $m\to1$. As a result, the definition (\ref{c46t}) implies
\be
V\sim
\wp\biggl(K(m)-\sqrt{\frac{6\pi^2 k}{c}}\biggr)+\cO\bigl((1{-}m)^2\bigr)
\qquad\text{as }m\to1.
\label{cc46t}
\ee
To simplify the right-hand side of this relation, we now use the addition formula (\ref{bb51}) and the asymptotic expansion (\ref{KLA}) of the $\wp$ function, adapted to the case where the real period $\omega_1=K(m)$ goes to infinity:
\be
\wp\biggl(K(m)-\sqrt{\frac{6\pi^2 k}{c}}\biggr)+\cO\bigl((1{-}m)^2\bigr)
\;\stackrel{m\to1}{\sim}\;
\frac{1}{3}+\biggl[
\cosh^2\bigl(\sqrt{6\pi^2k/c}\bigr)-\frac{2}{3}
\biggr](1-m)
+\cO\bigl((1{-}m)^2\bigr).
\nn
\ee
Upon plugging this back into (\ref{cc46t}), the announced result (\ref{POTO}) follows.

\end{document}